\documentclass[12pt,psamsfonts]{amsart}
\usepackage[margin=1.5in]{geometry}
\usepackage{wrapfig}
\usepackage{amssymb}
\usepackage{geometry}

\usepackage{float}

\usepackage{amsmath,amsfonts,amsthm,amstext}
\usepackage[latin1]{inputenc}
\usepackage{fancybox}
\usepackage{niceframe}
\usepackage{array}
\usepackage{newlfont}
\usepackage{verbatim}
\usepackage[dvips]{psfrag}
\usepackage{color}
\usepackage{float}
\usepackage{indentfirst}
\usepackage[normalem]{ulem}
\usepackage{pst-text}
\usepackage{graphicx}
\usepackage{graphics}
\usepackage{overpic}
\usepackage{hyperref}
\usepackage{pst-node}
\usepackage{tikz-cd} 
\usepackage{enumerate}

\graphicspath{{./Figures/}}
\usepackage{wrapfig}
\newcommand{\R}{\ensuremath{\mathbb{R}}}

\newcommand{\N}{\ensuremath{\mathbb{N}}}

\newcommand{\s}{\Sigma}

\newcommand{\e}{\varepsilon}

\newcommand{\Cr}{\mathcal{C}^{r}}
\newcommand{\Xr}{\chi^{r}}
\newcommand{\Or}{\Omega^{r}}
\newcommand{\rn}[1]{\mathbb{R}^{#1}}

\newcommand{\p}{\varphi}
\newcommand{\ag}{\alpha}

\newcommand{\cg}{\gamma}

\newcommand{\sgn}{\textrm{sgn}}

\newtheorem {theorem} {Theorem} 
\newtheorem {prop}  {Proposition}

\newtheorem {lemma}  {Lemma}
\newtheorem {definition}  {Definition}
\newtheorem {remark}  {Remark}

\makeatletter
\DeclareFontFamily{U}{tipa}{}
\DeclareFontShape{U}{tipa}{m}{n}{<->tipa10}{}
\newcommand{\arc@char}{{\usefont{U}{tipa}{m}{n}\symbol{62}}}%

\newcommand{\arc}[1]{\mathpalette\arc@arc{#1}}

\newcommand{\arc@arc}[2]{%
	\sbox0{$\m@th#1#2$}%
	\vbox{
		\hbox{\resizebox{\wd0}{\height}{\arc@char}}
		\nointerlineskip
		\box0
	}%
}
\makeatother

\definecolor{verde}{rgb}{0.0,0.5,0.0}
\definecolor{azul}{rgb}{0,0,128}
\definecolor{roxo}{rgb}{0.44,0.16,0.39}
\definecolor{vinho}{rgb}{0.5,0.0,0.13}
\definecolor{lilas1}{rgb}{0.6,0.33,0.73}
\definecolor{rosa}{rgb}{0.84,0.04,0.33}
\definecolor{mostarda}{rgb}{0.91,0.41,0.17}
\definecolor{mostarda2}{rgb}{1.0,0.66,0.07}

\newtheorem {mtheorem} {Theorem} 

\begin{document}

\title[Chains in $3D$ Filippov Systems]{Chains in $3D$ Filippov Systems: A Chaotic Phenomenon}

\author[O. M. L. Gomide]{Ot\'avio M. L. Gomide}
\address[OMLG]{Department of Mathematics, UFG, IME\\ Goi\^ania-GO, 74690-900, Brazil}
\email[Corresponding author]{otaviomarcal@ufg.br}

\author[M. A. Teixeira]{Marco A. Teixeira}
\address[MAT]{Department of Mathematics, Unicamp, IMECC\\ Campinas-SP, 13083-970, Brazil}
\email{teixeira@ime.unicamp.br}



\maketitle

\begin{abstract}
This work is devoted to the study of global connections between typical generic singularities, named $T$-singularities, in piecewise smooth dynamical systems. Such a singularity presents the so-called nonsmooth diabolo, which consists on a pair of invariant cones emanating from it. 

We analyze global features arising from the communication between the branches of a nonsmooth diabolo of a $T$-singularity and we prove that, under generic conditions, such communication leads to a chaotic behavior of the system. More specifically, we relate crossing orbits of a Filippov system presenting certain crossing self-connections to a $T$-singularity, with a Smale horseshoe of a first return map associated to the system. The techniques used in this work rely on the detection of transverse intersections between invariant manifolds of a hyperbolic fixed point of saddle type of such a first return map and the analysis of the Smale horseshoe associated to it. 

From the specific case discussed in our approach, we present a robust chaotic phenomenon for which its counterpart in the smooth case seems to happen only for highly degenerate systems.

\end{abstract}


\section{Introduction}

The theory of piecewise smooth dynamical systems has been extensively studied in recent years due its applications in modeling physical phenomenon involving some disruption in the motion (see \cite{AVK66, M62} for instance). In light of this, the comprehension of phenomena arising in that scenario seems to be increasingly necessary. 

Among all the approaches considered (see \cite{F,U77}), the notion of solutions of a piecewise smooth differential system 
\begin{equation}\label{eq:filippov}
Z(p)=\left\{\begin{array}{ll}
X(p),&f(p)>0,\\
Y(p),&f(p)<0,
\end{array}\right.
\end{equation}
with a regular {\it switching manifold} $\Sigma=f^{-1}(0)$ given by the Filippov's convention (see \cite{F} for more details) is the most considered one. In this context, special attention must be paid to some singularities lying on $\Sigma,$ known as $\Sigma$-singularities.

In planar Filippov system, the $\Sigma$-singularities have been extensively studied in \cite{BLS, GTS, K, KRG, ST,T1, T6,T2, T3} and references therein. In particular, the characterization of $\Sigma$-singularities of codimensions $0$ and $1$ was done in \cite{GTS,KRG}. 

Some of these $\Sigma$-singularities present non-standard local invariant manifolds, and therefore admit global connections. In \cite{NTZ18}, the bifurcation diagram of a planar Filippov system around a self-connection at a typical singularity is studied. In \cite{AGN19}, a methodology to study unfoldings of planar Filippov systems around global connections involving $\Sigma$-singularities is provided as well as the application of such ideas to describe bifurcation diagrams of closed global connections between $\Sigma$-singularities. 

In dimension $3$, the behavior of Filippov systems around $\Sigma$-singularities becomes strongly more complicated. In fact, the characterization of local structural stability of $3D$ Filippov systems at $\Sigma$-singularities was an open problem during more than $30$ years. In \cite{CJ1,CJ2,CJn}, they provided lots of studies on these $\Sigma$-singularities and a specific one, named T-singularity, caught their attention. The lack of comprehension of the local behavior around a T-singularity was the principal obstruction to characterize local structural stability in dimension $3$. In light of this, in \cite{J2}, the analysis of some specific models having such kind of $\Sigma$-singularity is provided. Recently, in \cite{GT} , the dynamics around a T-singularity is fully described and the locally structurally stable systems in dimension $3$ is completely characterized.

In particular, it is proved in \cite{GT} that a Filippov system $Z$ which is robust at a T-singularity $p\in\Sigma$ presents a pair of (local) nonsmooth invariant cones emanating from $p$ which are foliated by crossing orbits of $Z$, they are known as the \textit{nonsmooth diabolo} associated to $p$. One branch of this diabolo is attractive and the other one is repelling. In \cite{K08,KH11}, an analysis of some models of FIlippov systems presenting invariant cones is done. In analogy to the planar case, these local invariant manifolds can be globally extended for Filippov systems and might originate global connections at $p$ which should present a non-trivial dynamics.

In \cite{GT19}, it is considered a notion of (semi-local) structural stability in $3D$ Filippov systems which cares about what happens in a neighborhood of the whole switching manifold $\Sigma$, and it is proved that such kind of T-singularity appears naturally in semi-local structurally stable Filippov systems. The characterization of structural stability of $3D$ Filippov systems in the most global comprehensive way is one of the most relevant topics in the development of the theory of piecewise smooth dynamical systems, and to achieve this goal it is crucial to understand the behavior of the invariant manifolds emanating from a T-singularity. Another types of global connections in $3D$ Filippov systems involving $\Sigma$-singularities have already been studied (see \cite{Nshil15, G18,GTloop, KH19} for instance), nevertheless, as well as the authors know, global connections involving T-singularities have not been considered in dimension $3$ yet.

The main goal of this work is to describe the global behavior of a $3D$ Filippov system $Z$ having a robust T-singularity $p$ for which there exists a communication between their invariant cones (i.e. the stable and unstable invariant manifolds of $p$ intersect). In the present paper we prove that, in a generic scenario, the communication between these cones leads to a chaotic behavior of $Z$. More specifically, we show that there exists a (global) first return map in $\Sigma$ associated to $Z$ which presents a Smale horseshoe $\Lambda\subset \Sigma$. Moreover, every orbit of $Z$ passing through $\Lambda$ is a crossing orbit and therefore, $Z$ presents an infinity of closed crossing orbits, an infinity of non-closed orbits and a recurrent crossing orbit. We highlight that the generic conditions considered in this work gives rise to a robust behavior, i. e. the dynamics of small perturbations of $Z$ presents the same characteristics as the dynamics of $Z$.

Our methods rely on the extension of the local first return map in $\Sigma$ associated to the T-singularity $p$ to a global first return map in $\Sigma$ of $Z$ and the detection of transverse intersections between invariant manifolds of a hyperbolic fixed point of saddle type of such a global first return map. 

This paper is organized as follows. In Section \ref{basicsec} we introduce some basic concepts about Filippov systems. Section \ref{mainsec} is devoted to set the problem and state the main result of this work. In order to do this, the concept of global connections between some $\s$-singularities is formalized in Section \ref{tchainssec}, some generic conditions which allow such connections to be robust are given in Section \ref{robustsec} and the main result is stated in Section \ref{mainsec} as well as some of their consequences. Section \ref{ferraduraprova} is devoted to the proof of the main result, and a model presenting two robust connection between two T-singularities is given in Section \ref{model_sec}. Finally, in Section \ref{further_sec}, some further directions of this problem are pointed out.

\section{Basic Concepts}\label{basicsec}

Let $M$ be an open bounded  connected set of $\rn{3}$ and let $f:M\rightarrow \rn{}$ be a smooth function having $0$ as a regular value. Therefore, $\s=f^{-1}(0)$ is an embedded codimension one submanifold of $M$ which splits it in the sets $M^{\pm}=\{p\in M; \pm f(p)>0\}$.

A \textbf{germ of vector field} of class $\Cr$ at a compact set $\Lambda\subset M$ is an equivalence class $\widetilde{X}$ of $\Cr$ vector fields defined in a neighborhood of $\Lambda$. More specifically, two $\Cr$ vector fields $X_{1}$ and $X_{2}$ are in the same equivalence class if:
\begin{itemize}
	\item $X_{1}$ and $X_{2}$ are defined in neighborhoods $U_{1}$ and $U_{2}$ of $\Lambda$ in $M$, respectively;
	\item There exists a neighborhood $U_3$ of $\Lambda$ in $M$ such that $U_{3}\subset U_{1}\cap U_{2}$;
	\item $X_{1}|_{U_{3}}=X_{2}|_{U_{3}}$.
\end{itemize}
In this case, if $X$ is an element of the equivalence class $\widetilde{X}$, then $X$ is said to be a representative of $\widetilde{X}$.  
The set of germs of vector fields of class $\Cr$ at $\Lambda$ will be denoted by $\Xr(\Lambda)$, or simply $\Xr$. For the sake of simplicity, a germ of vector field $\widetilde{X}$ will be referred simply by its representative $X$.

Analogously, a \textbf{germ of piecewise smooth vector field} of class $\Cr$ at a compact set $\Lambda\subset M$ is an equivalence class  $\widetilde{Z}=(\widetilde{X},\widetilde{Y})$ of pairwise $\Cr$ vector fields defined as follows: $Z_{1}=(X_{1},Y_{1})$ and $Z_{2}=(X_{2},Y_{2})$ are in the same equivalence class if, and only if,
\begin{itemize}
	\item $X_{i}$ and $Y_{i}$ are defined in neighborhoods $U_{i}$ and $V_{i}$ of $\Lambda$ in $M$, respectively, $i=1,2$.
	\item There exist neighborhoods $U_{3}$ and $V_{3}$ of $\Lambda$ in $M$ such that $U_{3}\subset U_{1}\cap U_{2}$ and $V_{3}\subset V_{1}\cap V_{2}$.
	\item $X_{1}|_{U_{3}\cap \overline{M^{+}}}=X_{2}|_{U_{3}\cap\overline{M^{+}}}$ and $Y_{1}|_{V_{3}\cap \overline{M^{-}}}=Y_{2}|_{V_{3}\cap\overline{M^{-}}}$.
\end{itemize}
In this case,  if $Z=(X,Y)$ is an element of the equivalence class $\widetilde{Z}$, then $Z$ is said to be a representative of $\widetilde{Z}$.
The set of germs of piecewise smooth vector fields of class $\Cr$ at $\Lambda$ will be denoted by $\Or(\Lambda)$, or simply $\Or$.

We emphasize that the germ language is used due to its effectiveness to describe local and semi-local phenomena.

If $Z=(X,Y)\in\Or$ then a \textbf{piecewise smooth vector field} is defined in some neighborhood $V$ of $\Lambda$ in $M$ as
\begin{equation}\label{filippov}
Z(p)=F_1(p)+\sgn(f(p)) F_2(p),
\end{equation}
where $F_1(p)=\frac{X(p)+Y(p)}{2}$ and $F_2(p)=\frac{X(p)-Y(p)}{2}$.

The \textbf{Lie derivative} $Xf(p)$ of $f$ in the direction of the vector field $X\in \Xr$ at $p\in\Sigma$ is defined as $Xf(p)=\langle X(p), \nabla f(p)\rangle$. Accordingly, the \textbf{tangency set} between $X$ and $\s$ is given by $S_{X}=\{p\in\Sigma;\ Xf(p)=0\}$.

\begin{remark}
	Notice that the Lie derivative is well-defined for a germ $\widetilde{X}\in\Xr$ since all the elements in this class coincide in $\s$.
\end{remark}

For $X_{1},\cdots, X_{k}\in \Xr$, the higher order Lie derivatives of $f$ are defined recurrently as
$$X_{k}\cdots X_{1}f(p)=X_{k}(X_{k-1}\cdots X_{1}f)(p),$$
i.e. $X_{k}\cdots X_{1}f(p)$ is the Lie derivative of the smooth function $X_{k-1}\cdots X_{1}f$ in the direction of the vector field $X_{k}$ at $p$. In particular, $X^{k}f(p)$ denotes $X_{k}\cdots X_{1}f(p)$, where $X_{i}=X$, for $i=1,\cdots,k$.

For a piecewise smooth vector field $Z=(X,Y)$ the switching manifold $\Sigma$ is generically the closure of the union of the following three distinct open regions:
\begin{itemize}
	\item Crossing Region: $\Sigma^{c}(Z)=\{p\in \Sigma;\ Xf(p)Yf(p)>0\}.$
	\item Stable Sliding Region: $\Sigma^{ss}(Z)=\{p\in \Sigma;\ Xf(p)<0,\ Yf(p)>0\}.$
	\item Unstable Sliding Region: $\Sigma^{us}(Z)=\{p\in \Sigma;\ Xf(p)>0,\ Yf(p)<0\}.$			
\end{itemize}	

\begin{remark}
	If there is no misunderstanding, the dependence of these regions on $Z$ will be omitted. In addition, $\s$ can be denoted by $\s(Z)$, in order to distinguish the regions of $\s$ corresponding to $Z$, when necessary. 
\end{remark}

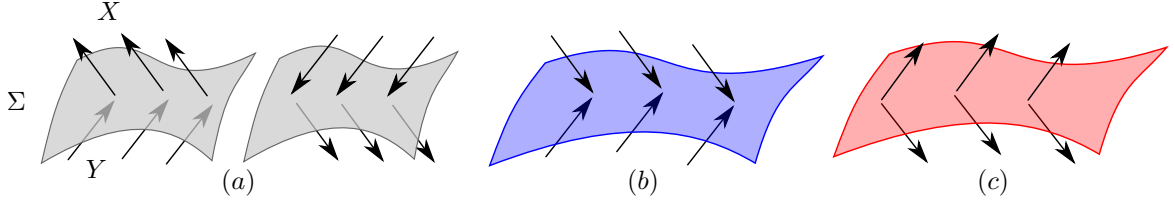
\begin{figure}[h!]
	\centering
	\bigskip
	\begin{overpic}[width=15cm]{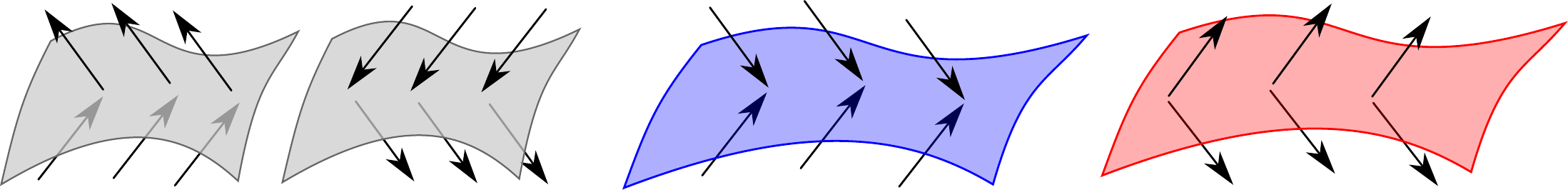}
		\put(16,-2){{\footnotesize $(a)$}}		
		\put(52,-2){{\footnotesize $(b)$}}		
		\put(83,-2){{\footnotesize $(c)$}}	
		\put(-3,5){{\footnotesize $\s$}}	
		\put(5,13){{\footnotesize $X$}}	
		\put(4,-1){{\footnotesize $Y$}}							
		
	\end{overpic}
	\bigskip
	\caption{Regions in $\s$: $\s^{c}$ in $(a)$, $\s^{ss}$ in $(b)$ and $\s^{us}$ in $(c)$.   }	\label{dicfig}
\end{figure}

The tangency set of $Z$ will be referred as $S_{Z}=S_{X}\cup S_{Y}$. Notice that $\s$ is the disjoint union $\s^{c}\cup \s^{ss}\cup \s^{us}\cup S_{Z}$. Herein, $\s^{s}=\s^{ss}\cup \s^{us}$ is called \textbf{sliding region} of $Z$. See Figure \ref{dicfig}.

The concept of solution of $Z$ follows the Filippov's convention (see, for instance, \cite{F,GTS,T8}). The local solution of $Z=(X,Y)\in \Or$ at $p\in \Sigma^{s}$ is given by the \textbf{sliding vector field}
\begin{equation}
F_{Z}(p)=\frac{1}{Yf(p)-Xf(p)}\left(Yf(p)X(p)-Xf(p)Y(p)\right).
\end{equation}
Notice that $F_{Z}$ is a $\Cr$ vector field tangent to $\s^{s}$. The critical points of $F_{Z}$ in $\s^{s}$ are called \textbf{pseudo-equilibria} of $Z$.

\begin{definition}\label{NSVF}
	We defined the \textbf{normalized sliding vector field} $F_{Z}^{N}$ of $Z$ by
	\begin{equation}
	F_{Z}^{N}(p)=Y f(p)X (p)-X f(p)Y(p),
	\end{equation}	
	for every $p\in\s^{s}$. 
\end{definition}

Notice that $F_Z^N$ is also a $\Cr$ vector field  tangent to $\s^{s}$.

\begin{remark}
	The normalized sliding vector field can be $\Cr$ extended beyond the boundary of $\s^s$. In addition, if $R$ is a connected component of $\s^{ss}$, then $F_{Z}^{N}$ is a re-parameterization of $F_{Z}$ in $R$, and so the phase portraits of both coincide. If $R$ is a connected component of $\s^{us}$, then $F_{Z}^{N}$ is a (negative) re-parameterization of $F_{Z}$ in $R$, then they have the same phase portrait, but the orbits are oriented in opposite direction.
\end{remark}

If $p\in\s^c$, then the orbit of $Z=(X,Y)\in\Or$ at $p$ is defined as the concatenation of the orbits of $X$ and $Y$ at $p$. Nevertheless, if $p\in\s\setminus\s^c$, then it may occur a lack of uniqueness of solutions. In this case,  the flow of $Z$ is multivalued and any possible trajectory passing through $p$ originated by the orbits of $X$, $Y$ and $F_Z$ is considered as a solution of $Z$. More details can be found in \cite{F,GTS}.

In the following definition, we introduce the so-called $\s$-singularities of a Filippov system.

\begin{definition}\label{sigmasing}
	Let $Z=(X,Y)\in\Or$, a point $p\in \s$ is said to be:
	\begin{enumerate}[i)]
		\item a \textbf{tangential singularity} of $Z$ provided that $Xf(p)Yf(p)=0$ and $X(p), Y(p)\neq 0$;
		\item a \textbf{$\s$-singularity} of $Z$ provided that $p$ is either a tangential singularity, an equilibrium of $X$ or $Y$, or a pseudo-equilibrium of $Z$. 
	\end{enumerate}
\end{definition}

\begin{remark}
	A point $p\in\s$ which is not a $\s$-singularity of $Z$ is also referred as a \textbf{regular-regular} point of $Z$.
\end{remark}

We say that $\gamma$ is a \textbf{regular orbit} of $Z=(X,Y)$ if it is a piecewise smooth curve such that $\gamma\cap M^{+}$ and $\gamma\cap M^{-}$ are unions of regular orbits of $X$ and $Y$, respectively, and  $\gamma\cap\s\subset\s^{c}$.

\begin{definition}\label{elementarydef}
	Let $Z=(X,Y)\in\Or$. A tangential singularity $p\in \s$ is said to be a \textbf{fold-fold singularity} if $Xf(p)= 0$, $X^{2}f(p)\neq0, Yf(p)= 0$, $Y^{2}f(p)\neq 0$ and $S_{X}\pitchfork S_{Y}$ at $p$. In addition, a fold-fold singularity $p$ of $Z=(X,Y)\in\Or$ is said to be:
	\begin{itemize}
		\item a \textbf{visible fold-fold} if $X^{2}f(p)>0$ and $Y^{2}f(p)<0$;
		\item an \textbf{invisible-visible fold-fold} if $X^{2}f(p)<0$ and $Y^{2}f(p)<0$;	
		\item a \textbf{visible-invisible fold-fold} if $X^{2}f(p)>0$ and $Y^{2}f(p)>0$;			
		\item an \textbf{invisible fold-fold} if $X^{2}f(p)<0$ and $Y^{2}f(p)>0$, in this case, $p$ is also called a \textbf{T-singularity}.				
	\end{itemize}
\end{definition}

\begin{figure}[h!]
	\centering
	\bigskip
	\begin{overpic}[width=15cm]{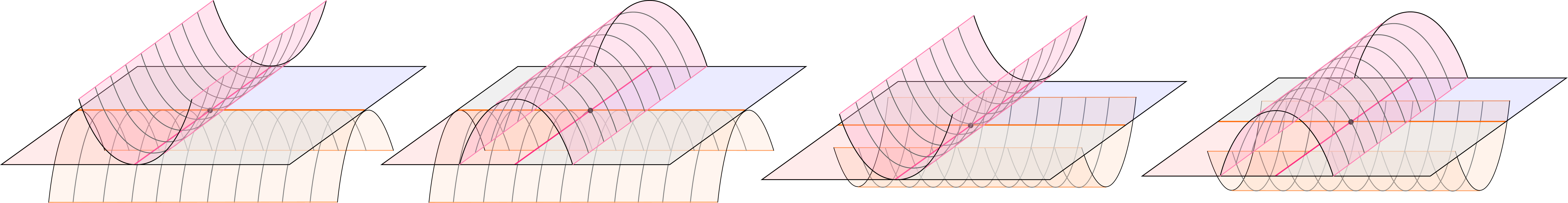}
		\put(11,-3){{\footnotesize (a)}}						
		\put(35,-3){{\footnotesize (b)}}
		\put(60,-3){{\footnotesize (c)}}	
		\put(84,-3){{\footnotesize (d)}}	
		\put(0,5){{\footnotesize $\s$}}		
		\put(7,11){{\footnotesize $X$}}	
		\put(1,0){{\footnotesize $Y$}}																					
	\end{overpic}
	\bigskip
	\caption{Fold-Fold Singularity: (a) Hyperbolic, (b,c) Parabolic and (d) Elliptic.}	\label{types_fig}
\end{figure}
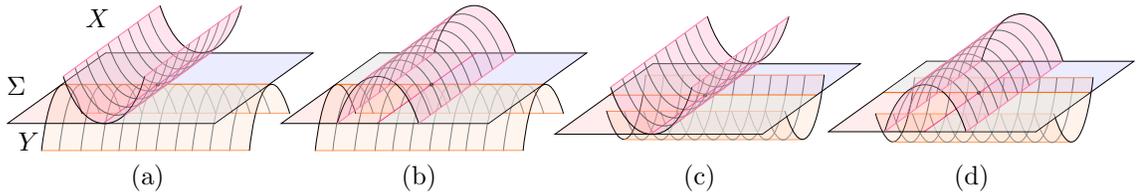

\begin{remark}
	Notice that the visible-invisible case can be obtained from the invisible-visible one by performing an orientation reversing change of coordinates. Also, we refer to a visible, invisible-visible/visible-invisible, invisible fold-fold point as a hyperbolic, parabolic, elliptic fold-fold point, respectively. See Figure \ref{types_fig}.
\end{remark}

\section{Setting the Problem and Main Result} \label{setsec}

This section is devoted to introduce the concept of T-chains of Filippov systems and to present the main result achieved in this work.

\subsection{T-chains}\label{tchainssec}
From \cite{GT}, we have that a Filippov system is locally structurally stable at certain specific types of fold-fold singularities and it presents local invariant manifolds at such points. As we have mentioned before, the knowledge of the dynamics of Filippov systems around global connections involving these points plays a crucial role in the attempt to characterize the structurally stable Filippov systems. In light of this, we formalize such connections in the following definition.

\begin{definition}
	Let $Z_0=(X_0,Y_0)\in\Or$ having fold-fold singularities $p_0,q_0\in\s$ ($p_0=q_0$ is also considered) and let $-\infty\leq a<b\leq \infty$. An oriented piecewise smooth curve $\Gamma: (a,b)\rightarrow M$ is said to be a \textbf{fold-fold connection} of $Z_0$ between $p_0$ and $q_0$ if it satisfies the following conditions.
	\begin{enumerate}[i)]
		\item $\mathrm{Im}(\Gamma)\cap M^+$ (resp. $\mathrm{Im}(\Gamma)\cap M^-$) is a union of orbits of $X_0$ (resp. $Y_0$).
		\item $\mathrm{Im}(\Gamma)\cap \s\subset \s^{c}$.
		\item $\displaystyle\lim_{t\rightarrow a} \Gamma(t)= p_0$ and $\displaystyle\lim_{t\rightarrow b} \Gamma(t)= q_0$.
	\end{enumerate}
\end{definition}

Motivated by the definition of polycycles in Filippov systems (see \cite{AGN19}), we introduce the following concept.

\begin{definition}
	Consider $Z_0\in\Or$ having a finite number of fold-fold singularities $p_i\in\s$, $i=1,\cdots,k$. We say that $\gamma\subset M$ is a \textbf{fold-fold chain of order $k$} of $Z_0$ if $$\gamma=\{p_1,\cdots,p_k\}\displaystyle\cup_{i=1}^k \mathrm{Im}(\Gamma_i),$$
	where $\Gamma_i$ is a fold-fold connection between $p_i$ and $p_{i+1}$, $i=1,\cdots, k$, where $p_{k+1}=p_1$, and either one of the following conditions is satisfied.
	\begin{enumerate}[i)]
		\item $\Gamma_i$ is an oriented piecewise smooth curve from $p_i$ to $p_{i+1}$, $i=1,\cdots,k$.
		
		\item $\Gamma_i$ is an oriented piecewise smooth curve from $p_{i+1}$ to $p_{i}$, $i=1,\cdots,k$.
	\end{enumerate}	
\end{definition}

\begin{figure}[h!]
	\centering
	\bigskip
	\begin{overpic}[width=12cm]{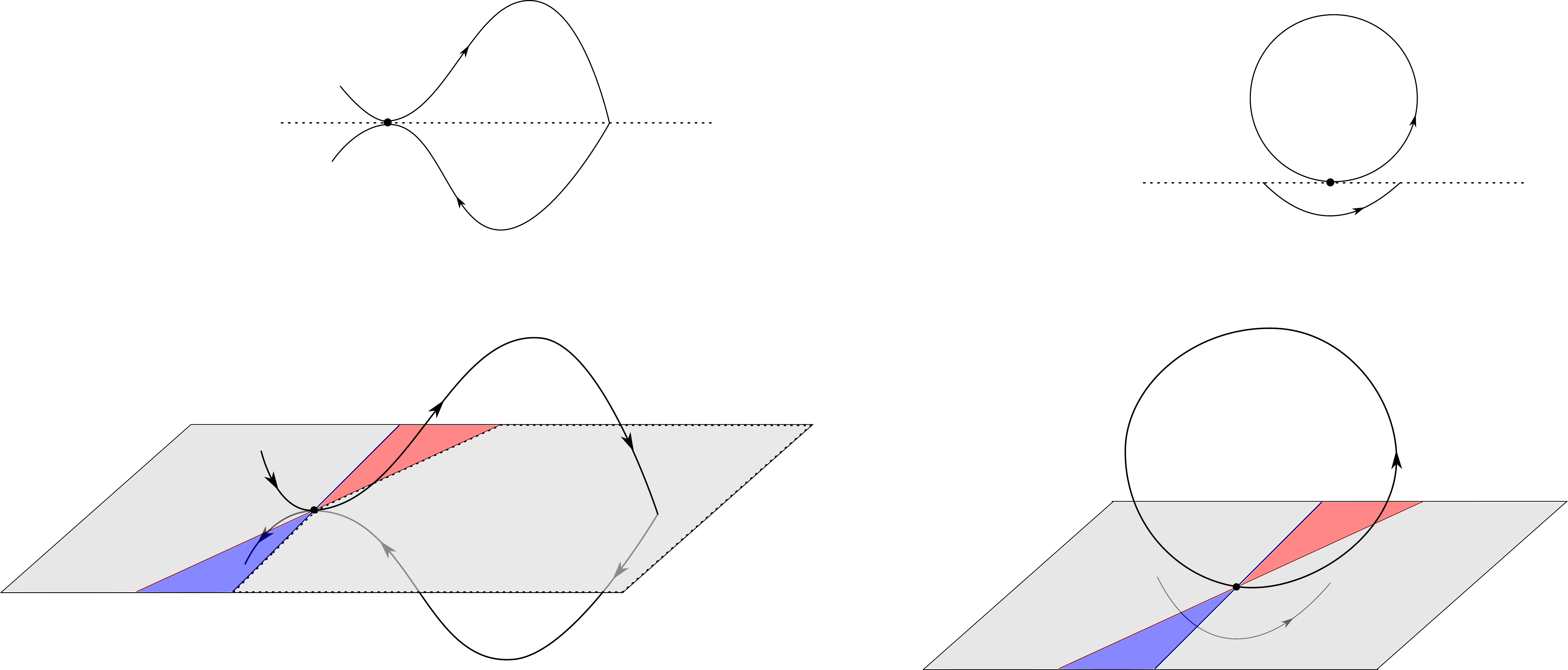}
		\put(4,12){{\footnotesize $\s$}}
		\put(10,2){{\footnotesize $\s^{ss}$}}
		\put(30,16){{\footnotesize $\s^{us}$}}	
		\put(16,37){{\footnotesize $\s$}}	
		\put(30,24){{\footnotesize $(a)$}}
		\put(83,24){{\footnotesize $(b)$}}	
		\put(25,-3){{\footnotesize $(c)$}}	
		\put(80,-3){{\footnotesize $(d)$}}				
		%
	\end{overpic}
	\bigskip
	\caption{Planar $\s$-polycycles passing through a visible-visible fold-fold singularity $(a)$ and a  visible-invisible fold-fold singularity $(b)$, and tridimensional fold-fold chains passing through a visible-visible fold-fold singularity $(c)$ and a  visible-invisible fold-fold singularity $(d)$.}	\label{chainsfig}
\end{figure}
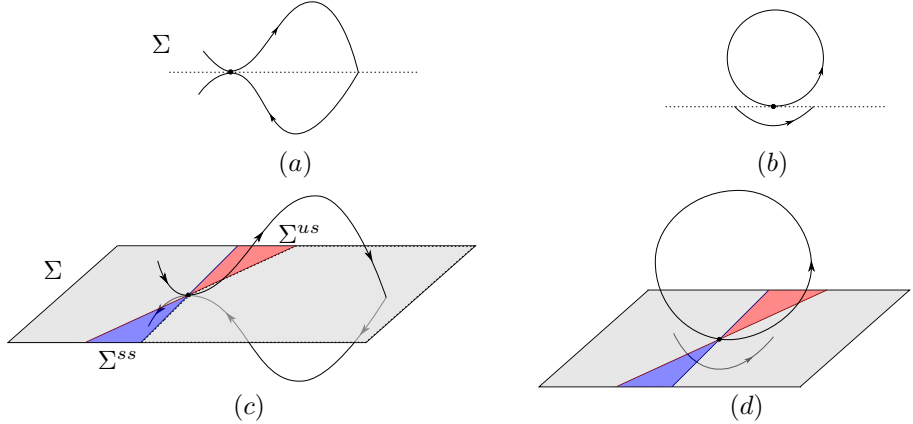

Notice that, the above definition of a fold-fold chain generalizes the concept of a $\s$-polycycle having only fold-fold singularities in the planar Filippov systems introduced in \cite{AGN19} to $3D$ Filippov systems. Figure \ref{chainsfig} illustrates some fold-fold chains.

As far as the authors know, there is a lack of results in the literature concerning this kind of object, maybe due to the difficult inherent to the problem. In fact, $3D$ Filippov systems exhibit a rich local dynamics at fold-fold singularities which is still hard to comprehend, and thus, the understanding of global phenomena involving such objects becomes even harder.

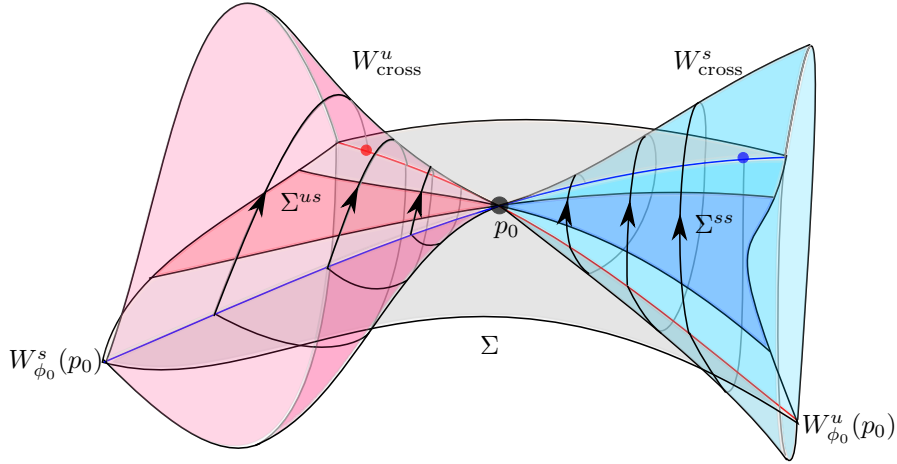
\begin{figure}[h!]
	\centering
	\bigskip
	\begin{overpic}[width=9.5cm]{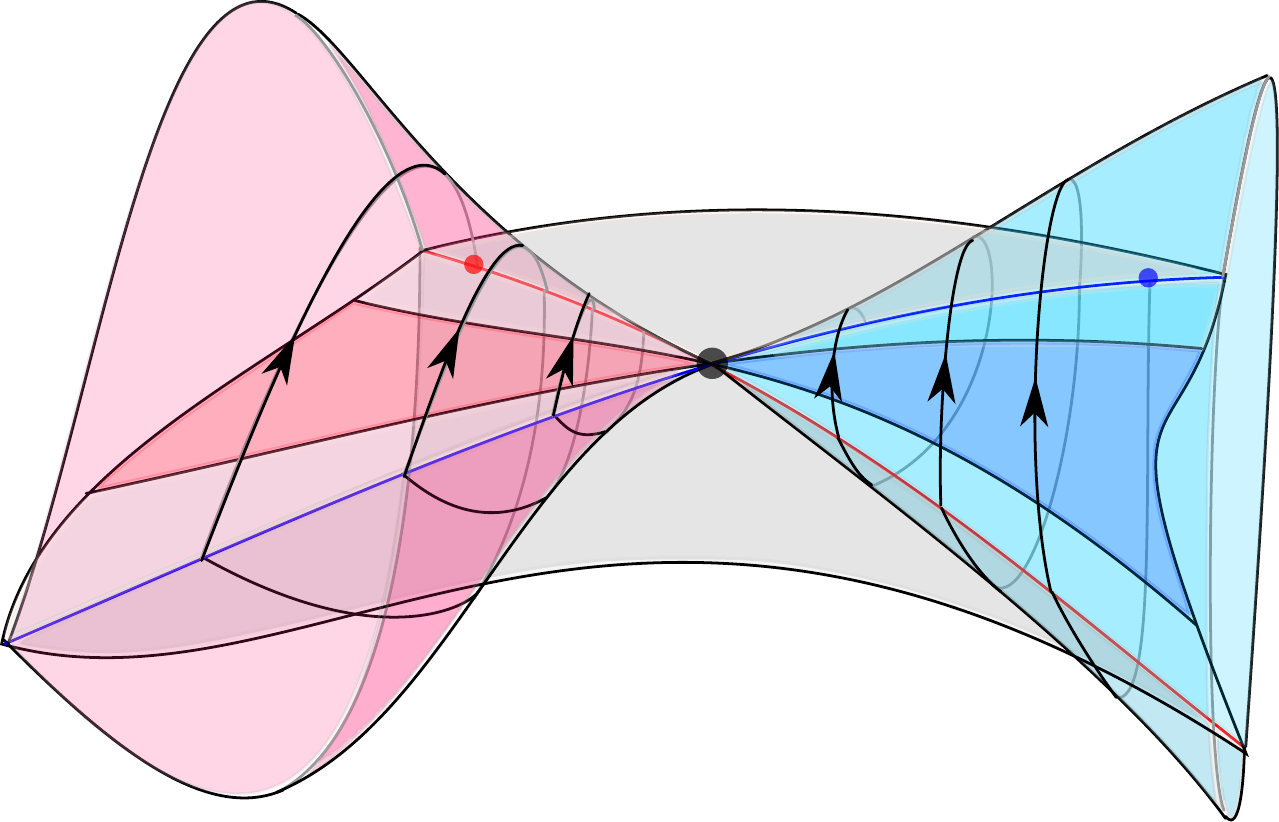}							
		\put(55,32){{\footnotesize $p_0$}}		
		\put(25,35){{\footnotesize $\s^{us}$}}		
		\put(83,32){{\footnotesize $\s^{ss}$}}	
		\put(53,15){{\footnotesize $\s$}}		
		\put(98,4){{\footnotesize $W_{\phi_0}^{u}(p_0)$}}		
		\put(-13,13){{\footnotesize $W_{\phi_0}^{s}(p_0)$}}	
		\put(35,55){{\footnotesize $W_{\mathrm{cross}}^{u}$}}
		\put(80,55){{\footnotesize $W_{\mathrm{cross}}^{s}$}}				
	\end{overpic}
	\bigskip
	\caption{Nonsmooth diabolo $\mathcal{N}(p_0)$ at a stable $T$-singularity $p_0$ of $Z_0$.}
	\label{nonsmoothdiabfig}	
\end{figure}

One of the most challenging types of fold-fold singularity is the invisible one, also known as T-singularity. In this case, if $Z_0=(X_0,Y_0)\in\Or$ has a T-singularity at $p_0$ (see Figure \ref{nonsmoothdiabfig}), then there are two involutions $\s\subset\rn{2}$, $\phi_{X_0},\phi_{Y_0}:(\s,p_0)\rightarrow (\s,p_0)$, which are  induced by the orbits of $X_0$ and $Y_0$ near $p_0$. It means that the points $p\in\Sigma$ and $\phi_{X_0}(p)\in\Sigma$ (resp. $\phi_{Y_0}(p)$) are connected by a segment of orbit of $X_0$ (resp. $Y_0$) which is contained in $M^+$ (resp. $M^-$). These involutions give rise to a $\Cr$ germ of first return map $\phi_0:(\s,p_0)\rightarrow(\s,p_0)$ given by $\phi_0=\phi_{X_0}\circ\phi_{Y_0}$. See \cite{GT} for more details.

For simplicity, we say that $p_0$ is a \textbf{stable T-singularity} of $Z_0$ if, and only if, $p_0$ is a T-singularity for which $\phi_0$ has a hyperbolic fixed point of saddle type at $p_0$ with both local invariant manifolds $W^{u,s}_{\phi_0}(p_0)$ of $\phi_0$ at $p_0$ contained in $\s^c$. In \cite{GT}, it is proved that $Z_0\in\Or$ is locally structurally stable at a T-singularity $p_0$ if, and only if, $p_0$ is a stable T-singularity of $Z_0$.  

Also, if $Z_0$ has a stable T-singularity at $p_0$, then there exists a (local) invariant cone $\mathcal{N}(p_0)$ with vertex at $p_0$ which is filled up with crossing orbits of $Z_0$. In addition, $\mathcal{N}(p_0)$ is piecewise smooth and $\mathcal{N}(p_0)\cap\s=W^{u}_{\phi_0}(p_0)\cup W^{s}_{\phi_0}(p_0)$. Denote the stable and unstable branches of $\mathcal{N}(p_0)$ by $W_{\textrm{cross}}^s$ and $W_{\textrm{cross}}^u$, respectively. The existence of such cone $\mathcal{N}(p_0)$ has been proved in \cite{GT}, and it is also referred as the nonsmooth diabolo associated to $Z_0$ at $p_0$ (see \cite{CJ1}). See Figure \ref{nonsmoothdiabfig}.

In light of this discussion, we have seen that a Filippov system has local crossing invariant manifolds (stable and unstable) at a stable T-singularity. Therefore, a natural question arises in such scenario: what kind of dynamics is originated from the global extension of these local invariant manifolds?

\begin{definition}
	Consider $Z_0\in\Or$. We say that $\gamma\subset M$ is a \textbf{T-chain} of $Z_0$ if $\gamma$ is a fold-fold chain of order $1$ of $Z_0$ having a unique stable T-singularity of $Z_0$.
\end{definition}

In this work, we  study the dynamics of Filippov systems around T-chains through a semi-local analysis at this global connection. We highlight that T-chains are the simplest fold-fold chains having stable T-singularities and we restrict our study to this case because, even in this situation, a Filippov system displays a very complicated dynamics in the presence of such an object.

\subsection{Robustness Conditions}\label{robustsec}
\begin{figure}[h!]
	\centering
	\bigskip
	\begin{overpic}[width=8cm]{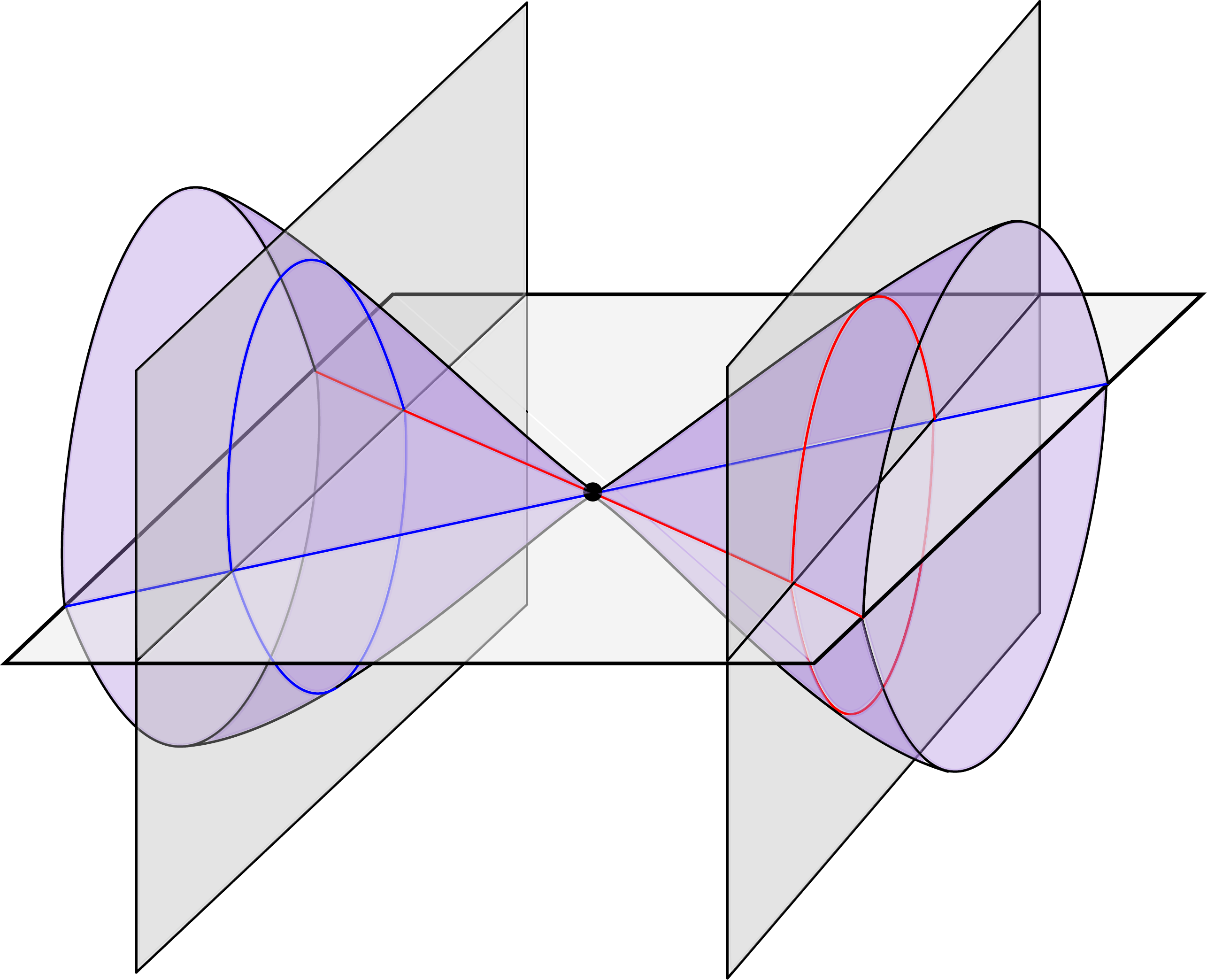}							
		\put(47,36){{\footnotesize $p_0$}}		
		\put(90,60){{\footnotesize $W_{\mathrm{cross}}^{u}$}}
		\put(0,65){{\footnotesize $W_{\mathrm{cross}}^{s}$}}				
		\put(29,59){{\footnotesize $\mathcal{C}^s$}}
		\put(70,58){{\footnotesize $\mathcal{C}^u$}}
		\put(45,77){{\footnotesize $\tau^s$}}
		\put(87,77){{\footnotesize $\tau^u$}}
		\put(45,20){{\footnotesize $\s$}}												
	\end{overpic}
	\bigskip
	\caption{Sections $\tau^u$ and $\tau^s$.}
	\label{transsec}	
\end{figure}
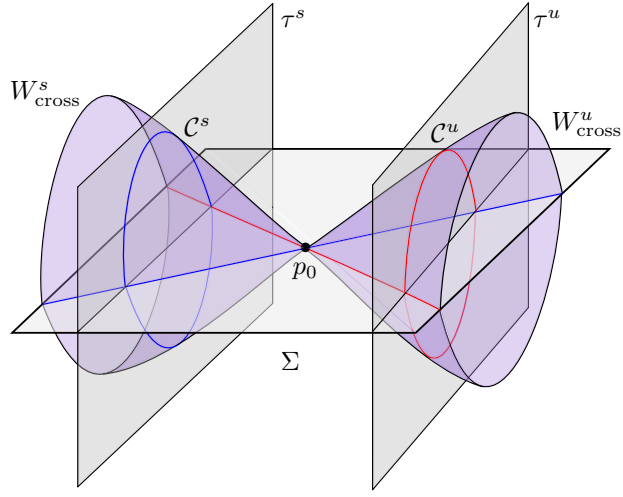

Let $Z_0=(X_0,Y_0)\in\Or$ having a stable T-singularity at $p_0$. For $\star=u,s$, let $\tau^{\star}$ be a  section  such that $\tau^{\star}\cap W_{\textrm{cross}}^{\star}=\mathcal{C}^{\star}$ is a piecewise smooth closed curve homotopic to a circle which is nonsmooth only at (the two) points belonging to $\mathcal{C}^{\star}\cap\s$. Assume that $\tau^{\star}$ is transverse to the flow of $Z_0$ at the points of $\mathcal{C}^{\star}$ and that $\tau^{\star}$, $W_{\textrm{cross}}^{\star}$ and $\s$ are in general position (see Figure \ref{transsec}). Also, consider that $\tau^{\star}$ is contained in a neighborhood $V_{3D}$ of $p_0$ in $M$, for which the local first return map $\phi_0:V\rightarrow \s$ associated to $Z_0$ at $p_0$ is defined in $V=V_{3D}\cap\s$.

\begin{remark}
	By saying that $\tau^{\star}$ is transverse to the flow of $Z_0$ at points of $\mathcal{C}^{\star}$, we mean that $X_0$ (resp. $Y_0$) is transverse to $\tau^{\star}$ at each point $q\in W_{\textrm{cross}}^{\star}\cap \overline{M^+}$ (resp. $q\in W_{\textrm{cross}}^{\star}\cap \overline{M^-}$).
\end{remark}

Assume that $Z_0$ satisfies the following \textbf{(TC) conditions} :

\begin{enumerate}
	\item[($TC_1$)] $Z_0$ has a stable T-singularity at $p_0\in\s$;\vspace{0.2cm}
	\item[($TC_2$)] There exists a germ of diffeomorphism $\mathcal{D}: \tau^u\rightarrow \tau^s$ at $\mathcal{C}^u$,  induced by orbits of $X_0$ and $Y_0$ such that, for each $q\in \mathrm{Dom}(\mathcal{D})$ (domain of $\mathcal{D}$), $q$ and $\mathcal{D}(q)$ are connected by a crossing orbit of $Z_0$ and $\mathcal{D}(\mathcal{C}^u)=\widehat{\mathcal{C}^u}$ is a topological circle contained in $\tau^s$;	\vspace{0.2cm}
	
	\item[($TC_3$)] There exists a $Z_0$-invariant topological cylinder $\mathcal{R}$ ($2$-dimensional) connecting $\mathcal{C}^u$ and $\widehat{\mathcal{C}^u}$, which is filled up with crossing orbits of $Z_0$. Assume that $\mathcal{R}\cap\s$ is given by two compact distinct curves $\mathcal{R}^u$, $\mathcal{R}^s$ which contains the points $W^u_{\phi_0}(p_0)\cap\tau^u$ and $W^s_{\phi_0}(p_0)\cap\tau^u$, respectively. Also, consider that each crossing orbit contained in $\mathcal{R}$ does not intersect $\Sigma$ in $\mathcal{R}^{\star}$ consecutively, for $\star=u,s$.	
\end{enumerate}	

\begin{figure}[h!]
 	\centering
 	\bigskip
 	\begin{overpic}[width=14cm]{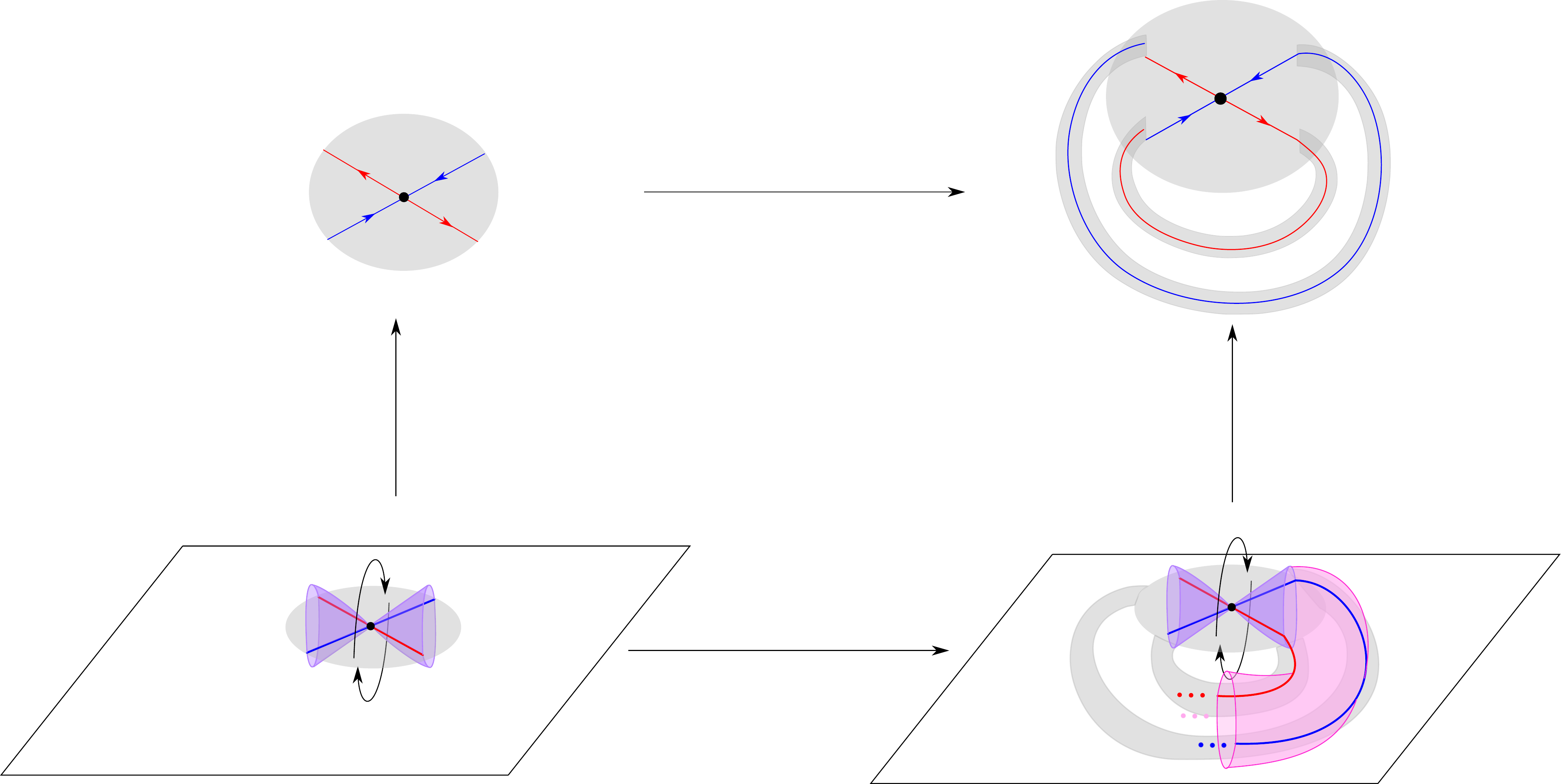}							
 		\put(22.6,8.5){{\tiny $p_0$}}
 		\put(25,36){{\tiny $p_0$}}	
 		 		\put(77,42){{\tiny $p_0$}}
 		\put(77.7,9.2){{\tiny $p_0$}}	
 		\put(30,41){{\scriptsize \textcolor{blue}{$W_{\phi_0}^{s}$}}}
 		\put(88,45){{\scriptsize \textcolor{blue}{$W_{\Phi_0}^{s}$}}}
\put(30,32){{\scriptsize \textcolor{red}{$W_{\phi_0}^{u}$}}}
\put(78,36){{\scriptsize \textcolor{red}{$W_{\Phi_0}^{u}$}}}
 		\put(19,14){{\scriptsize $\phi_{X_0}$}}
 				\put(74,16){{\scriptsize $\Phi_{X_0}$}}
 					\put(74,7){{\scriptsize $\Phi_{Y_0}$}}
 			\put(25,4){{\scriptsize $\phi_{Y_0}$}}
 		\put(4,10){{\footnotesize $\s$}}
 				\put(10,40){{\footnotesize $\s$}}
 				\put(94,40){{\footnotesize $\s$}}
 			\put(100,10){{\footnotesize $\s$}}	
 			\put(15,10){{\scriptsize $V$}}	
 			\put(17,37){{\scriptsize $V$}}	
 				\put(65,7){{\scriptsize $W$}}	
 					\put(64,40){{\scriptsize $W$}}	
 					\put(40,6){{\scriptsize Global extension of}}		
 					\put(40,3){{\scriptsize invariant cones}}	
 					\put(5,25){{\scriptsize First Return Map}}	
 					\put(5,22){{\footnotesize $\phi_0=\phi_{Y_0}\circ\phi_{X_0}$}}
 					\put(80,25){{\scriptsize First Return Map}}	
 					\put(80,22){{\footnotesize $\Phi_0=\Phi_{Y_0}\circ\Phi_{X_0}$}}	
 					 					\put(41,42){{\scriptsize Global extension of}}		
 					\put(41,39){{\scriptsize first return map}}									
 	\end{overpic}
 	\bigskip
 	\caption{Extension of the local first return map $\phi_0:V\rightarrow \s$ into the global first return map $\Phi_0:W\rightarrow \s$.}
 	\label{extensionfig}	
 \end{figure}
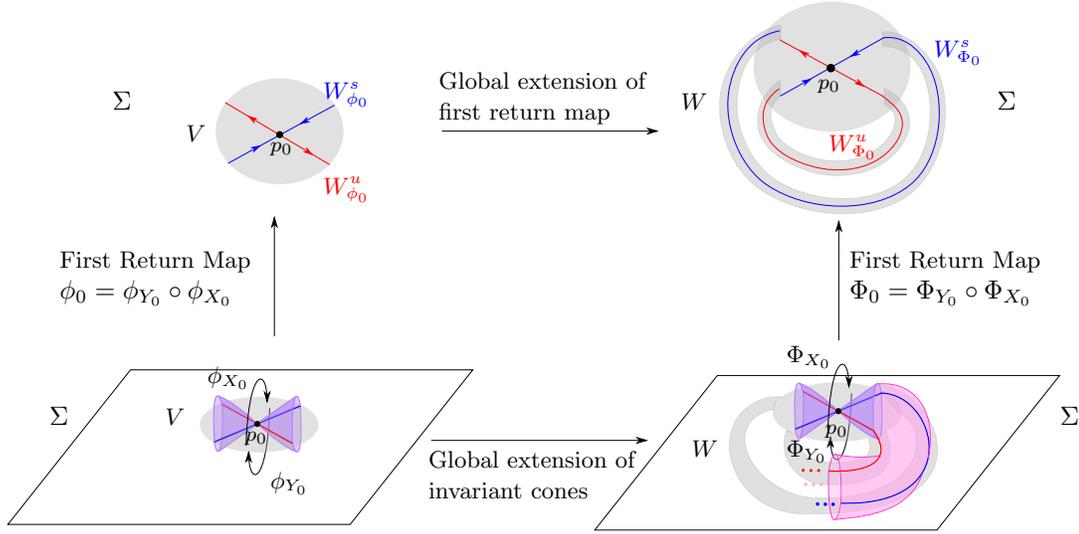

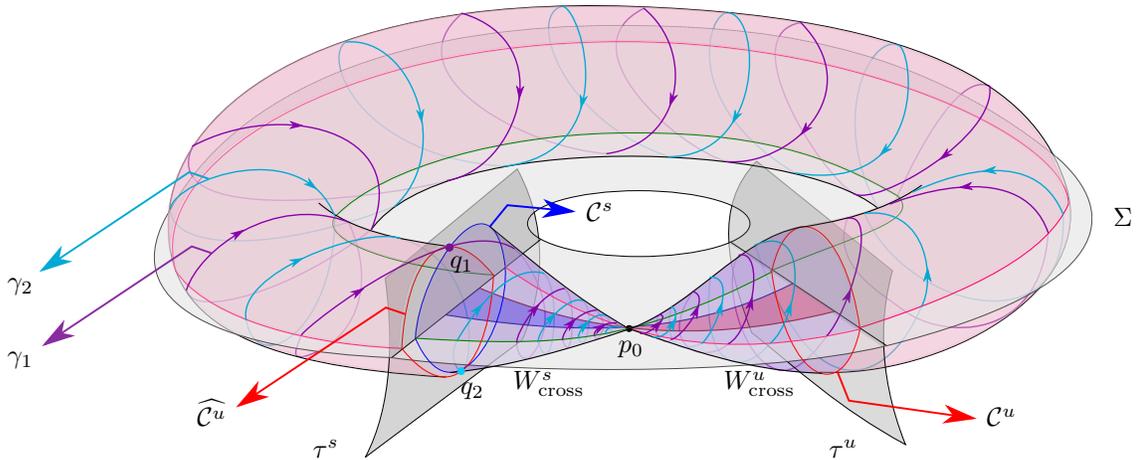
\begin{figure}[h!]
	\centering
	\bigskip
	\begin{overpic}[width=14cm]{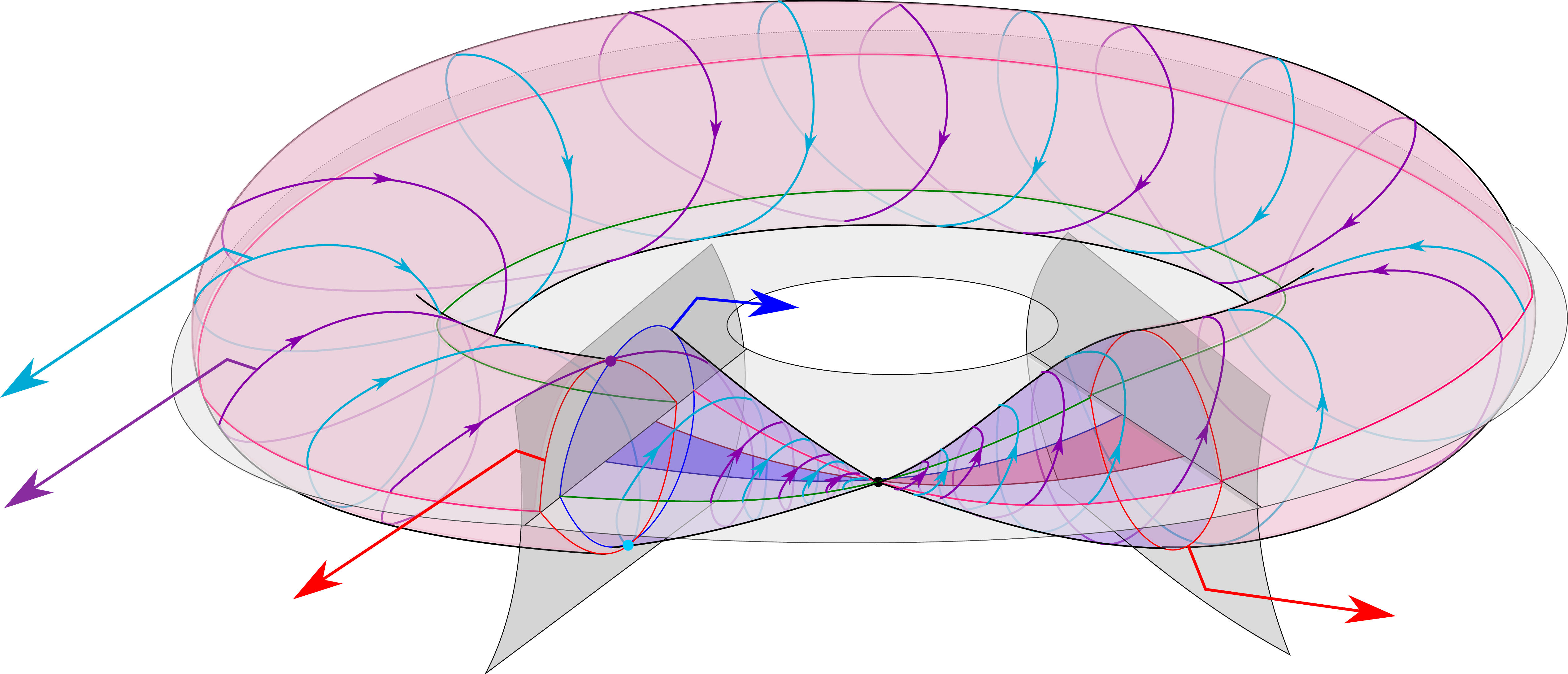}
		\put(26,0){{\footnotesize $\tau^s$}}		
		\put(75,0){{\footnotesize $\tau^u$}}		
		\put(55,10){{\footnotesize $p_0$}}			
		\put(39,18){{\footnotesize $q_1$}}		
		\put(40,6){{\footnotesize $q_2$}}
		\put(45,6.5){{\footnotesize $W^s_{\textrm{cross}}$}}
		\put(65,6.5){{\footnotesize $W^u_{\textrm{cross}}$}}				
		\put(52,22.5){{\footnotesize $\mathcal{C}^s$}}	
		\put(90,3){{\footnotesize $\mathcal{C}^u$}}	
		\put(15,3){{\footnotesize $\widehat{\mathcal{C}^u}$}}			
		\put(-3,9){{\footnotesize $\gamma_1$}}
		\put(-3,16){{\footnotesize $\gamma_2$}}				
		\put(102,22){{\footnotesize $\s$}}									
	\end{overpic}
	\bigskip
	\caption{A Filippov system $Z_0$ satisfying $(TC)$ and $(R)$ conditions having two T-chains $\gamma_1$ and $\gamma_2$ passing through $q_1$ and $q_2$, respectively.}
	\label{tcadeia}	
\end{figure}

The above set of hypotheses $(TC)$ allows us to extend the crossing invariant manifold $W^u_{\textrm{cross}}(p_0)$ of $Z_0$ through a cylinder $\mathcal{R}$ in such a way that it intersects the section $\tau^s$ at a topological circle $\widehat{C^u}$ (see Figure \ref{tcadeia}). Next lemma show that such conditions allow us to extend the local first return map $\phi_0$ of $Z_0$ at $p_0$ into a first return map $\Phi_0$ in $\s$ around  $\{p_0\}\cup\mathcal{R}^u\cap \mathcal{R}^s$, in such way that the local invariant manifolds $W^u_{\phi_0}$ and $W^s_{\phi_0}$ are extended by $\mathcal{R}^u$ and $\mathcal{R}^s$ into global invariant manifolds $W^u_{\Phi_0}$ and $W^s_{\Phi_0}$ of $\Phi_0$.  See Figure \ref{extensionfig}. 
 
\begin{lemma}[Extension] \label{extlema}Let $Z_0=(X_0,Y_0)\in\Or$ satisfying $(TC)$ conditions  and let $\phi_0=\phi_{X_0}\circ\phi_{Y_0}:V\rightarrow\s$ be its local first return map at the stable $T$-singularity $p_0$. There exists a small connected neighborhood $W$ of $\{p_0\}\cup\mathcal{R}^u\cup\mathcal{R}^s$ in $\s$ such that
	\begin{enumerate}[i)]
		\item $V\subset W$ and $W\setminus V\subset\s^c$;
		\item There exists an involution $\Phi_{Y_0}:W\rightarrow W$ induced by orbits of $Y_0$, i.e., for each $p\in W$, $p$ and $\Phi_{Y_0}(p)$ are connected through an orbit of $Y_0$ contained in $\overline{M^-}$;
		\item There exists an involution $\Phi_{X_0}:W\rightarrow \s$ induced by orbits of $X_0$;
		\item $\Phi_0=\Phi_{X_0}\circ \Phi_{Y_0}$ is a reversible mapping which is an extension of $\phi_0$. In addition, $\Phi_0$ has a unique hyperbolic fixed point at $p_0$ which is of saddle type.
	\end{enumerate}
	Furthermore, for $\star=u,s$, the global invariant manifold $W^\star_{\Phi_0}$ of $\Phi_0$ at $p_0$ is an extension of $W^{\star}_{\phi_0}(p_0)$ and contains the curve $\mathcal{R}^\star$.
\end{lemma}

The proof of Lemma \ref{extlema} is done in Appendix \ref{appendixA}.

Notice that $T$-chains of $Z_0$ at $p_0$ are well characterized as intersections between the topological circles $\widehat{C^u}$ and $\mathcal{C}^s$. 

\begin{prop}\label{tiposTchain}
	Let $Z_0\in\Or$ satisfying $(TC)$ conditions. The following statements hold.
	\begin{enumerate}[i)]
		\item If $\widehat{C^u}\cap\mathcal{C}^s=\emptyset$, then $Z_0$ has no $T$-chains at $p_0$;
		\item If $\widehat{C^u}=\mathcal{C}^s$, then $Z_0$ has an invariant (piecewise smooth) pinched torus at $p_0$ foliated by $T$-chains at $p_0$;
		\item If $\widehat{C^u}\cap\mathcal{C}^s={q_1,\cdots, q_K}\subset M^+\cup M^-$  and $\widehat{\mathcal{C}^u}\pitchfork \mathcal{C}^s$ at $q_i$, $i\in{1,\cdots,K}$, then $Z_0$ has $K$ distinct $T$-chains at $p_0$ and $K=2k$, for some $k\in\N$. 
	\end{enumerate}
\end{prop}

The proof of Proposition \ref{tiposTchain} is straightforward and it will be omitted.
\begin{remark}
	Notice that the topological circles $\mathcal{C}^s$ and $\widehat{\mathcal{C}^u}$ are smooth at the points $q_i,$ $1\leq i\leq 2k$, since $q_i\notin\s$. Therefore, the notion of transversality is well-defined in condition $(iii)$ of Proposition \ref{tiposTchain}.
\end{remark}

We notice that if item $(ii)$ of Proposition \ref{tiposTchain} is satisfied then the (reversible) first return map $\Phi_0$ obtained in Lemma \ref{extlema} has a homoclinic connection at $p_0$. Clearly, such a situation is not robust, since a small perturbation breaks the condition $\widehat{\mathcal{C}^u}=\mathcal{C}^s$. It is worth mentioning that results on bifurcation of reversible maps around homoclinic orbits can be used to understand what happens with these manifolds under small perturbations, nevertheless, this situation is highly degenerated and thus, it can give rise to very complicated phenomena. In \cite{DGGLS12}, the authors have studied bifurcations of homoclinic orbits of some planar reversible maps.

In order to avoid further degeneracies, we consider the following \textbf{robustness condition} on $Z_0$:

\textbf{(R) }  $\mathcal{C}^s\cap \widehat{\mathcal{C}^u}=\{q_1,\cdots,q_{2k}\}$, for some $k\in\N$, where $q_i\notin\s$ and $\mathcal{C}^s\pitchfork \widehat{\mathcal{C}^u}$ at $q_i$, $i\in{1,\cdots,2k}$. 

Without loss of generality, we consider that $k=1$ throughout this work. Also, we highlight that the condition $q_i\notin\s$ in $(R)$ and item $(iii)$ of Proposition \ref{tiposTchain} is only technical and can be dropped by extending the notion of transversality of $\mathcal{C}^s$ and $\widehat{\mathcal{C}^u}$ at points of $\s$.

Therefore, if $Z_0$ satisfies $(TC)$ and $(R)$ conditions, then $Z_0$ has two distinct $T$-chains $\gamma_1$ and $\gamma_2$ at $p_0$. We notice that, in this case, each $T$-chain is a crossing homoclinic orbit of $Z_0$, since it reaches $p_0$ only at infinite time (see Figure \ref{tcadeia}). In addition, $(TC)$ and $(R)$ conditions are persistent under small perturbations of $Z_0$, thus we have that such $T$-chains of $Z_0$ are robust in $\Or$ (i.e. they can not be destroyed for $Z$ near $Z_0$).

\subsection{Main Result}\label{mainsec}

Let $Z_0\in\Or$ satisfying $(TC)$ and $(R)$ conditions. From Lemma \ref{extlema}, we have that $Z_0$ is associated to a first return map $\Phi_0=\Phi_{X_0}\circ\Phi_{Y_0}$ induced by orbits of $X_0$ and $Y_0$. Recall that $\Phi_{X_0}$ and $\Phi_{Y_0}$ describes the foliation generated by $X_0$ and $Y_0$ in $\overline{ M^+}$ and $\overline{ M^-}$, respectively, in the sense that $x$ and $\phi_{X_0}(x)$ (resp. $\phi_{Y_0}(x)$) are connected by an orbit of $X_0$ (resp. $Y_0$) contained in $\overline{ M^+}$ (resp. $\overline{ M^-}$), for every $x\in W$. It follows that the foliation generated by all orbits of $Z_0$ is well described by $\Phi_0$.

Hence, the behavior of the foliation generated by $Z_0$ around the $T$-chains $\gamma_1$ and $\gamma_2$ at $p_0$, can be determined by the dynamics of the first return map $\Phi_0$. Notice that the map $\Phi_0$ does not care about how two pieces of orbits of $X_0$ and $Y_0$ are concatenated. In fact, the oriented itinerary of a point $p\in\Sigma$ through $\Phi_0$  (i.e. the points $\Phi_0^k(p),\ k\in\N$) can be related to a piecewise smooth curve having pieces of orbits of $X_0$ and $Y_0$ which are concatenated in opposite directions. In light of this, we introduce the following definition.

\begin{definition}\label{pseudo}
	We say that a piecewise smooth curve $\gamma$ is a \textbf{pseudo-orbit} of $Z_0=(X_0,Y_0)$ if it satisfies the following conditions
	\begin{enumerate}[i)]
		\item $\gamma\cap \overline{ M^+}$ is tangent to $X_0$;
		\item $\gamma\cap \overline{ M^-}$ is tangent to $Y_0$;
		\item There exists at least a point $p\in\gamma\cap\s$ such that $X_0f(p)Y_0f(p)<0$.
	\end{enumerate}
\end{definition}

Hence, the orbits of $\Phi_0$ are associated to crossing orbits and pseudo-orbits of $Z_0$, and vice-versa. Also, notice that if $\left(\Phi_0^n(x)\right)_{n\in\N}$ corresponds to a crossing orbit of $Z_0 $, then the evolution of $x$ through $\Phi_0$ might not coincide with the evolution in time of the corresponding orbit of $Z_0$.

Finally, we state the main result of this paper. 

\begin{mtheorem}\label{ferradura}
	Let $Z_0\in\Or$ satisfying $(TC)$ and $(R)$ conditions and let $\gamma_1$ and $\gamma_2$ be the two distinct $T$-chains of $Z_0$ at $p_0$ (provided by these conditions). Let $U$ be an arbitrarily small neighborhood of the stable T-singularity $p_0$ in $\s$. Then:
	\begin{enumerate}
		\item There exists $n\in\N$ such that $\Phi_0^{n}$ admits a Smale horseshoe $\Delta$ in $U$.
		\item The hyperbolic invariant set $\Lambda$ in the horseshoe $\Delta$ always contains a point $\widehat{q_i}\in\gamma_i\cap\s$, for $i=1,2$.
		\item Every orbit of $Z_0$ passing though a point of $\Lambda$ is a crossing orbit.
	\end{enumerate}
\end{mtheorem}

The proof of Theorem \ref{ferradura} is done in Section \ref{ferraduraprova}.

\begin{remark}
	In \cite{KATOK}, one finds a detailed description of Smale horseshoes for a diffeomorphism and some basic properties. Also, in \cite{Wiggins}, the author provides an elucidative construction of Smale horseshoes.
\end{remark}

Theorem \ref{ferradura} shows us that, if $Z_0$ satisfies $(TC)$ and $(R)$ conditions, then the dynamics given by the crossing orbits is \textbf{chaotic} (see \cite{Wiggins} for more details) and thus $Z_0$ has positive entropy. A direct consequence of Theorem \ref{ferradura} is stated below.

\begin{mtheorem}\label{conseq}
	Let $Z_0\in\Or$ satisfying $(TC)$ and $(R)$ conditions and let $\Lambda$ be the hyperbolic set given by Theorem \ref{ferradura}. The following statements hold.
	\begin{enumerate}
		\item There exists an infinity of closed crossing orbits  $\Gamma$ of $Z_0$, which are of saddle type (i.e. the first return map of $Z_0$ associated to $\Gamma$ has a hyperbolic fixed point of saddle type);
		\item There exists an infinity of non-closed crossing orbits $\Gamma$ of $Z_0$;
		\item There exists a crossing orbit $\Gamma_d$ of $Z_0$ such that $\Gamma_d\cap\Lambda$ is dense in $\Lambda$.
	\end{enumerate}
\end{mtheorem}	

The proof of Theorem \ref{conseq} follows directly from Theorem \ref{ferradura} and Theorem $2.1.4$ of \cite{Wiggins}.

\section{Proof of Theorem \ref{ferradura}}\label{ferraduraprova}

First we discuss about the local structure of the stable $T$-singularity $p_0$. To fix thoughts and without loss of generality, we consider the following assumptions below.

\begin{itemize}
	\item The switching manifold is given by $\s=\{(x,y,z);\ z=0\}$ ($\s=\{z=0\}$ for short) and $p_0=(0,0,0)$.
	\item The sections $\tau^u$ and $\tau^s$ are contained in the planes $\{y=\e\}$ and $\{y=-\e\}$, for some $\e>0$ sufficiently small.
	\item $S_{X_0}\cap V$ and $S_{Y_0}\cap V$ are contained in the lines $x=K_1 y$ and $x=K_2 y$, respectively, for some coefficients $K_1<0$ and $K_2>0$.
	\item The orbits of $Y_0$ in $V_{3D}$ go from $\{x<K_2 y\}$ to $\{x>K_2 y\}$ and the orbits of $X_0$ in $V_{3D}$ goes from $\{x>K_1 y\}$ to $\{x<K_1 y\}$. 
\end{itemize}

Recall that the origin is a stable $T$-singularity of $Z_0$ and a hyperbolic fixed point of saddle type of $\phi_0$. Call $W^{u,s}_{\phi_0}(0,0)$ the local invariant manifolds of $\phi_0$ at $(0,0)$. Without loss of generality, we assume that these invariant manifolds are contained in the union of lines $\{x=K_3 y\}\cup \{x=K_4 y\}$, where $K_3<K_1$ and $K_4>K_2$. 

\begin{remark}
	From Lemma \ref{extlema}, we have that the local first return map $\phi_0$ is globally extended to a first return map $\Phi_0$ defined in an open set of $\s$ and the local invariant manifolds  $W^{u,s}_{\phi_0}(0,0)$ are extended to global invariant manifolds $W^{u,s}_{\Phi_0}$ of $\Phi_0$. In light of this, we study the local invariant manifolds $W^{u,s}_{\phi_0}(0,0)$ having in mind that they exist globally.
\end{remark}

It follows from the orientation of the orbits of $X_0$ and $Y_0$ and the position of $\tau^{u,s}$  that
$$W^{s}_{\phi_0}(0,0)\subset\{x=K_4 y\} \quad and \quad W^{u}_{\phi_0}(0,0)\subset\{x=K_3y\}.$$

\begin{figure}[h!]
	\centering
	\bigskip
	\begin{overpic}[width=10cm]{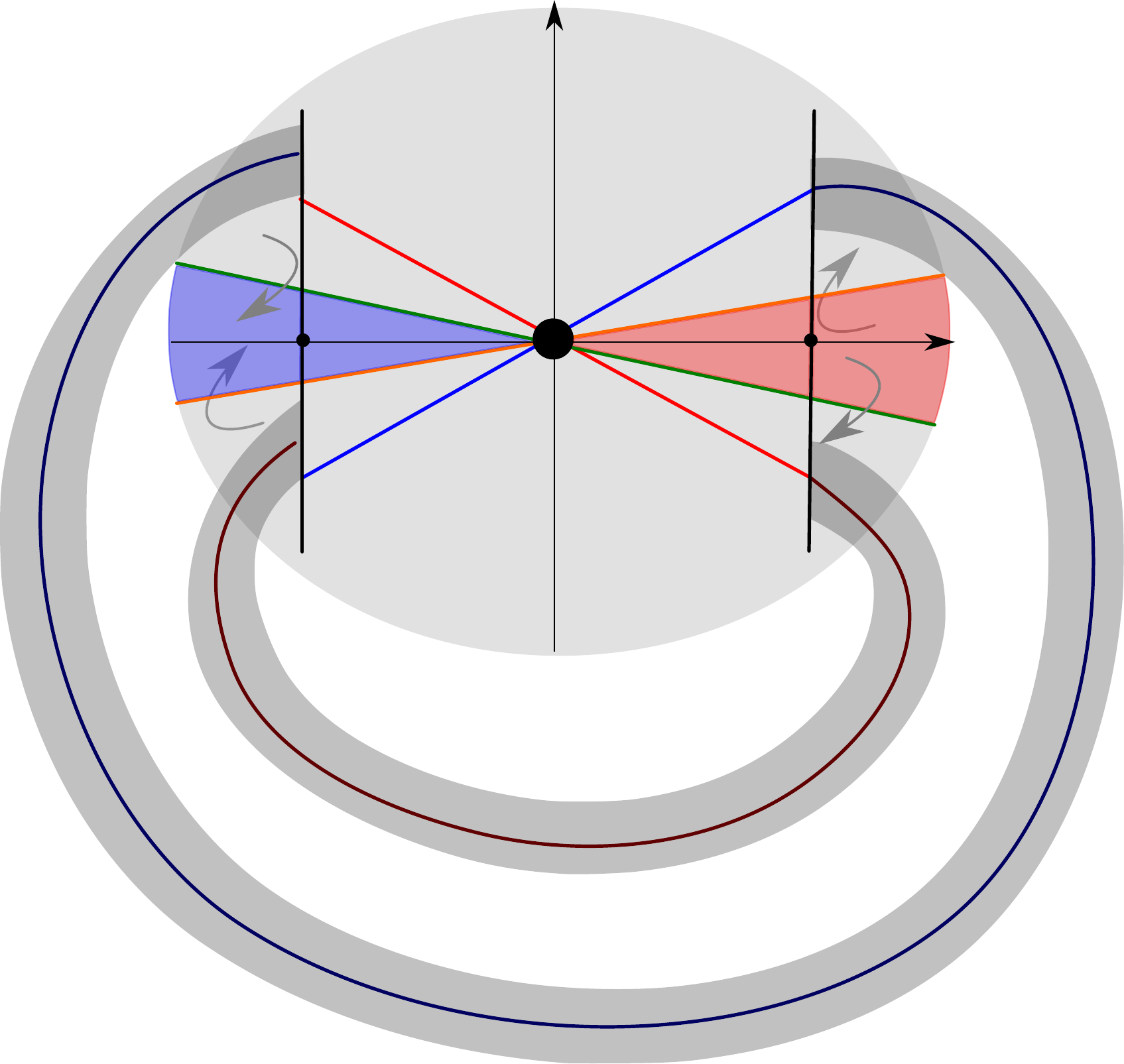}
		\put(86,63){{\footnotesize $y$}}
		\put(52,92){{\footnotesize $x$}}		
		\put(57,75){{\footnotesize $W^s_{\phi_0}$}}
		\put(35,75){{\footnotesize $W^u_{\phi_0}$}}
		\put(30,69){{\footnotesize $S_{X_0}$}}
		\put(30,58){{\footnotesize $S_{Y_0}$}}	
		\put(20,62){{\footnotesize $-\e$}}	
		\put(73,62){{\footnotesize $\e$}}
		\put(28,82){{\footnotesize $\tau^s$}}
		\put(68,82){{\footnotesize $\tau^u$}}	
		\put(90,71){{\footnotesize $\mathcal{R}^s$}}
		\put(81,40){{\footnotesize $\mathcal{R}^u$}}
		\put(15,66){{\footnotesize $\s^{ss}$}}	
		\put(78,65){{\footnotesize $\s^{us}$}}
		\put(5,85){{\footnotesize $\s$}}										
	\end{overpic}
	\bigskip
	\caption{Switching manifold of $Z_0$: Position of the invariant manifolds and tangency sets.}
	\label{posicaofig}	
\end{figure}
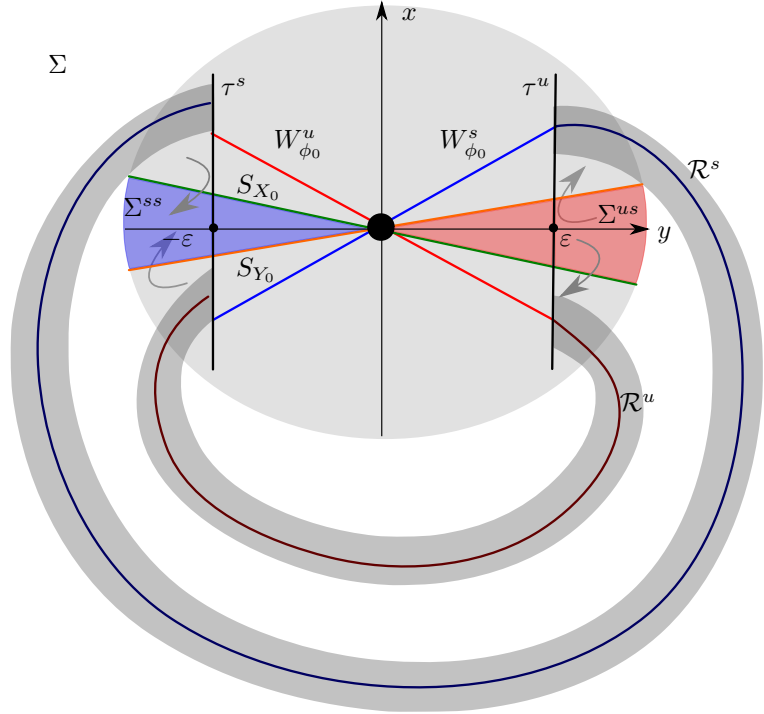

Figure \ref{posicaofig} illustrates the situation considered above. 

Such assumptions imply that, if $p\in \{x<K_3 y\}\cap \{x<K_4 y\}$, then the orbit $(\phi_0^n(p))_{n\in\N}$ of the local first return map $\phi_0:V\rightarrow \s$ represents a crossing orbit of $Z_0$ and its evolution through time coincides with the order generated by $(\phi_0^n(p))_{n\in\N}$.

\begin{remark}
	Notice that, if $x>0$, then the orientation of the crossing orbits through a point $p$ of $W^{u,s}_{\phi_0}(0,0)$ is reverse with respect to the order given by the orbit $(\phi_0^n(p))_{n\in\N}$ of $\phi_0$ through $p$.
\end{remark}

Using the Extension Lemma \ref{extlema}, we obtain the first return map $\Phi_0:W\rightarrow \s$ induced by orbits of $X_0$ and $Y_0$, which extends $\phi_0:V\rightarrow \s$.

From $(TC)$ condition, we have that $X_0$ (resp. $Y_0$) is transverse to $\tau^s$ at points of $(\mathcal{C}^s\cup \widehat{\mathcal{C}^u})\cap \overline{M^+}$ (resp. $(\mathcal{C}^s\cup \widehat{\mathcal{C}^u})\cap \overline{M^-}$). Also, $\pi_2(X_0(p)),\pi_2(Y_0(p))>0$ in such points.

\begin{remark}
	Notice that $S_{X_0}$ and $S_{Y_0}$ must intersect $\tau^s$ at points lying in the interior of the bounded regions of $\{y=-\e\}$ delimited by the circles $\mathcal{C}^s$ and $\widehat{\mathcal{C}^u}$.
\end{remark}

\begin{figure}[H]
	\centering
	\bigskip
	\begin{overpic}[width=7cm]{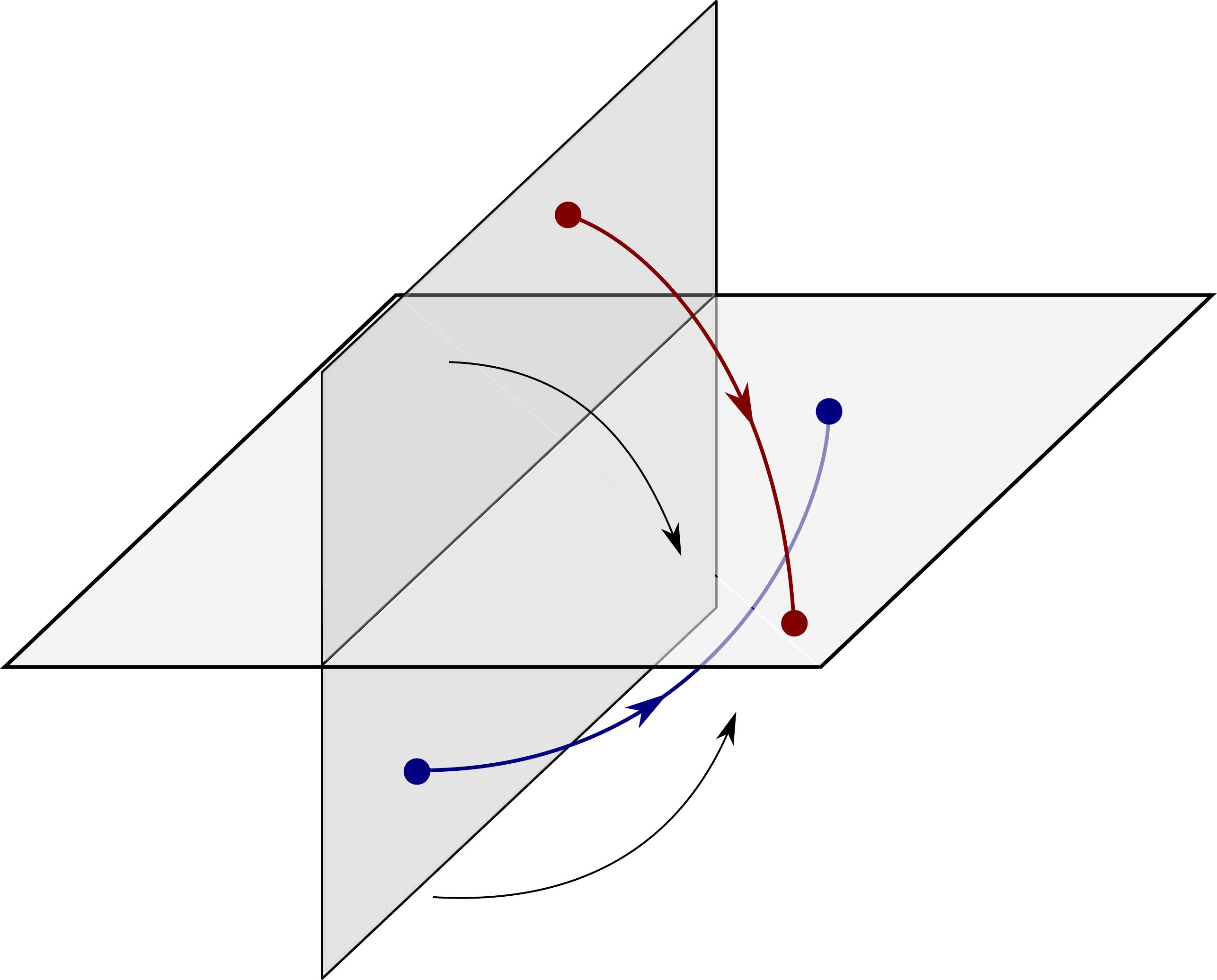}							
		\put(49,65){{\footnotesize $x$}}
		\put(30,19){{\footnotesize $x$}}
		\put(50,28){{\footnotesize $\p_+(x)$}}	
		\put(46,40){{\footnotesize $\p_+$}}
		\put(70,49){{\footnotesize $\p_-(x)$}}	
		\put(55,7){{\footnotesize $\p_-$}}	
		\put(60,77){{\footnotesize $\tau^s$}}
		\put(95,45){{\footnotesize $\s$}}												
	\end{overpic}
	\bigskip
	\caption{Sections $\tau^u$ and $\tau^s$.}
	\label{figg}	
\end{figure}

Recall that near a T-singularity, for any point $p$ of $M^+$ (resp. $M^-$),  there exists a unique orbit of $X_0$ (resp. $Y_0$) contained in $M^+$ (resp. $M^-$) passing through $p$ which connects two points of $\s$. Now, since $(\mathcal{C}^s\cup \widehat{\mathcal{C}^u})\cap \overline{M^+}$  is a compact set, and $\pi_2(X_0(p))>0$ for every $p\in (\mathcal{C}^s\cup \widehat{\mathcal{C}^u})\cap \overline{M^+}$, it follows from the Implicit Function Theorem that there exist an open neighborhood $N^+$ of $(\mathcal{C}^s\cup \widehat{\mathcal{C}^u})\cap \overline{M^+}$ in the plane $\{y=-\e\}$ and a $\Cr$ diffeomorphism $\p_{+}: N^+\cap M^+\rightarrow \s$ such that, for each $x\in N^+\cap M^+$, $x$ and $\p_{+}(x)$ are connected by a unique piece of orbit of $X_0$ contained in $M^+$ which is oriented from $x$ to $\p_{+}(x)$. Analogously, we obtain a $\Cr$ diffeomorphism $\p_{-}: N^-\cap M^-\rightarrow \s$ defined in a neighborhood $N^-$ of $(\mathcal{C}^s\cup \widehat{\mathcal{C}^u})\cap \overline{M^-}$ in the plane $\{y=-\e\}$, such that, for every $x\in N^-\cap M^-$, there exists a unique piece of orbit of $Y_0$ contained in $M^-$ connecting $x$ and $\p_-(x)$ oriented from $x$ to $\p_-(x)$. See Figure \ref{figg}.

Let $N=(N^+\cap M^+)\cup(N^-\cap M^-)$ and $\p:N\rightarrow \s$ be the $\Cr$ diffeomorphism defined by
$$\p(p)=\left\{\begin{array}{l}
\p_{+}(p),\textrm{ if }p\in N^+\cap M^+,\\
\p_{-}(p),\textrm{ if }p\in N^-\cap M^-.
\end{array}
\right.$$

\begin{figure}[h!]
	\centering
	\bigskip
	\begin{overpic}[width=8cm]{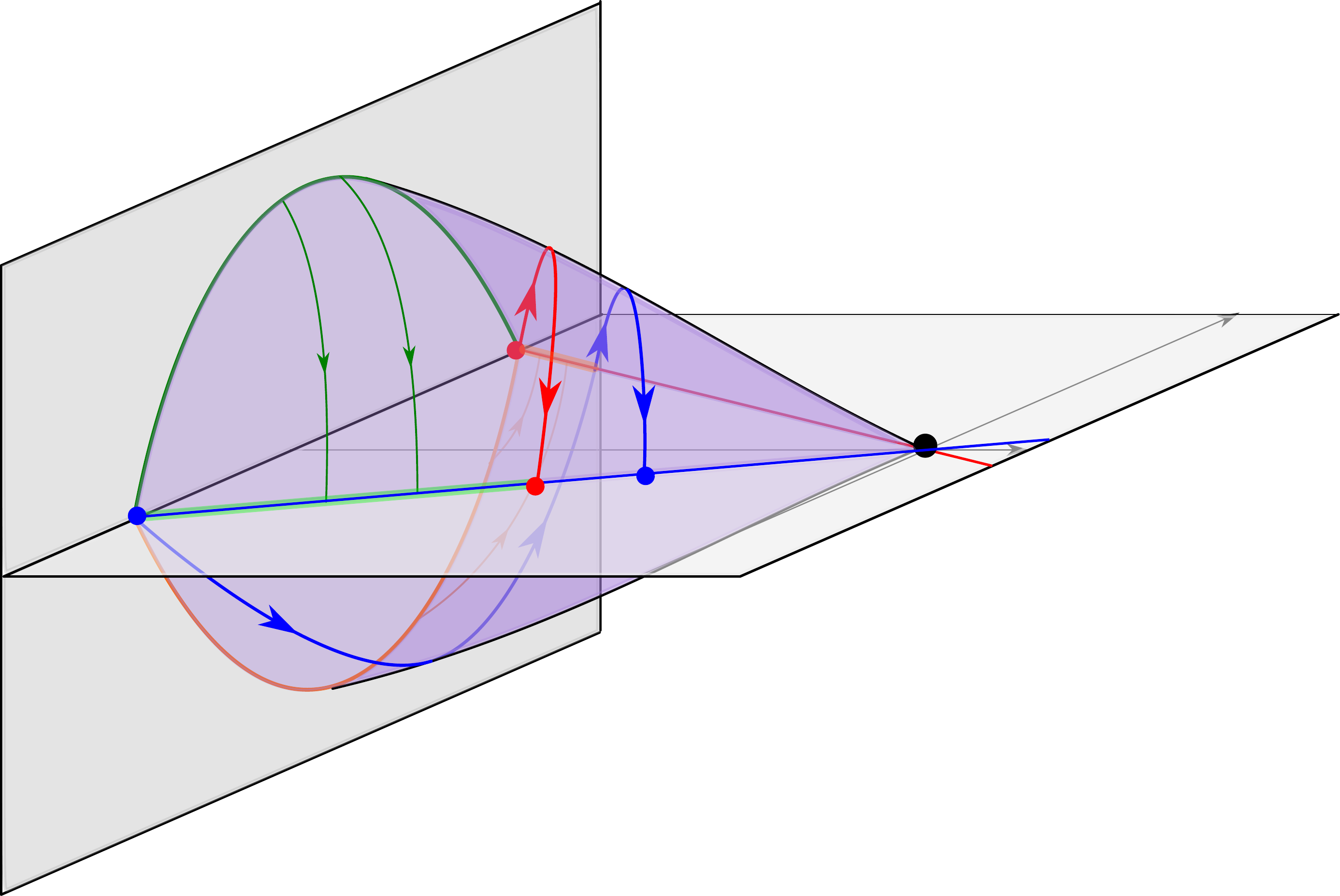}							
		\put(67,38){{\footnotesize $p_0$}}		
		\put(3,30){{\scriptsize $p^s_{\mathrm{loc}}$}}
		\put(31,42){{\scriptsize $p^u_{\mathrm{loc}}$}}				
		\put(31,54){{\footnotesize $\mathcal{C}^s$}}
		\put(47,60){{\footnotesize $\tau^s$}}
		\put(92,45){{\footnotesize $x$}}
		\put(77,31){{\footnotesize $y$}}
		\put(95,37){{\footnotesize $\s$}}	
		\put(-15,23){{\footnotesize $y=-\e$}}
		\put(45,20){{\footnotesize $y=0$}}	
		\put(33,27){{\tiny $\Phi_{X_0}(p^u_{\mathrm{loc}})$}}	
		\put(49,28){{\tiny $\Phi_{0}(p^s_{\mathrm{loc}})$}}											
	\end{overpic}
	\bigskip
	\caption{Evolution of the flow of $Z_0$ through the points $p^{u}_{\mathrm{loc}}$ and $p^{s}_{\mathrm{loc}}$.}
	\label{fig}	
\end{figure}
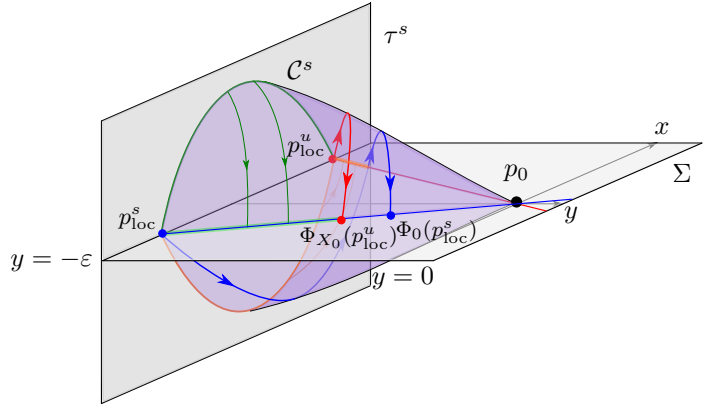

Now, for $\star=u,s$, let $p^{\star}_{\mathrm{loc}}$ be the unique point of $\mathcal{C}^s$ contained in the local invariant manifold $W_{\phi_0}^{\star}(p_0)$. Recall that, there exists a crossing orbit of $Z_0$ from $p^s_{\mathrm{loc}}$ to $\Phi_0(p^s_{\mathrm{loc}})\in W_{\phi_0}^s(p_0)$ and $\pi_2(\Phi_0(p^s_{\mathrm{loc}}))>-\e$. Also we have that $-\e< \pi_2(\Phi_{X_0}(p^u_{\mathrm{loc}}))<\pi_2(\Phi_0(p^s_{\mathrm{loc}}))$. See Figure \ref{fig}.

From the definition of $\Phi_0$, it follows that:

\begin{itemize}
	\item if $p\in \mathcal{C}^s\cap M^+$, then $\p(p)\in W_{\phi_0}^s(p_0)$;
	\item if $p\in \mathcal{C}^s\cap M^-$, then $\p(p)\in W_{\phi_0}^u(p_0)$. 
\end{itemize}

Recall that, for each $p\in W^u_{\phi_0}(0,0)\cap \{x<0\}$, $(\phi_0^n(p))_{n\in\N}$ represents a crossing orbit of $Z_0$ which is oriented in the order given by the iterations of  $p$ through $\phi_0$. Since $\mathcal{R}^u$ extends $W^u_{\phi_0}(0,0)\cap \{x<0\}$, the same property holds for $p\in\mathcal{R}^u$. 

\begin{figure}[h!]
	\centering
	\bigskip
	\begin{overpic}[width=12cm]{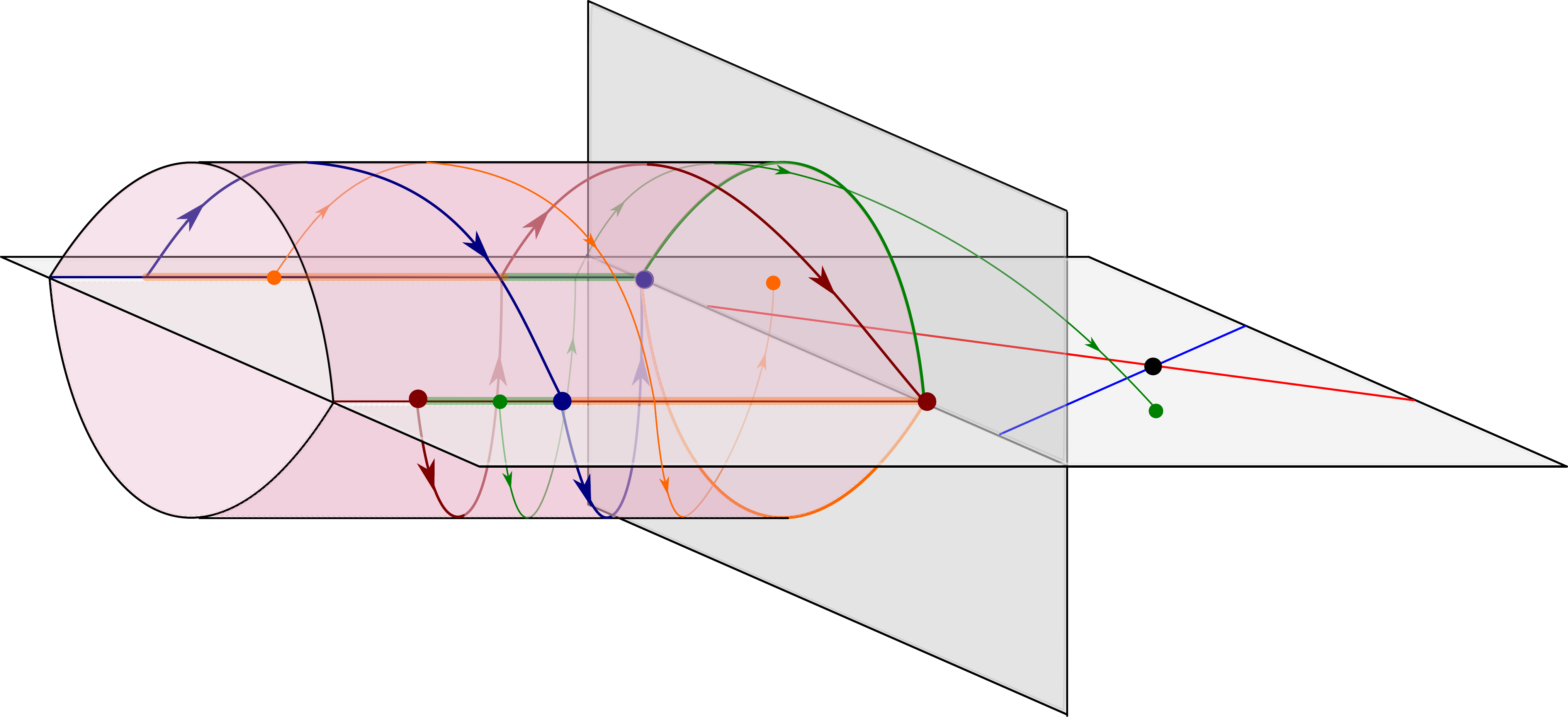}							
		\put(72,24){{\footnotesize $p_0$}}		
		\put(79,26){{\footnotesize $W^s_{\phi_0}$}}		
		\put(92,20){{\footnotesize $W^u_{\phi_0}$}}		
		\put(42.5,28){{\scriptsize $p^s_{\mathcal{R}}$}}
		\put(61,20){{\scriptsize $p^u_{\mathcal{R}}$}}
		\put(46,37){{\footnotesize $\widehat{\mathcal{C}^u}$}}
		\put(34,42){{\footnotesize $\tau^s$}}
		\put(24,29){{\footnotesize $\mathcal{R}^s$}}
		\put(44,21){{\footnotesize $\mathcal{R}^u$}}
		\put(95,12){{\footnotesize $\s$}}	
		\put(33,18){{\tiny $\Phi_{Y_0}(p^s_{\mathcal{R}})$}}	
		\put(23,22){{\tiny $\Phi_{0}^{-1}(p^u_{\mathcal{R}})$}}											
	\end{overpic}
	\bigskip
	\caption{Evolution of the flow of $Z_0$ through the points $p^{u}_{\mathcal{R}}$ and $p^{s}_{\mathcal{R}}$.}
	\label{fig222}	
\end{figure}
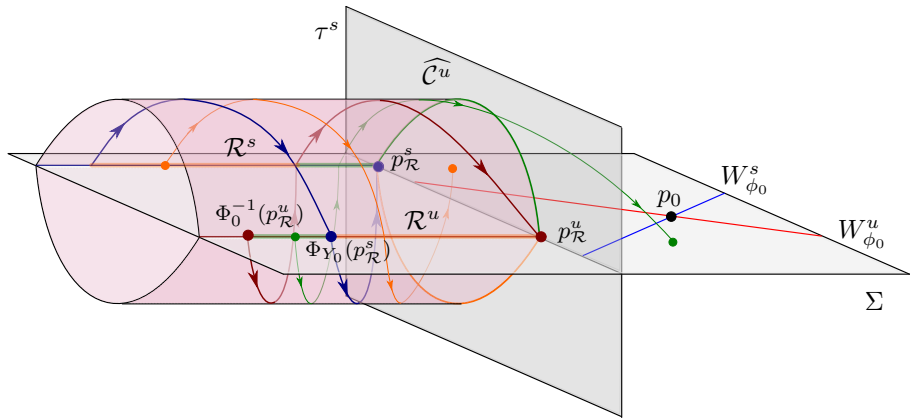

For $\star=u,s$, let $p^{\star}_{\mathcal{R}}$ be the unique point of  $\widehat{\mathcal{C}^u}$ contained in the curve $\mathcal{R}^{\star}$. There exists a crossing orbit of $Z_0$ from $\Phi_0^{-1}(p^u_{\mathcal{R}})\in\mathcal{R}^u$ to $p^u_{\mathcal{R}}$, and $\pi_2(\Phi_0^{-1}(p^u_{\mathcal{R}}))<-\e$. Also we have that $\pi_2(\Phi_0^{-1}(p^u_{\mathcal{R}}))<\pi_2(\Phi_{Y_0}(p^{s}_{\mathcal{R}}))<-\e$. See Figure \ref{fig222}. Therefore,
\begin{itemize}
	\item if $p\in \widehat{\mathcal{C}^u}\cap M^+$, then $\p(p)=\Phi_0(\widetilde{p})$, for some $\widetilde{p}\in W_{\Phi_0}^u$;
	\item if $p\in \widehat{\mathcal{C}^u}\cap M^-$, then $\p(p)=\Phi_0^{-1}(\widetilde{p})$, for some $\widetilde{p}\in W_{\Phi_0}^s$. 
\end{itemize}

Finally, from $(R)$ condition, we have that $\mathcal{C}^s\cap \widehat{\mathcal{C}^u}=\{q_1,q_2\}$, where $q_i\notin\s$ and $\mathcal{C}^s\pitchfork \widehat{\mathcal{C}^u}$ at $q_i$, $i=1,2$. Without loss of generality, assume that $q_1\in M^+$. Since $\p$ is a diffeomorphism, $\p(\widehat{\mathcal{C}^u}\cap M^+)\subset W_{\Phi_0}^u$ and $\p(\mathcal{C}^s\cap M^+)\subset W_{\Phi_0}^s$, it follows that the invariant manifolds $W_{\Phi_0}^s$ and $W_{\Phi_0}^u$ intersect transversally at the point $\widehat{q_1}=\p(q_1)$. Analogously, we have that $W_{\Phi_0}^s$ and $W_{\Phi_0}^u$ intersect transversally at the point $\widehat{q_2}=\p(q_2)$. Also, $(\Phi_0^n(\widehat{q_1}))_{n\in\N}$ and $(\Phi_0^n(\widehat{q_2}))_{n\in\N}$ define two distinct orbits of $\Phi_0$.

	Therefore,  the existence of the Smale horseshoe for $\Phi_0$ follows straightly from Theorem $6.5.5$ from \cite{KATOK}, which is stated below.
	
	\begin{theorem}[Theorem $6.5.5$ of \cite{KATOK}]
		Let $M$ be a smooth manifold, $U\subset M$ open, $f: U\rightarrow M$ an embedding, and $p\in U$ a hyperbolic fixed point with a transverse homoclinic point $q$. Then in an arbitrarily small neighborhood of $p$ there exists a horseshoe for some iterate of $f$. Furthermore the hyperbolic invariant set in this horseshoe contains an iterate of $q$. 
\end{theorem} 

This proves statements $(1)$ and $(2)$ of Theorem \ref{ferradura}. Now, we will prove statement $(3)$ of Theorem \ref{ferradura}. In order to do this, we will show that the hyperbolic invariant set $\Lambda$ of the horseshoe $\Delta$ is contained in $\s^c$, and as a consequence of the construction, each orbit of $Z_0$ which intersects $\Lambda$ is a crossing orbit.

We use topological arguments intrinsic to the structure of the problem to construct a horseshoe $\Delta$ with the points $\widehat{q_1}$ and $\widehat{q_2}$, such as the set $\Lambda$ associated to $\Delta$. Consider the following steps:

\textbf{Step $1$:} First we use $(TC)$ and $(R)$ conditions to deduce the pattern associated to the horseshoe $\Delta$.

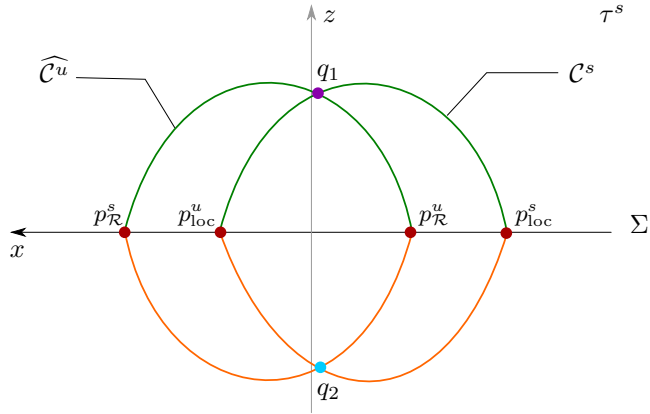
\begin{figure}[h!]
	\centering
	\bigskip
	\begin{overpic}[width=8cm]{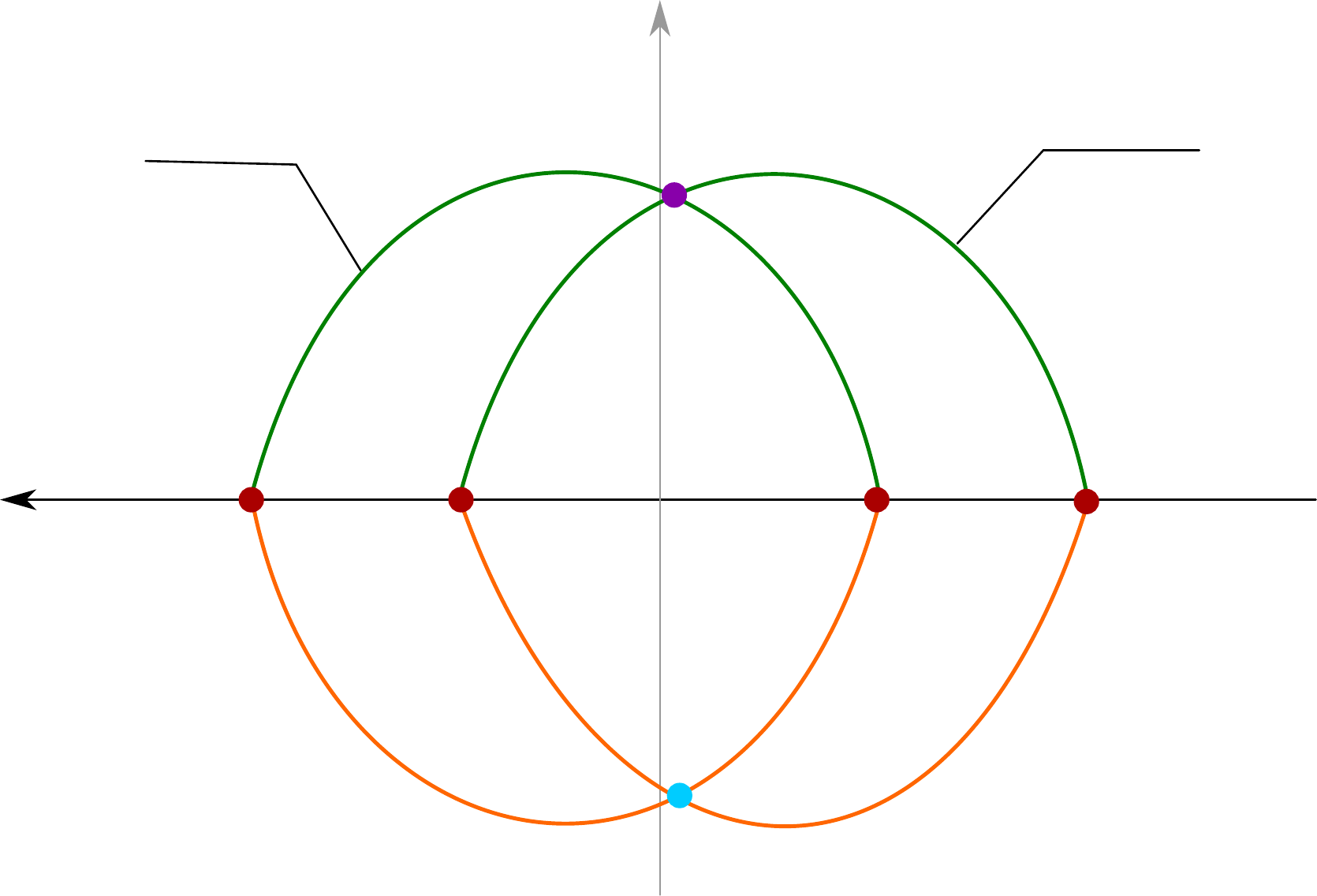}							
		\put(5,55){{\footnotesize $\widehat{\mathcal{C}^u}$}}		
		\put(93,55){{\footnotesize $\mathcal{C}^s$}}		
		\put(51,56){{\footnotesize $q_1$}}		
		\put(51,3){{\footnotesize $q_2$}}			
		\put(14,32){{\scriptsize $p^s_{\mathcal{R}}$}}
		\put(68,32){{\scriptsize $p^u_{\mathcal{R}}$}}
		\put(28,32){{\scriptsize $p^u_{\mathrm{loc}}$}}			
		\put(84,32){{\scriptsize $p^s_{\mathrm{loc}}$}}					
		\put(0,26){{\footnotesize $x$}}
		\put(52,65){{\footnotesize $z$}}
		\put(98,65){{\footnotesize $\tau^s$}}		
		\put(103,30){{\footnotesize $\s$}}										
	\end{overpic}
	\bigskip
	\caption{Configuration of the circles $\mathcal{C}^s$ and $\widehat{\mathcal{C}^u}$ in the section $\tau^s$.}
	\label{fig1}	
\end{figure}

Without loss of generality, assume that the points $p^\star_{\mathcal{R}}$ and $p^{\star}_{\mathrm{loc}}$, $\star=u,s$ given by $\widehat{\mathcal{C}^u}\cap \s$ and $\mathcal{C}^s\cap\s$, respectively, satisfy the following order in the line $\tau^s\cap\s=\{z=0,y=-\e\}$:
$$p^s_{\mathcal{R}}>p^{u}_{\mathrm{loc}}>0>p^u_{\mathcal{R}}>p^{s}_{\mathrm{loc}}.$$
Also, assume that $q_1\in M^+$ and $q_2\in M^-$. Hence we have the following situation illustrated in Figure \ref{fig1}.

\begin{figure}[h!]
	\centering
	\bigskip
	\begin{overpic}[width=7cm]{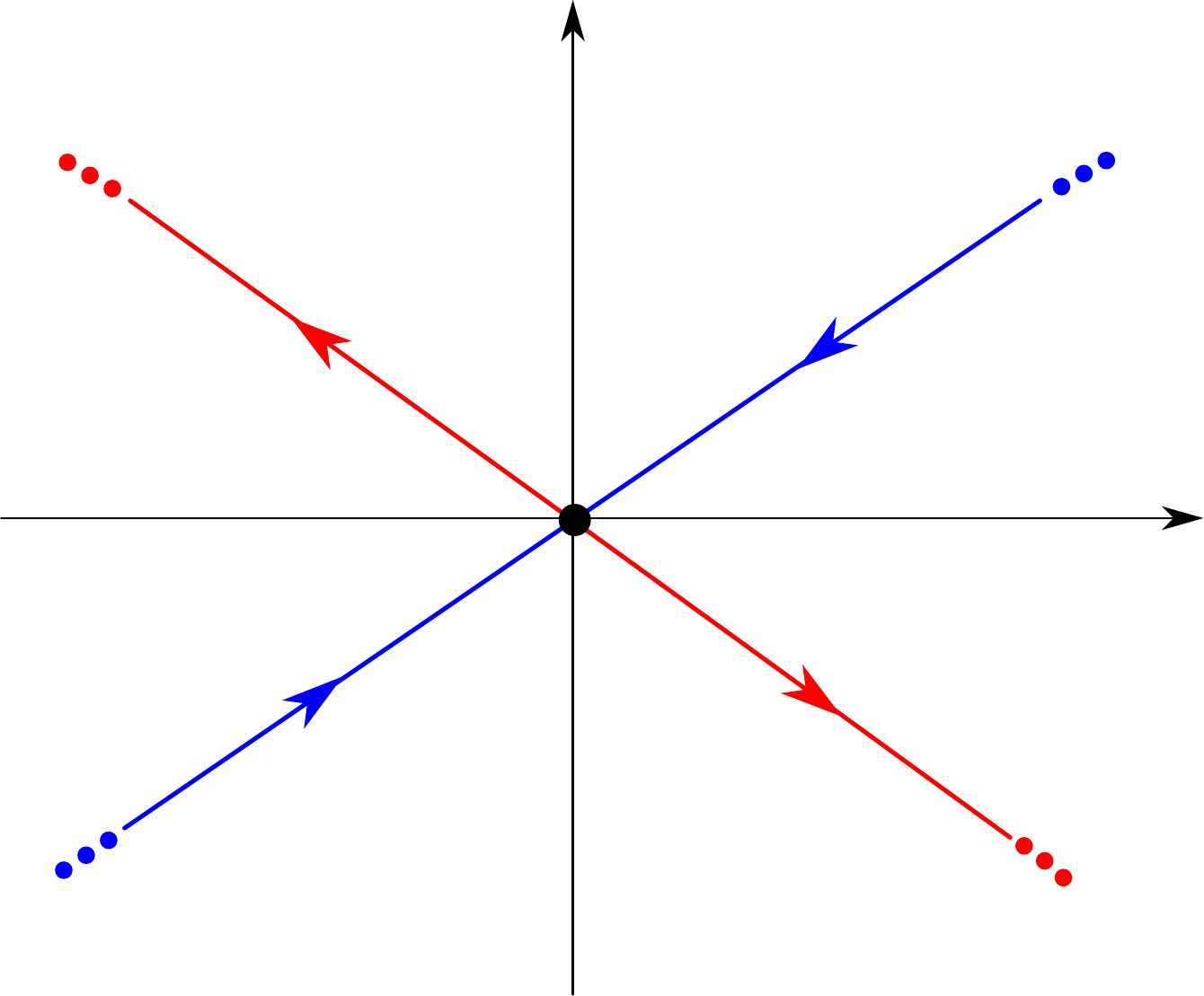}							
		\put(97,42){{\footnotesize $y$}}
		\put(50,80){{\footnotesize $x$}}
		\put(72,65){{\footnotesize $W^{s,+}_{\Phi_0}$}}	
		\put(15,65){{\footnotesize $W^{u,+}_{\Phi_0}$}}		
		\put(15,12){{\footnotesize $W^{s,-}_{\Phi_0}$}}	
		\put(67,12){{\footnotesize $W^{u,-}_{\Phi_0}$}}					
		\put(95,80){{\footnotesize $\s$}}										
	\end{overpic}
	\bigskip
	\caption{Invariant manifolds $W^{\star,\pm}_{\Phi_0},\ \star=u,s$ in $\s$.}
	\label{fig2}	
\end{figure}
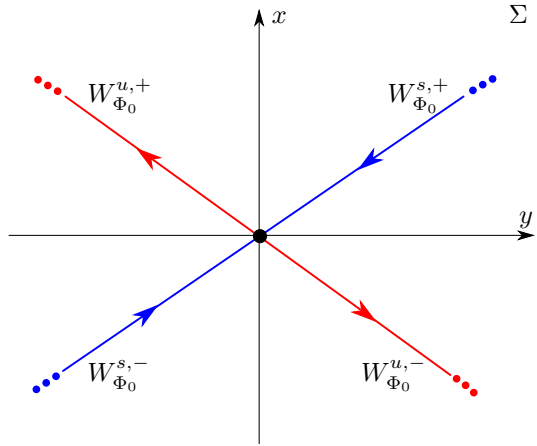

In this case, we have that $\widehat{q_1}\in W^s_{\phi_0}(0,0)$ and $\widehat{q_2}\in W^u_{\phi_0}(0,0)$. In order to clarify the notation, we distinguish $W^\star_{\Phi_0}$ between $W^{\star,+}_{\Phi_0}$ and $W^{\star,-}_{\Phi_0}$ in such a way that $W^{\star,\pm}_{\Phi_0}$ corresponds to the branch of $W^\star_{\Phi_0}$  which contains (the local invariant manifold) $W^\star_{\phi_0}(0,0)\cap\{\pm x>0\}$. See Figure \ref{fig2}.

As we have seen before, using $\p$ and $\Phi_{X_0}$, we can see that the circle $\mathcal{C}^s$ can be brought to a fundamental domain $F^s$ of $\Phi_0$ for the invariant manifold  $W^{s,-}_{\Phi_0}$. Analogously, we have that the circle $\widehat{\mathcal{C}^u}$ is a fundamental domain $F^u$ of $\Phi_0$ for $W^{u,-}_{\Phi_0}$. See Figures \ref{fig} and \ref{fig222}.

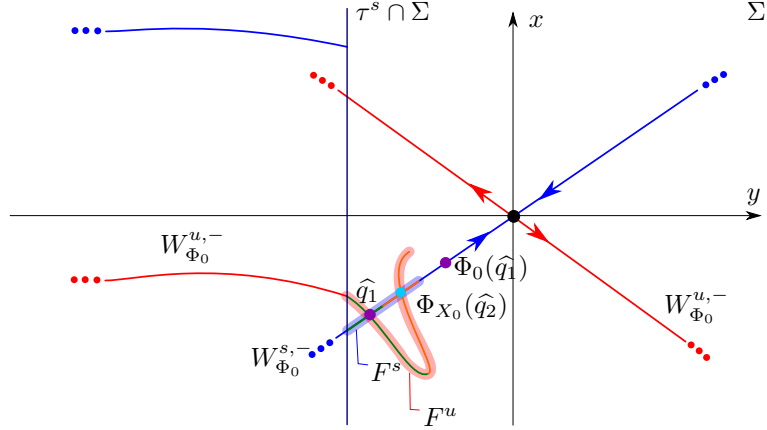
\begin{figure}[h!]
	\centering
	\bigskip
	\begin{overpic}[width=10cm]{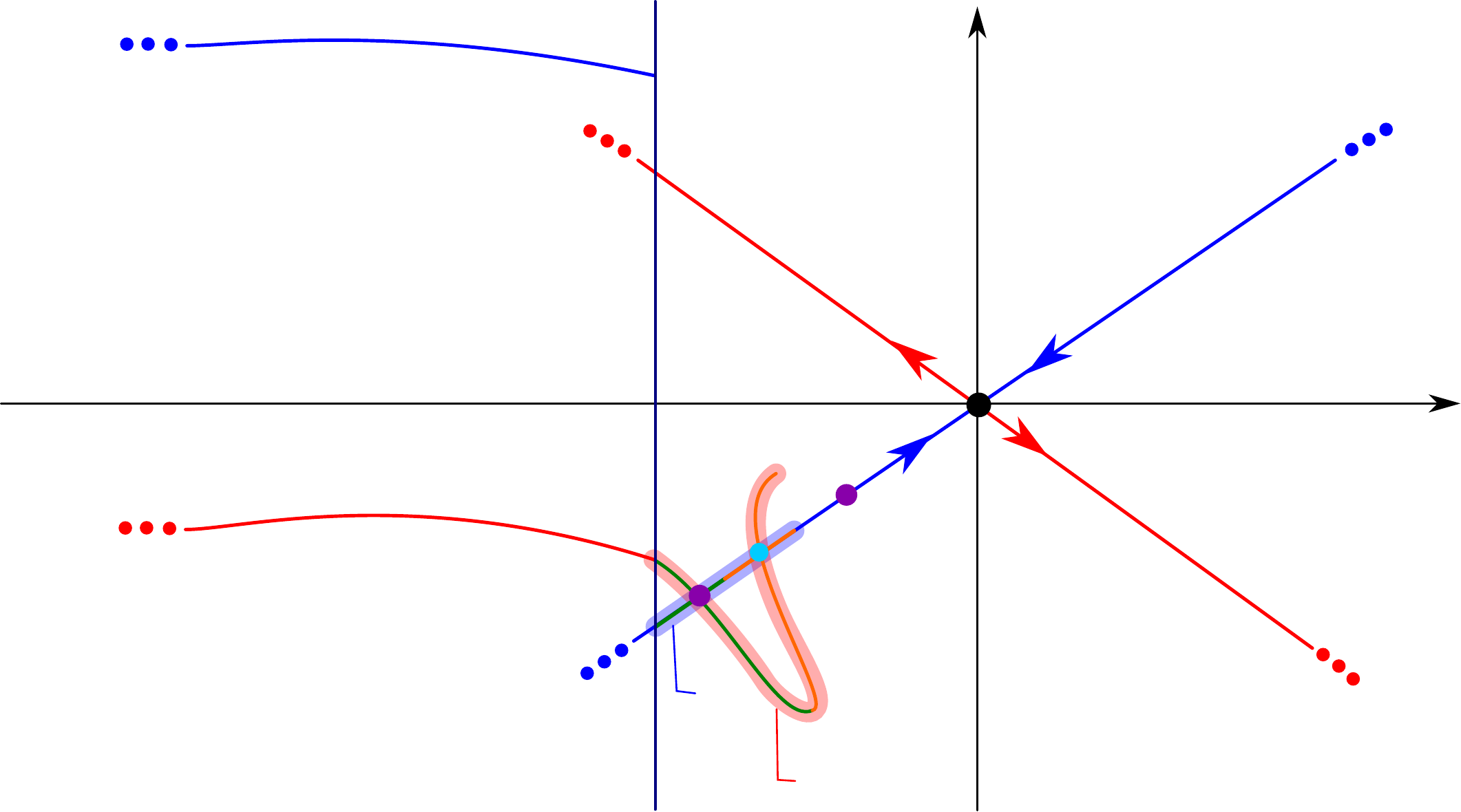}							
		\put(48,6){{\footnotesize $F^s$}}
		\put(55,0){{\footnotesize $F^u$}}
		\put(46,17){{\footnotesize $\widehat{q_1}$}}	
		\put(54,15){{\footnotesize $\Phi_{X_0}(\widehat{q_2})$}}		
		\put(59,20){{\footnotesize $\Phi_{0}(\widehat{q_1})$}}	
		\put(87,15){{\footnotesize $W^{u,-}_{\Phi_0}$}}				
	\put(20,23){{\footnotesize $W^{u,-}_{\Phi_0}$}}		
	\put(32,8){{\footnotesize $W^{s,-}_{\Phi_0}$}}								
		\put(69,53){{\footnotesize $x$}}	
		\put(98,30){{\footnotesize $y$}}	
		\put(98,54){{\footnotesize $\s$}}	
			\put(46,54){{\footnotesize $\tau^s\cap\s$}}														
	\end{overpic}
	\bigskip
	\caption{Fundamental domains $F^u$ and $F^s$ of $W^{u,-}_{\Phi_0}$ and $W^{s,-}_{\Phi_0}$ given by the projection of the circles $\widehat{\mathcal{C}^u}$ and $\mathcal{C}^s$ into $\s$ through crossing orbits of $Z_0$, respectively.}
	\label{fig3}	
\end{figure}

From $(TC)$ and $(R)$ conditions and the considerations above, we have that these fundamental domains of $W^{u,-}_{\Phi_0}$ and $W^{s,-}_{\Phi_0}$ intersect only at the two points $\widehat{q_1}$ and $\phi_{X_0}(\widehat{q_2})$ and present the configuration illustrated in the Figure \ref{fig3}

Now, since the considered fundamental domains are contained in the neighborhood $V$ of the origin for which the local involutions associated with the T-singularity are defined, then we can use the local behavior of the orbits of $X_0$ and $Y_0$ near the origin to obtain the following pattern of the intersections between the invariant manifolds $W^{u,-}_{\Phi_0}$ and $W^{s,-}_{\Phi_0}$. 

\begin{figure}[h!]
	\centering
	\bigskip
	\begin{overpic}[width=7cm]{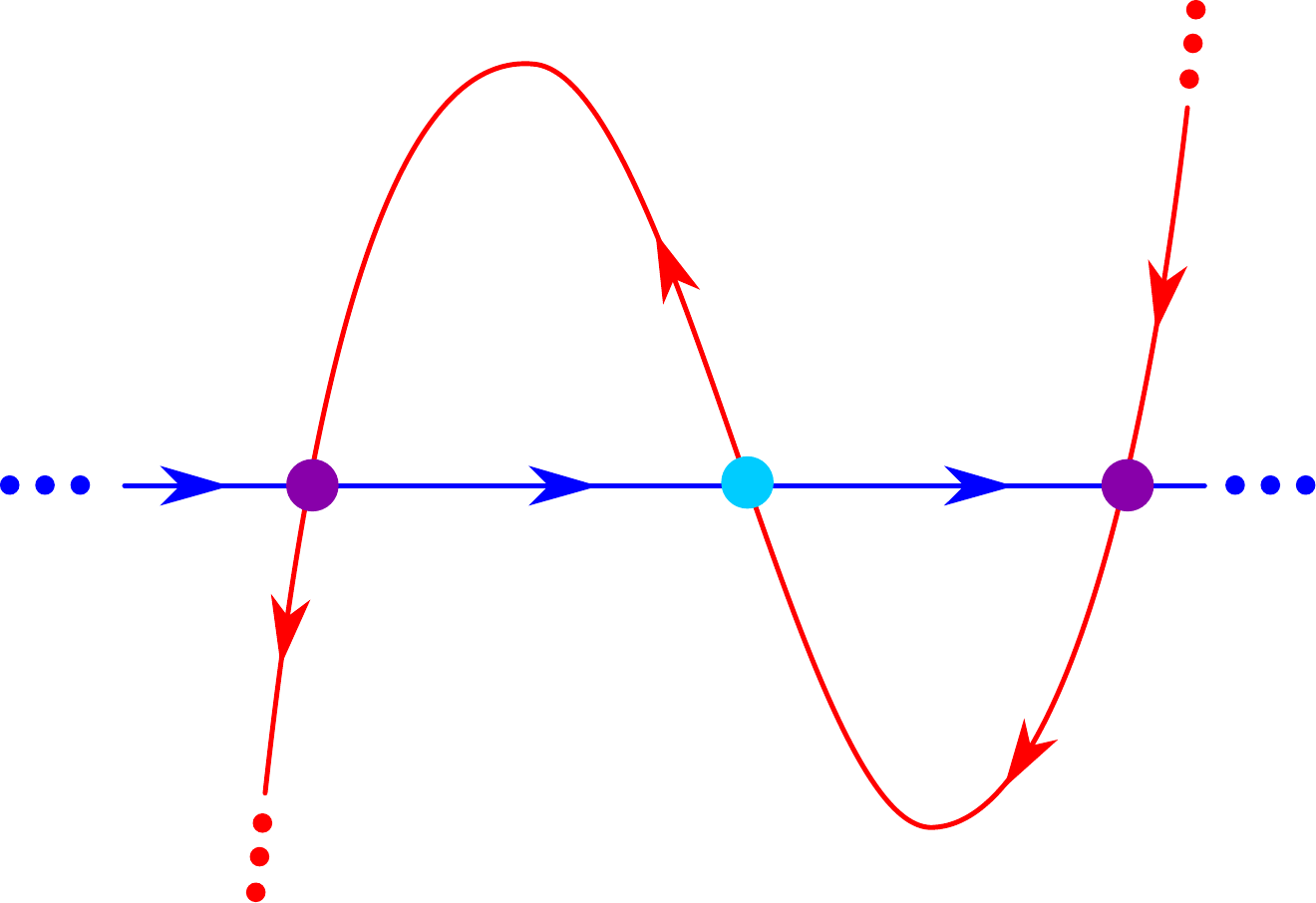}							
		\put(32,24){{\footnotesize $W^{s,-}_{\Phi_0}$}}	
		\put(52,55){{\footnotesize $W^{u,-}_{\Phi_0}$}}					
		\put(17,35){{\footnotesize $\widehat{q_1}$}}	
		\put(72,35){{\footnotesize $\Phi_0(\widehat{q_1})$}}	
		\put(37,35){{\footnotesize $\Phi_{X_0}(\widehat{q_2})$}}														
	\end{overpic}
	\bigskip
	\caption{Pattern of the intersections between the invariant manifolds $W^{u,-}_{\Phi_0}$ and $W^{s,-}_{\Phi_0}$. }
	\label{fig4}	
\end{figure}
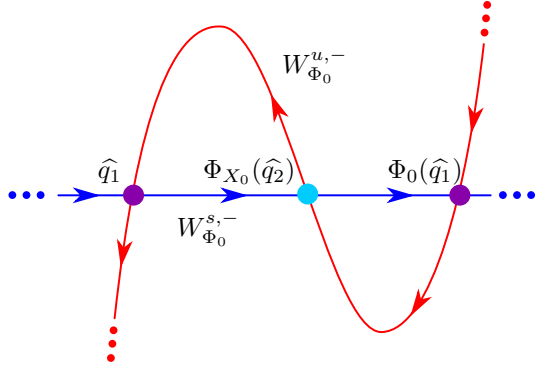

 \textbf{Pattern:} $W^{u,-}_{\Phi_0}$ intersects the segment of $W^{s,-}_{\Phi_0}$ between $\widehat{q_1}$ and $\Phi_0(\widehat{q_1})$ only at the point  $\Phi_0(\widehat{q_2})$. Furthermore, $W^{u,-}_{\Phi_0}\pitchfork W^{s,-}_{\Phi_0}$ at $\phi_{X_0}(\widehat{q_2})$. See Figure \ref{fig4}.

\begin{remark}
	The reversibility of the problem can be used to obtain similar properties for the invariant manifolds $W^{u,+}_{\Phi_0}$ and $W^{s,+}_{\Phi_0}$. In fact, all the constructions below can be done for these invariant manifolds, nevertheless it does not generate new dynamic features for the system, and for this reason we consider only  the invariant manifolds $W^{u,-}_{\Phi_0}$ and $W^{s,-}_{\Phi_0}$.
\end{remark}	

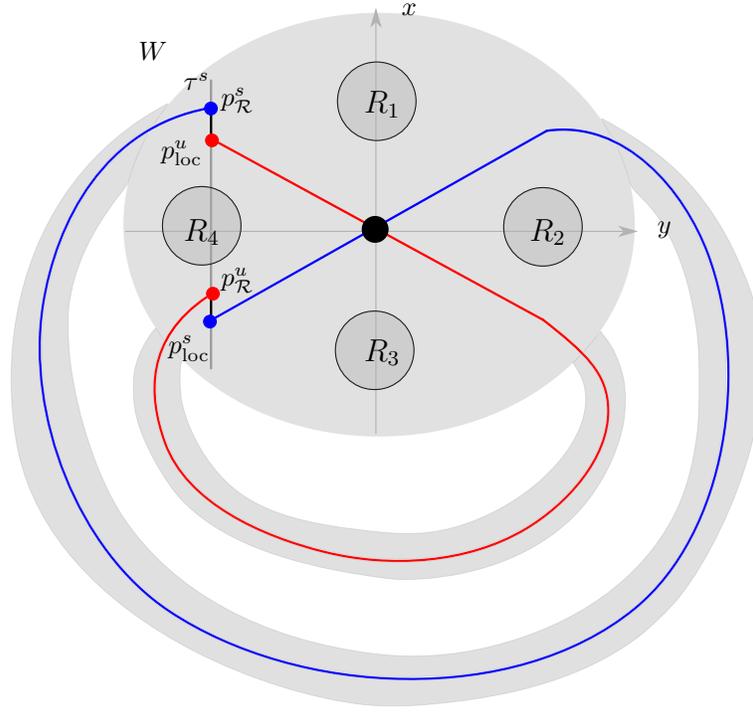
\begin{figure}[h!]
	\centering
	\bigskip
	\begin{overpic}[width=10cm]{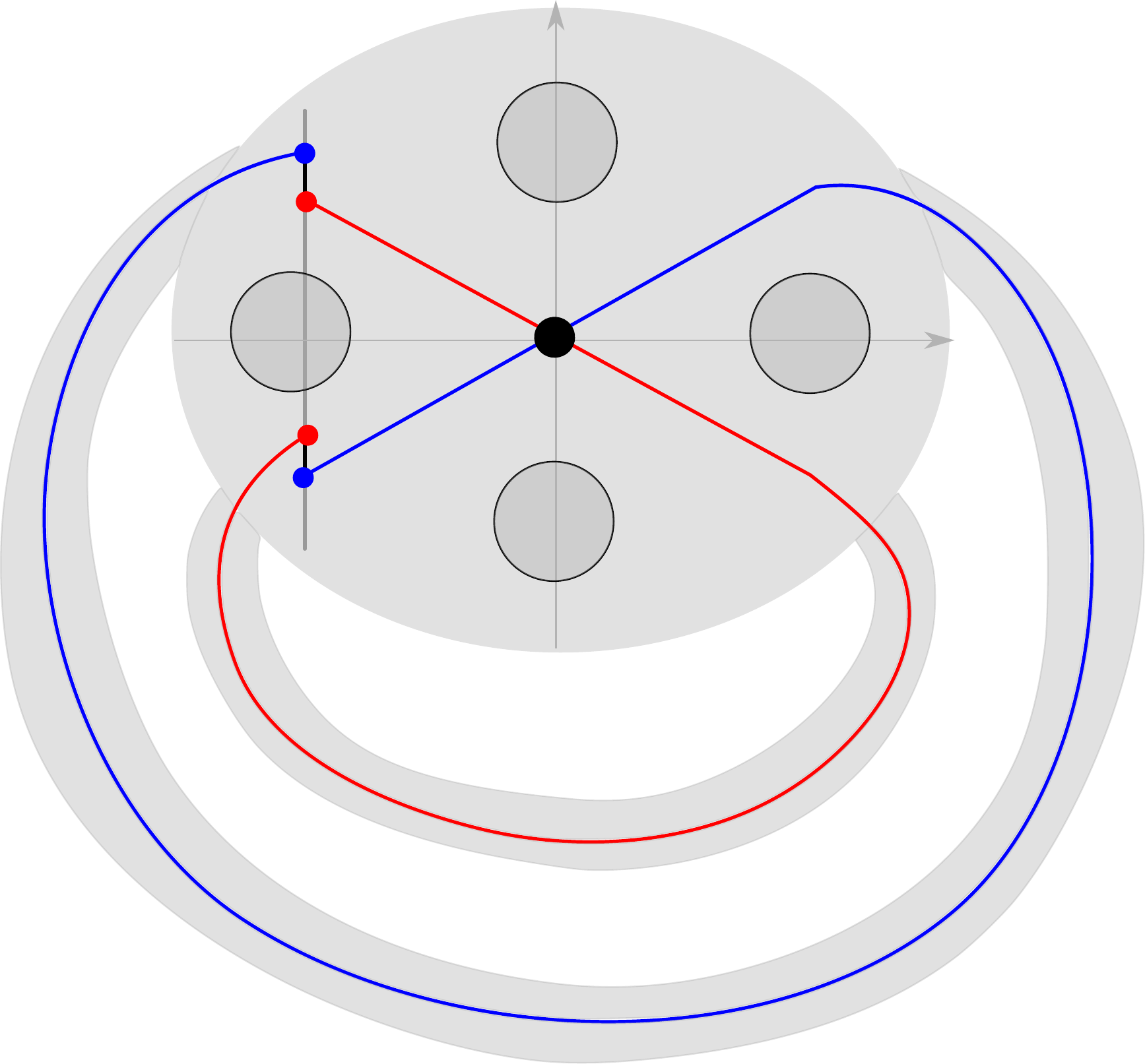}
		\put(86,63){{\footnotesize $y$}}
		\put(52,92){{\footnotesize $x$}}		
		\put(23,82){{\footnotesize $\tau^s$}}	
			\put(17,86){{\footnotesize $W$}}	
	\put(28,80){{\footnotesize $p^s_{\mathcal{R}}$}}
	\put(20,73){{\footnotesize $p^u_{\mathrm{loc}}$}}	
	\put(28,56){{\footnotesize $p^u_{\mathcal{R}}$}}
\put(21,47){{\footnotesize $p^s_{\mathrm{loc}}$}}
\put(47,79){$R_1$}			
\put(47,46){$R_3$}
\put(69,62){$R_2$}	
\put(23,62){$R_4$}													
	\end{overpic}
	\bigskip
	\caption{Regions $R_i$, $i=1,2,3,4$, of the neighborhood $W$.}
	\label{regionsfig}	
\end{figure}

\textbf{Step $2$:} Consider the orbits connecting $(0,0)$ to $p^u_{\mathrm{loc}},p^s_{\mathrm{loc}}, p^u_{\mathcal{R}}, p^s_{\mathcal{R}}$, and the segments of the section $\tau^s$ connecting $p^u_{\mathcal{R}}$ to $p^s_{\mathrm{loc}}$ and $p^s_{\mathcal{R}}$ to $p^u_{\mathrm{loc}}$, and notice that they split the neighborhood $W$ into four regions $R_i$, $i=1,\cdots, 4$, as it is illustrated in Figure \ref{regionsfig}.

Now, we consider a square $Q$ such that:
\begin{enumerate}[$i)$]
	\item $Q\subset R_3$
	\item There exists a side of $Q$, say it $L$, which is contained in $W^{s,-}_{\Phi_0}$;
	\item $\partial L=\{l_1,l_2\}$, where $l_1$ is a point of  $W^{s,-}_{\Phi_0}$ between $p^s_{\mathrm{loc}}$ and $\widehat{q_1}$, and $l_2$ is a point of  $W^{s,-}_{\Phi_0}$ between $\phi_{X_0}(\widehat{q_2})$ and $\Phi_0(\widehat{q_1})$;
	\item There exists a segment $A_1$ of $W^{u,-}_{\Phi_0}$ connecting $\widehat{q_1}$ and a point of the interior of the opposite side of $L$ which is contained in the interior of $Q$.
	\item There exists a segment $A_2$ of $W^{u,-}_{\Phi_0}$ connecting $\phi_{X_0}(\widehat{q_2})$ and a point of the interior of the opposite side of $L$ which is contained in the interior of $Q$.	
	\item $A_1$ and $A_2$ are transverse to $\partial Q$.
	\item $Q\subset\s^c$.
\end{enumerate}

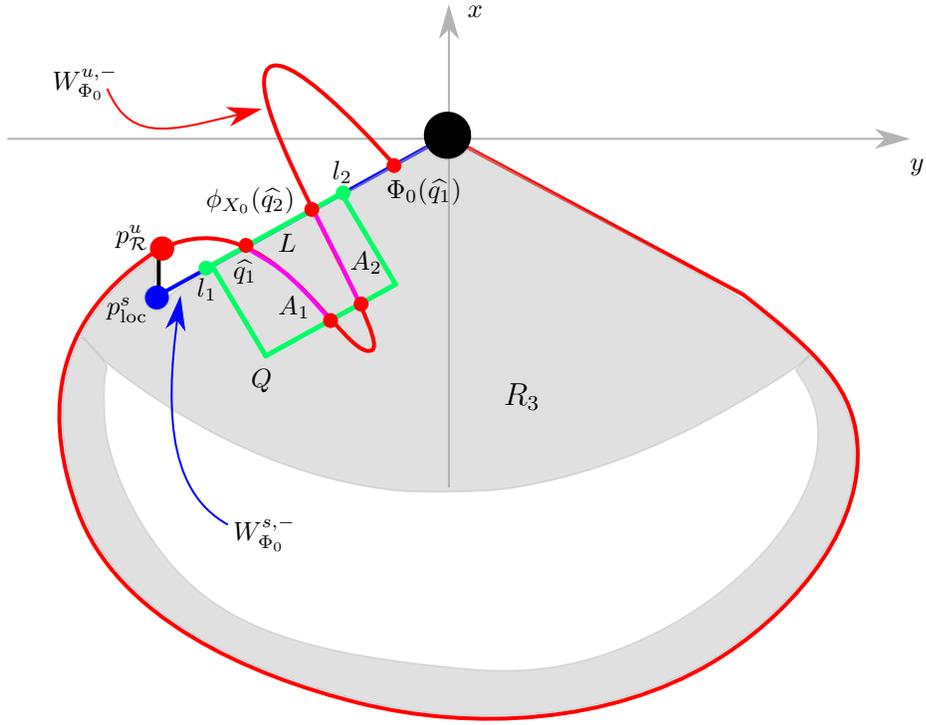
\begin{figure}[h!]
	\centering
	\bigskip
	\begin{overpic}[width=12cm]{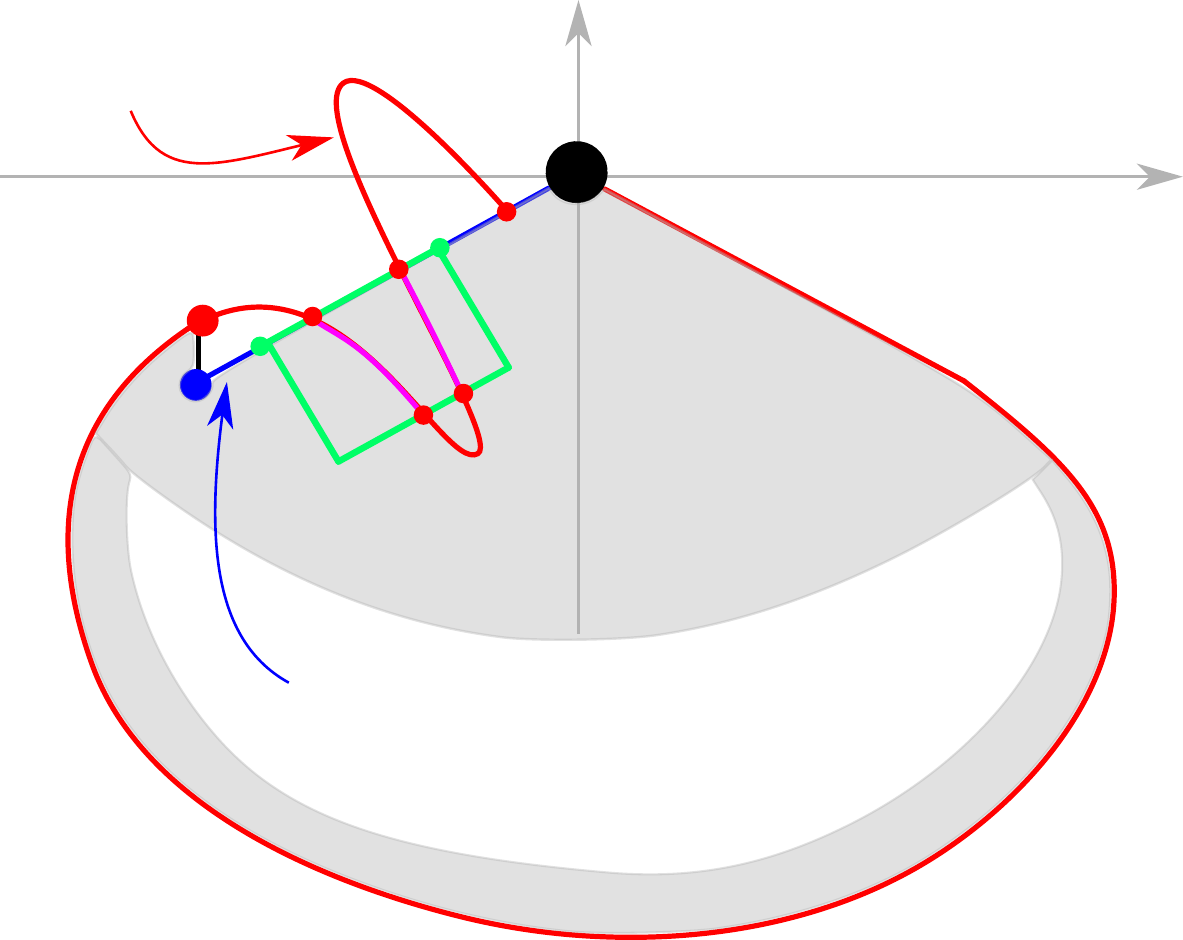}
			\put(100,61){{\footnotesize $y$}}
	\put(51,78){{\footnotesize $x$}}
		\put(25,20){{\footnotesize $W^{s,-}_{\Phi_0}$}}	
		\put(5,70){{\footnotesize $W^{u,-}_{\Phi_0}$}}
		\put(21,47){{\footnotesize $l_1$}}
		\put(36,60){{\footnotesize $l_2$}}		
			\put(30,52){{\footnotesize $L$}}
				\put(12,53){{\footnotesize $p^u_{\mathcal{R}}$}}
			\put(11,45){{\footnotesize $p^s_{\mathrm{loc}}$}}
			\put(25,49){{\footnotesize $\widehat{q_1}$}}
				\put(22,57){{\footnotesize $\phi_{X_0}(\widehat{q_2})$}}
						\put(42,58){{\footnotesize $\Phi_0(\widehat{q_1})$}}
							\put(30,45){{\footnotesize $A_1$}}
							\put(38,50){{\footnotesize $A_2$}}
								\put(27,37){{\footnotesize $Q$}}
								\put(55,35){$R_3$}
	\end{overpic}
\bigskip
\caption{Square $Q$ in the region $R_3$.}
\label{qfig}	
\end{figure}

In this case, we have the situation illustrated in Figure \ref{qfig}.

\begin{figure}[h!]
	\centering
	\bigskip
	\begin{overpic}[width=8cm]{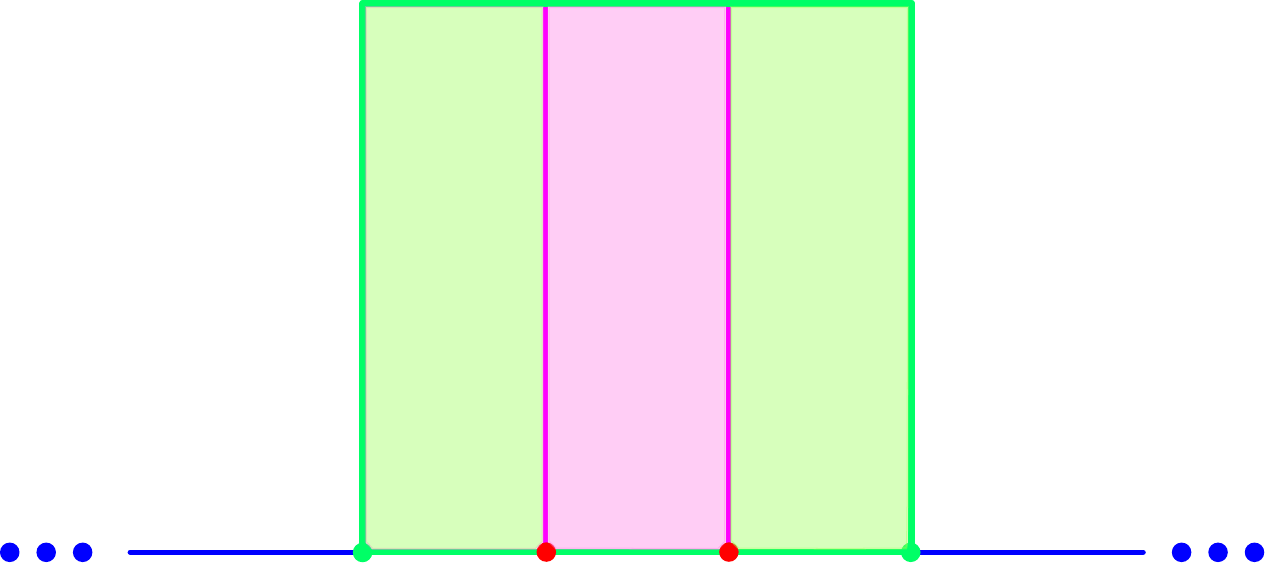}
		\put(10,-5){{\footnotesize $W^{s,-}_{\Phi_0}$}}	
		\put(28,-5){{\footnotesize $l_1$}}
		\put(72,-5){{\footnotesize $l_2$}}		
		\put(49,-7){{\footnotesize $L$}}
		\put(42,-5){{\footnotesize $\widehat{q_1}$}}
		\put(54,-5){{\footnotesize $\phi_{X_0}(\widehat{q_2})$}}
		\put(44,12){{\footnotesize $A_1$}}
		\put(58,12){{\footnotesize $A_2$}}
		\put(20,40){{\footnotesize $Q$}}
		\put(34,22){$Q_L$}
				\put(47,22){$Q_C$}
						\put(60,22){$Q_R$}
	\end{overpic}
	\vspace{0.8cm}
	
	\caption{Strips $Q_R$, $Q_C$ and $Q_L$ in the square $Q$.}
	\label{qstripfig}	
\end{figure}
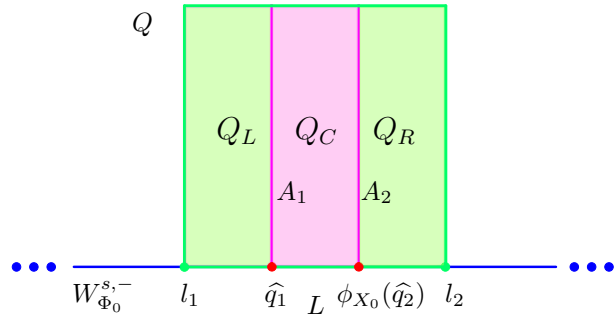

Notice that the segments $A_1$ and $A_2$ splits $Q$ into three strips. Denote the strip containing the segment between $l_1$ and $\widehat{q_1}$, $\widehat{q_1}$ and $\phi_{X_0}(\widehat{q_2})$, $\phi_{X_0}(\widehat{q_2})$ and $l_2$ by $Q_L$, $Q_C$ and $Q_R$, respectively. See Figure \ref{qstripfig}.

Now, we analyze the positive iterates of $Q$ by $\Phi_0$. First, notice that, since $Q$ is contained in the region $R_3$, it follows from the local behavior around a $T$-singularity that $\phi_{X_0}^n(Q)\cap V\subset\s^c$ and $\Phi_0^n(Q)\cap V\subset\s^c$, for every $n\in \N$, it means that the orbit of $Z_0$ connecting a point of $Q$ and some point of $(\phi_{X_0}^n(Q)\cup \Phi_0^n(Q))\cap V$, for some $n\in\N$, is a crossing orbit of $Z_0$.

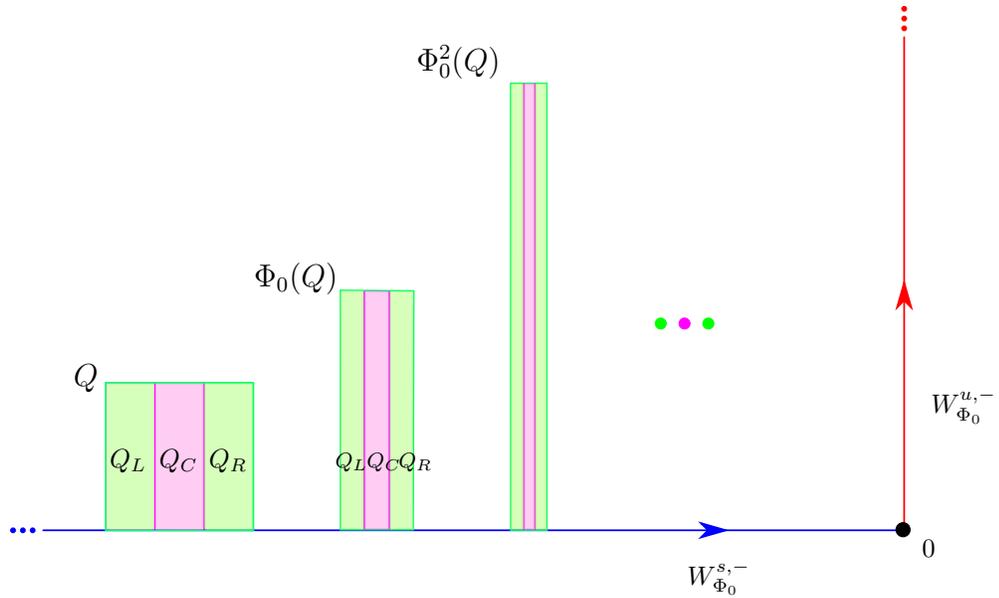
\begin{figure}[h!]
	\centering
	\bigskip
	\begin{overpic}[width=12cm]{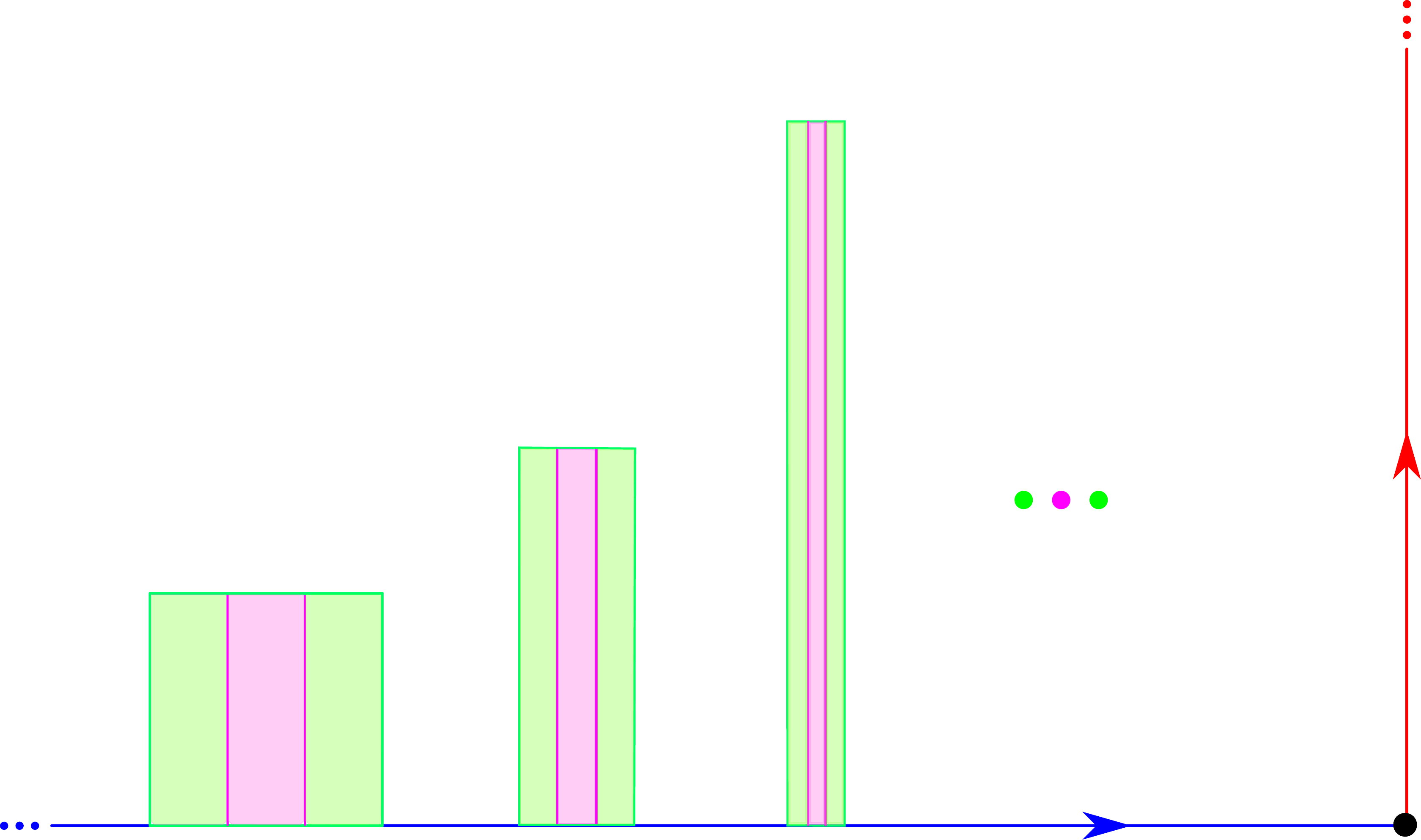}
				\put(101,-2){{\footnotesize $0$}}
		\put(75,-5){{\footnotesize $W^{s,-}_{\Phi_0}$}}	
			\put(102,14){{\footnotesize $W^{u,-}_{\Phi_0}$}}
		\put(11,8){{\footnotesize $Q_L$}}
		\put(16.5,8){{\footnotesize $Q_C$}}		
		\put(22,8){{\footnotesize $Q_R$}}
		\put(36,8){{\tiny $Q_L$}}
				\put(39.5,8){{\tiny $Q_C$}}
						\put(43,8){{\tiny $Q_R$}}					
		\put(7,17){$Q$}
		\put(27,28){$\Phi_0(Q)$}
				\put(45,52){$\Phi_0^2(Q)$}
	\end{overpic}
	\vspace{0.8cm}
	\caption{Positive iterates of $Q$ by $\Phi_0$.}
	\label{int1fig}	
\end{figure}

Take $W^{s,-}_{\Phi_0}$ as the horizontal direction and $W^{u,-}_{\Phi_0}$ as the vertical direction. Since $L$ is contained in $W^{s,-}_{\Phi_0}$, it follows from the $\lambda$-lemma that, when we apply $\Phi_0$, $Q$ is contracted in the horizontal direction and expanded in the vertical direction. See Figure \ref{int1fig}.

Now, we analyze the negative iterates of $Q$ by $\Phi_0$. 

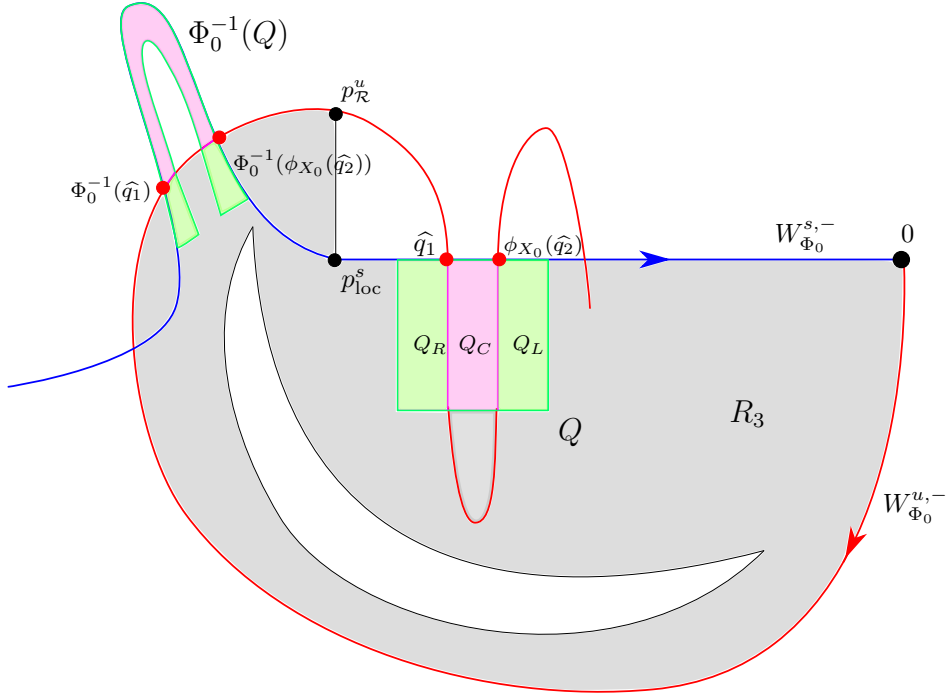
\begin{figure}[h!]
	\centering
	\bigskip
	\begin{overpic}[width=12cm]{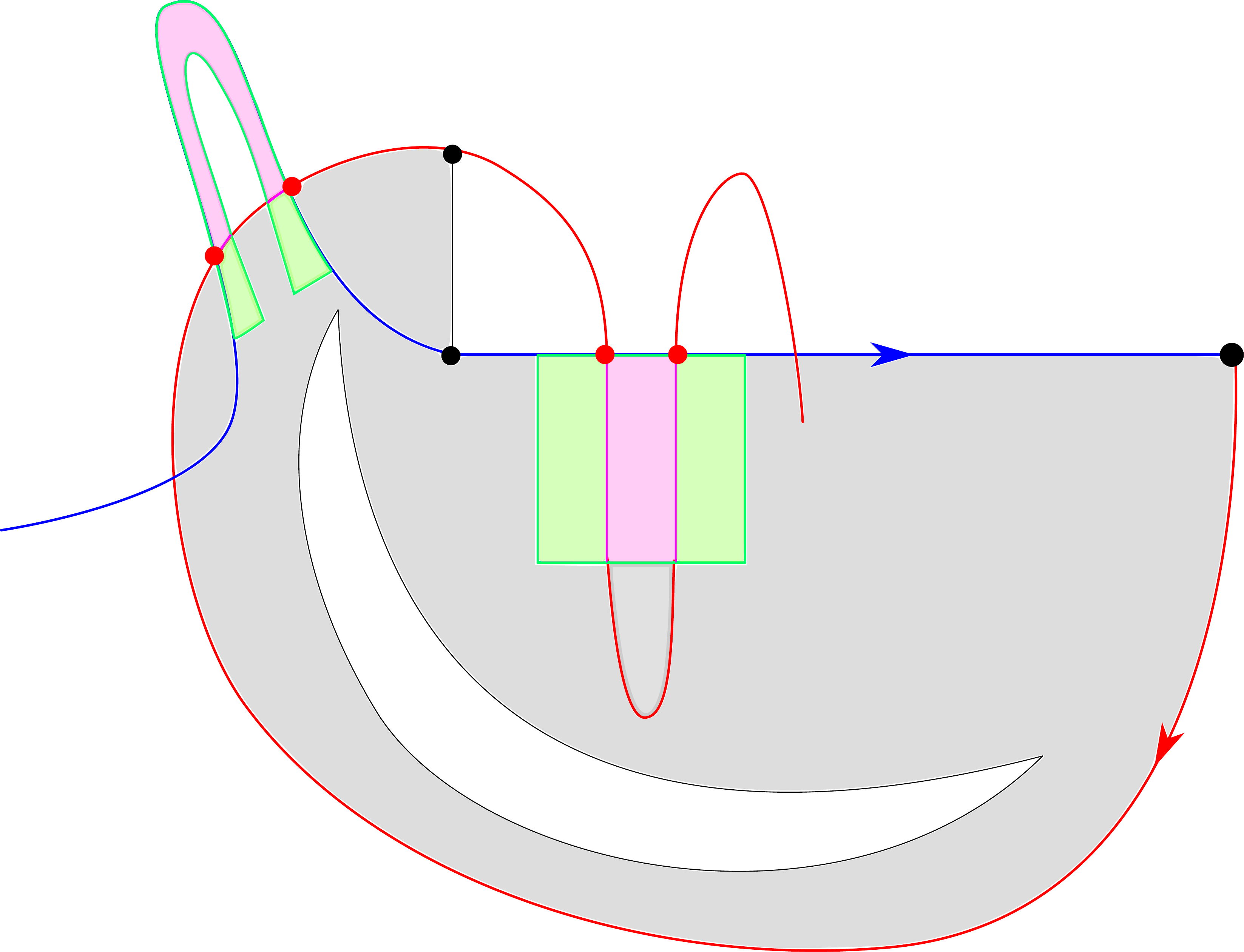}
							\put(37,66){{\scriptsize $p^u_{\mathcal{R}}$}}
							\put(37,45){{\footnotesize $p^s_{\mathrm{loc}}$}}
		\put(45,49){{\footnotesize $\widehat{q_1}$}}
		\put(99,50){{\footnotesize $0$}}
						\put(55,49){{\tiny $\phi_{X_0}(\widehat{q_2})$}}
								\put(61,28){$Q$}
								\put(80,30){$R_3$}
								\put(20,72){$\Phi_0^{-1}(Q)$}
								\put(25,58){{\tiny $\Phi_0^{-1}(\phi_{X_0}(\widehat{q_2}))$}}
								\put(7,55){{\tiny $\Phi_0^{-1}(\widehat{q_1})$}}
		\put(85,50){{\footnotesize $W^{s,-}_{\Phi_0}$}}	
		\put(97,20){{\footnotesize $W^{u,-}_{\Phi_0}$}}
		\put(45,38){{\tiny $Q_R$}}
		\put(50,38){{\tiny $Q_C$}}
		\put(56,38){{\tiny $Q_L$}}					
	\end{overpic}
	\bigskip
	\caption{Representation of $\Phi_0^{-1}(Q)$.}
	\label{bendedfig}	
\end{figure}

Notice that $W^{s,-}_{\Phi_0}$ can not cross $\tau^s$ in a point different from $p^s_{\mathrm{loc}}$, because $\tau^s$ is a transversal section. Hence, since the intersections between $W^{s,-}_{\Phi_0}$ and $W^{u,-}_{\Phi_0}$ are transversal, it follows that, all the points of $W^{s,-}_{\Phi_0}$ between $\widehat{q_1}$ and $\Phi_0^{-1}(\phi_{X_0}(\widehat{q_2}))$ are contained in the region $R_3$ and all the points between $\Phi_0^{-1}(\phi_{X_0}(\widehat{q_2}))$ and $\Phi_0^{-1}(\widehat{q_1})$ are outside $R_3$. Inductively, we show that the points of $W^{s,-}_{\Phi_0}$ between $\Phi_0^{-n}(\widehat{q_1})$ and $\Phi_0^{-(n+1)}(\phi_{X_0}(\widehat{q_2}))$ are contained in the region $R_3$ and all the points between $\Phi_0^{-n}(\phi_{X_0}(\widehat{q_2}))$ and $\Phi_0^{-n}(\widehat{q_1})$ are outside $R_3$. From this reasoning, it follows that $Q$ is bended into a horseshoe by $\Phi_0^{-1}$ in such a way that the strips $Q_L$ and $Q_R$ are contained in $R_3$ and $Q_C$ goes outside $R_3$. See Figure \ref{bendedfig}.

The same reasoning holds for all the negative iterations of $Q$ by $\Phi_0$. Therefore, $\Phi_0^{-n}(Q)$ is composed by two strips $Q_R$ and $Q_L$ contained in $R_3$ and a strip $Q_C$ outside $Q_3$. 

\begin{figure}[h!]
	\centering
	\bigskip
	\begin{overpic}[width=12cm]{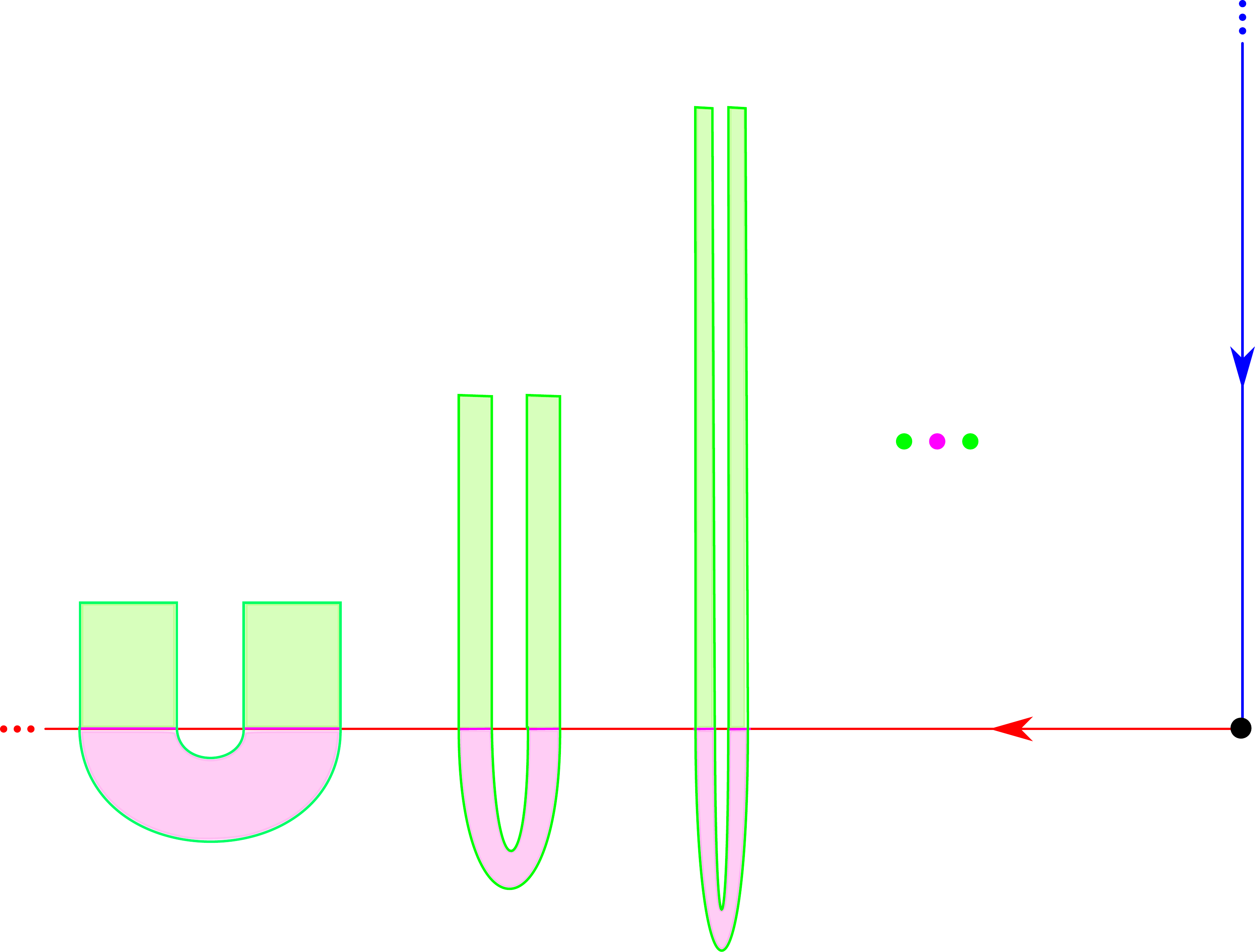}
		\put(101,15){{\footnotesize $0$}}
		\put(75,12){{\footnotesize $W^{u,-}_{\Phi_0}$}}	
		\put(102,32){{\footnotesize $W^{s,-}_{\Phi_0}$}}
		\put(8,22){{\footnotesize $Q_L$}}
		\put(15,11){{\footnotesize $Q_C$}}		
		\put(22,22){{\footnotesize $Q_R$}}
		\put(36.5,25){{\tiny $Q_L$}}
		\put(39,6){{\tiny $Q_C$}}
		\put(42,25){{\tiny $Q_R$}}					
		\put(5,30){$\Phi_0^{-1}(Q)$}
	\put(32,46){$\Phi_0^{-2}(Q)$}
	\put(52,69){$\Phi_0^{-3}(Q)$}
	\end{overpic}
	\bigskip
	\caption{Positive iterates of $Q$ by $\Phi_0$.}
	\label{int2fig}	
\end{figure}
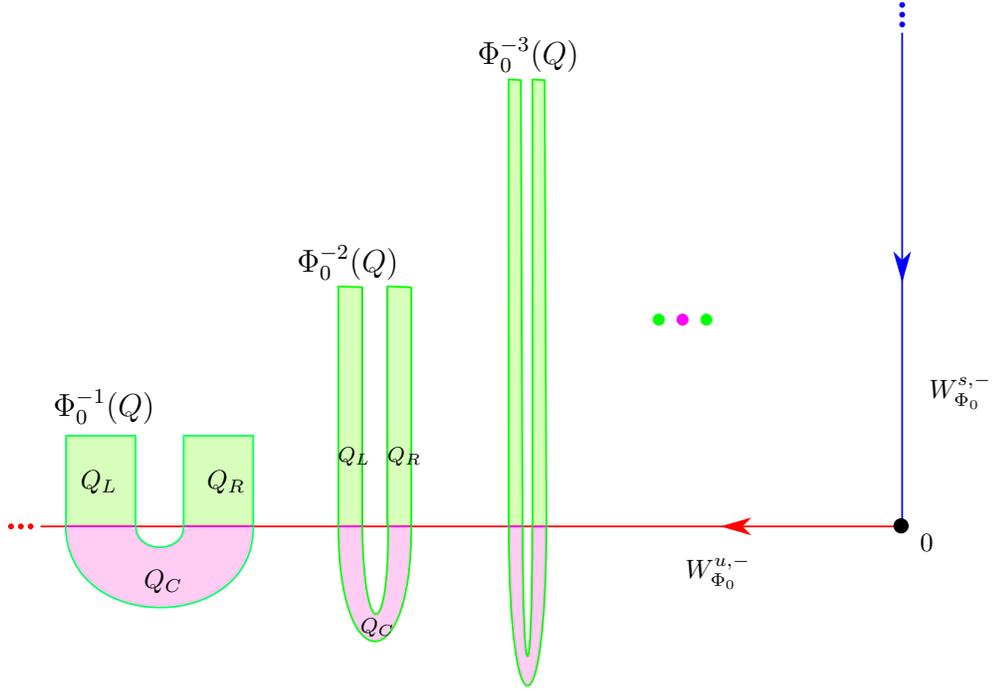

If we take $W^{u,-}_{\Phi_0}$ as the horizontal direction and  $W^{s,-}_{\Phi_0}$ as the vertical direction, using the $\lambda$-lemma, we obtain that, when we apply $\Phi_0^{-1}$ the strips $Q_R$ and $Q_L$ are contained in $R_3$ and they are contracted in the horizontal direction and expanded in the vertical direction. See Figure \ref{int2fig}.

 Notice that, the strips $Q_R$ and $Q_L$ are outside the cylinder $\mathcal{R}$, and thus it follows from the dynamics induced in the cylinder and the local dynamics of the $T$-singularity that, $\phi_{X_0}^{-n}(Q_{\star})\cap V\subset\s^c$ and $\Phi_0^{-n}(Q_{\star})\cap V\subset\s^c$, $\star=L,R$, for every $n\in\N$. Therefore, it follows that any orbit of $Z_0$ connecting a point of $Q_L\cup Q_R$ and $\phi_{X_0}^{-n}(Q_{L}\cup Q_R)$, $n\in\N$, is a crossing orbit of $Z_0$.

On the other hand, given $n\in\N$, an orbit of $Z_0$ passing through a point of $\Phi_0^{-n}(Q_{C})\cap V\subset\s^c$ reaches the sliding region at finite time.

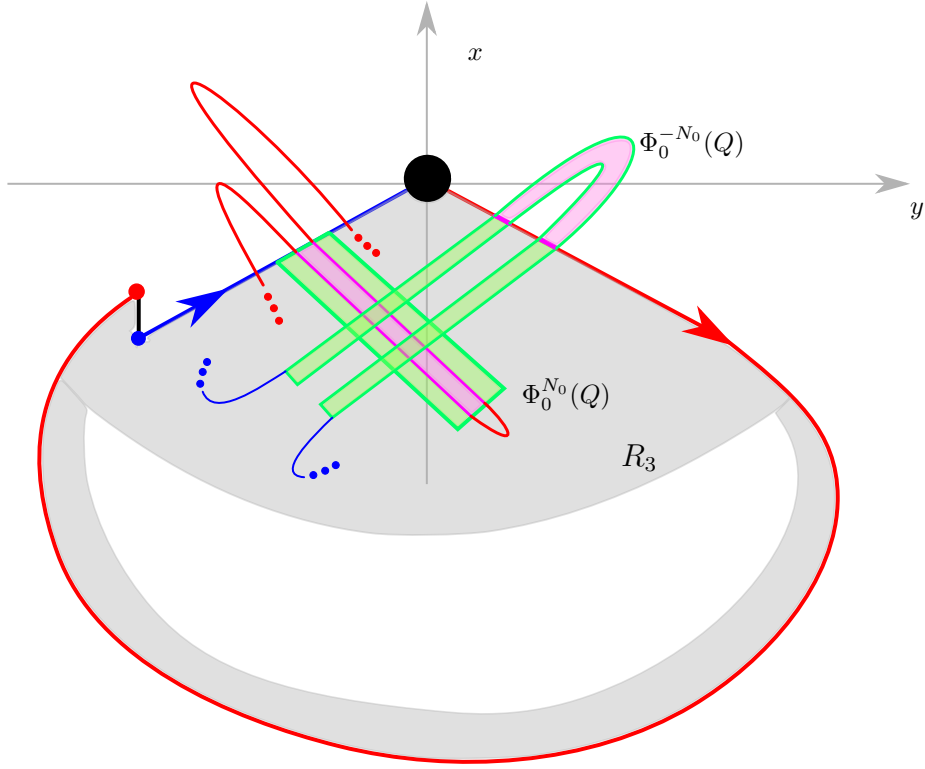
\begin{figure}[h!]
	\centering
	\bigskip
	\begin{overpic}[width=12cm]{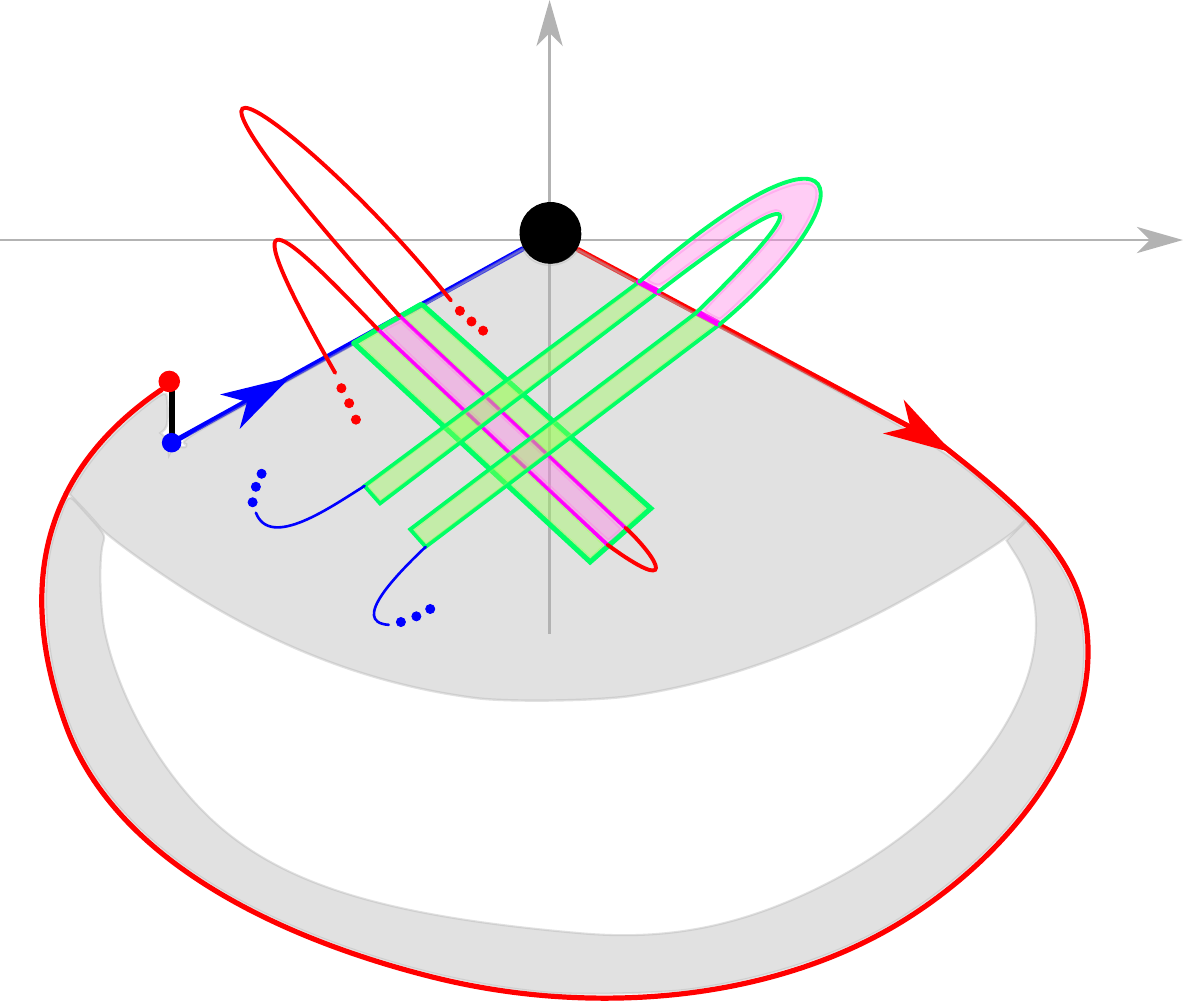}
		\put(100,61){{\footnotesize $y$}}
		\put(51,78){{\footnotesize $x$}}
		\put(57,40){{\footnotesize $\Phi_0^{N_0}(Q)$}}
			\put(70,68){{\footnotesize $\Phi_0^{-N_0}(Q)$}}
		\put(68,33){$R_3$}
	\end{overpic}
	\bigskip
	\caption{Action of the horseshoe map $\Delta$ on the square $Q$.}
	\label{horsefig}	
\end{figure}

Hence, there exists $N_0\in\N$ such that $\Delta=\Phi_0^{N_0}: Q\rightarrow \rn{2}$ is a horseshoe map acting on $Q$ (see Figure \ref{horsefig}). Finally, using the strips $Q_R$ and $Q_L$, it follows that
the invariant cantor set $\Lambda=\displaystyle\cup_{n=-\infty}^{\infty} \Delta^n(Q_R\cup Q_L)\cap V$ given by Theorem $6.5.5$ of \cite{KATOK} is contained in $\s^c$ and every orbit of $Z_0$ passing through $\Lambda$ is a crossing orbit of $Z_0$. This concludes the proof of item $(3)$ of Theorem \ref{ferradura}.

\begin{remark}
	Notice that the square $Q$ used in this construction can be considered in any neighborhood $U$ of the origin  (which is a T-singularity) in a similar way. Therefore, for any neighborhood  $U$ of the origin, we can construct a horseshoe map $\Delta$ acting on $Q\subset U$ such that, the hyperbolic set $\Lambda$ associated to $\Delta$ satisfies item $(3)$ of Theorem \ref{ferradura}.

\end{remark}

\section{A Model Presenting a $T$-connection}\label{model_sec}

In this section we construct a Filippov system having two robust fold-fold connections. It is worth mentioning that such model can be used to produce examples of $T$-chains satisfying $(TC)$ and $(R)$ conditions.

\subsection{Filippov System $Z_1$}
Consider the Filippov system

\begin{equation}\label{z1}
Z_1(x,y,z)=\left\{
\begin{array}{lcl}
X_1(x,y,z)=(1,-1,y),& \textrm{ if }z>0,\vspace{0.2cm}\\
Y_1(x,y,z)=(-1,2,-x),& \textrm{ if }z<0.
\end{array}
\right.
\end{equation}

\begin{figure}[h!]
	\centering
	\bigskip
	\begin{overpic}[width=6cm]{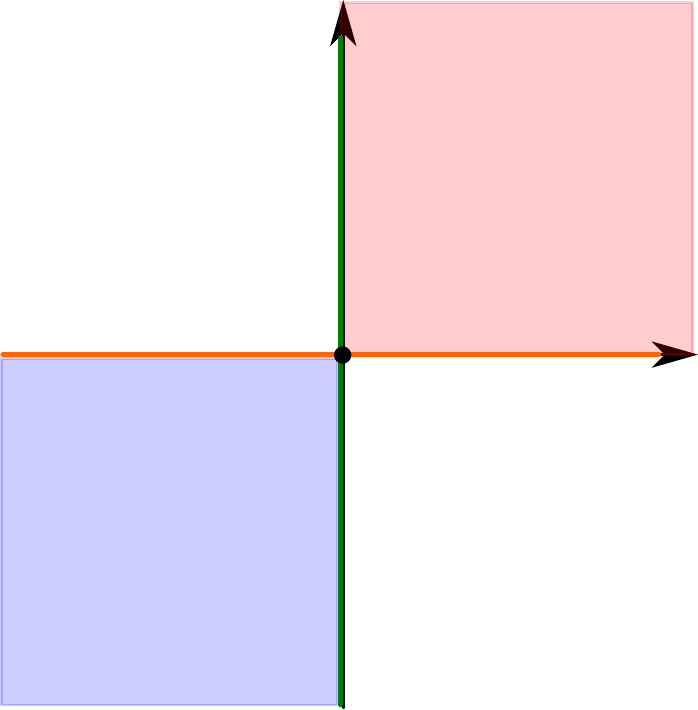}
		\put(101,49){{\footnotesize $x$}}		
		\put(49,102){{\footnotesize $y$}}		
		\put(80,42){{\footnotesize $S_{X_1}$}}	
		\put(35,85){{\footnotesize $S_{Y_1}$}}	
		\put(70,70){{\footnotesize $\s^{us}$}}	
		\put(20,20){{\footnotesize $\s^{ss}$}}														
	\end{overpic}
	\bigskip
	\caption{Switching manifold associated to $Z_1$.} \label{s1fig}	
\end{figure}
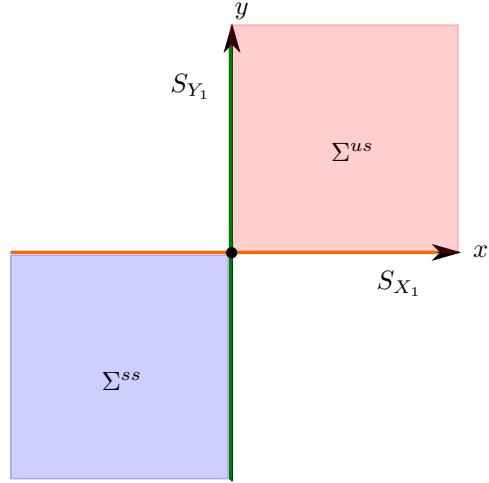

Notice that $(0,0,0)$ is a $T$-singularity and $S_{X_1}=\{y=0\}$ and $S_{Y_1}=\{x=0\}$ are fold lines which divide the switching manifold $\s=\{z=0\}$ in four quadrants (see Figure \ref{s1fig}). In addition
\begin{itemize}
	\item $\s^{ss}=\{(x,y,0);\ x,y<0\};\vspace{0.2cm}$
	\item $\s^{us}=\{(x,y,0);\ x,y>0\};\vspace{0.2cm}$
	\item $\s^c=\{(x,y,0);\ xy<0\};\vspace{0.2cm}$		
\end{itemize}

A straightforward computation shows that the flow of $X_1$ and $Y_1$ are given by
\begin{equation}
\p_{X_1}(t;(x,y,z))=\left(\begin{array}{c}
x+t\\
y-t\\
z-\dfrac{(y-t)^2}{2}+\dfrac{y^2}{2}
\end{array}\right),
\end{equation}
and 
\begin{equation}
\p_{Y_1}(t;(x,y,z))=\left(\begin{array}{c}
x-t\\
y+2t\\
z+\dfrac{(x-t)^2}{2}-\dfrac{x^2}{2}
\end{array}\right),
\end{equation}
respectively. It allows us to see that 
$$\p_{X_1}(t;(x,y,0))\in\s\Longleftrightarrow t=0 \textrm{ or } t=2y,$$
and
$$\p_{Y_1}(t;(x,y,0))\in\s\Longleftrightarrow t=0 \textrm{ or } t=2x.$$

Thus, we can associate involutions $\phi_{X_1},\phi_{Y_1}:\s\rightarrow\s$ to the vector fields $X_1$ and $Y_1$ respectively, which are given by

\begin{equation}
\phi_{X_1}(x,y)=\left(\begin{array}{c}
x+2y\\
-y
\end{array}\right), \textrm{ and } \phi_{Y_1}(x,y)=\left(\begin{array}{c}
-x\\
4x+y
\end{array}\right).
\end{equation}

Now, we use the involutions to construct the first return map $\phi_1:\s\rightarrow\s$ given by
\begin{equation}
\phi_1(x,y)=\phi_{Y_1}\circ\phi_{X_1}(x,y)=\left(\begin{array}{c}
-x-2y\\
4x+7y
\end{array}\right).
\end{equation}

Notice that $\phi_1$ is globally defined on $\s$. The eigenvalues of $\phi_1$ are given by
\begin{equation}
\lambda_1^{\pm}=3\pm 2\sqrt{2},
\end{equation}
and their respective eigenvectors are
\begin{equation}
v_1^{\pm}=(-1\pm \sqrt{2}/2,1).
\end{equation}

Since $\phi_1$ is linear, the lines
$$D_1^{\pm}=\{\ag(-1\pm \sqrt{2}/2,1);\ \ag\in\R \},$$
are global invariant manifolds of $\phi_1$. Furthermore, 
$$\phi_{X_1}(D_1^-)=D_1^+ \textrm{ and } \phi_{Y_1}(D_1^+)=D_1^-.$$

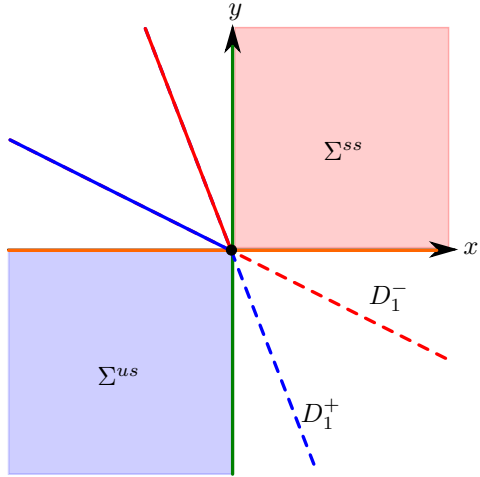
\begin{figure}[h!]
	\centering
	\bigskip
	\begin{overpic}[width=6cm]{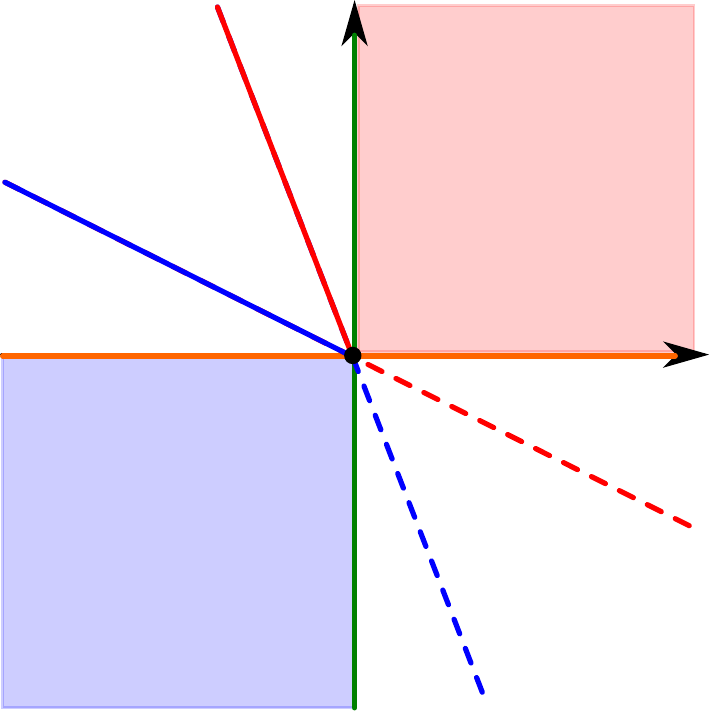}
		\put(101,49){{\footnotesize $x$}}		
		\put(49,102){{\footnotesize $y$}}		
		\put(80,38){{\footnotesize $D_1^-$}}	
		\put(65,12){{\footnotesize $D_1^+$}}	
		\put(70,70){{\footnotesize $\s^{ss}$}}	
		\put(20,20){{\footnotesize $\s^{us}$}}														
	\end{overpic}
	\bigskip
	\caption{Illustration of $\mathcal{N}_1\cap \s$.}\label{s2fig}	
\end{figure}

These facts ensure the existence of global crossing invariant manifolds in the form of a nonsmooth diabolo $\mathcal{N}_1$ which intersects $\s$ in the lines $D_1^{\pm}$. See Figure \ref{s2fig}.

\subsection{Filippov System $Z_2$}
Consider the Filippov system

\begin{equation}\label{z2}
Z_2(x,y,z)=\left\{
\begin{array}{lcl}
X_2(x,y,z)=(1,-1,y-2),& \textrm{ if }z>0,\vspace{0.2cm}\\
Y_2(x,y,z)=(-1,3,-(x-2)),& \textrm{ if }z<0.
\end{array}
\right.
\end{equation}

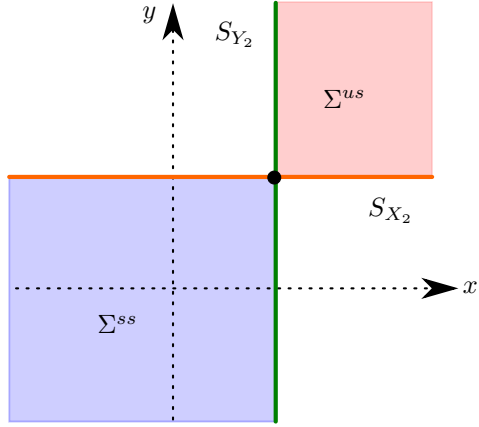
\begin{figure}[h!]
	\centering
	\bigskip
	\begin{overpic}[width=6cm]{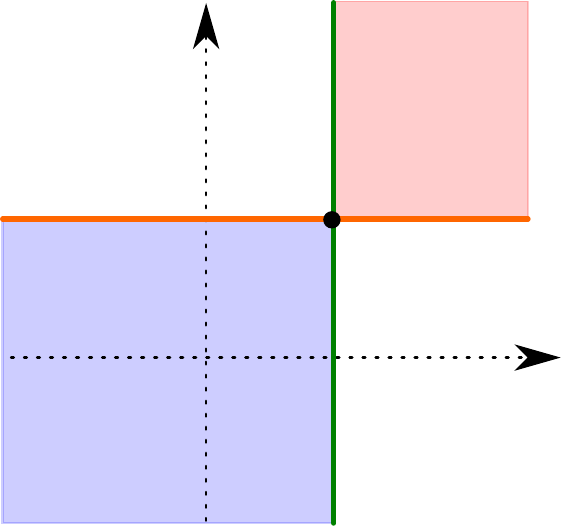}
		\put(101,29){{\footnotesize $x$}}		
		\put(30,90){{\footnotesize $y$}}		
		\put(80,46){{\footnotesize $S_{X_2}$}}	
		\put(46,85){{\footnotesize $S_{Y_2}$}}	
		\put(70,70){{\footnotesize $\s^{us}$}}	
		\put(20,20){{\footnotesize $\s^{ss}$}}														
	\end{overpic}
	\bigskip
	\caption{Switching manifold associated to $Z_2$.} \label{s3fig}	
\end{figure}

Notice that $(2,2,0)$ is a $T$-singularity and $S_{X_2}=\{y=2\}$ and $S_{Y_2}=\{x=2\}$ are fold lines which divide the switching manifold $\s=\{z=0\}$ in four quadrants (see Figure \ref{s3fig}). In addition
\begin{itemize}
	\item $\s^{us}=\{(x,y,0);\ x,y>2\};\vspace{0.2cm}$
	\item $\s^{ss}=\{(x,y,0);\ x,y<2\};\vspace{0.2cm}$
	\item $\s^c=\{(x,y,0);\ (x-2)(y-2)<0\};\vspace{0.2cm}$		
\end{itemize}

A straightforward computation shows that the flows of $X_2$ and $Y_2$ are given by
\begin{equation}
\p_{X_2}(t;(x,y,z))=\left(\begin{array}{c}
x+t\\
y-t\\
z-\dfrac{(y-t)^2}{2}+\dfrac{y^2}{2}- 2t
\end{array}\right),
\end{equation}
and 
\begin{equation}
\p_{Y_2}(t;(x,y,z))=\left(\begin{array}{c}
x-t\\
y+3t\\
z+\dfrac{(x-t)^2}{2}-\dfrac{x^2}{2}+2t
\end{array}\right),
\end{equation}
respectively. It allows us to see that 
$$\p_{X_2}(t;(x,y,0))\in\s\Longleftrightarrow t=0 \textrm{ or } t=2(y-2),$$
and
$$\p_{Y_2}(t;(x,y,0))\in\s\Longleftrightarrow t=0 \textrm{ or } t=2(x-2).$$

As before, we can associate involutions $\varphi_{X_2},\varphi_{Y_2}:\s\rightarrow\s$ associated to the vector fields $X_2$ and $Y_2$, respectively, which are given by

\begin{equation}
\phi_{X_2}(x,y)=\left(\begin{array}{c}
x+2(y-2)\\
2-(y-2)
\end{array}\right), \textrm{ and } \phi_{Y_2}(x,y)=\left(\begin{array}{c}
2-(x-2)\\
y+6(x-2)
\end{array}\right).
\end{equation}

Now, use the involutions to construct the first return map $\phi_2:\s\rightarrow\s$ given by
\begin{equation}
\phi_2(x,y)=\phi_{Y_2}\circ\phi_{X_2}(x,y)=\left(\begin{array}{c}
2\\
2
\end{array}\right)+\left(\begin{array}{cc}
-1&-2\\
6&11
\end{array}\right)\left(\begin{array}{c}
x-2\\
y-2
\end{array}\right).
\end{equation}

The eigenvalues of $\phi_2$ at the fixed point $(2,2)$ are given by
\begin{equation}
\lambda_2^{\pm}=5\pm 2\sqrt{6},
\end{equation}
and their respective eigenvectors are
\begin{equation}
v_2^{\pm}=(-1\pm \sqrt{6}/3,1).
\end{equation}

Since $\phi_2$ is linear, the lines
$$D_2^{\pm}=\{(2,2)+\ag(-1\pm \sqrt{6}/3,1);\ \ag\in\R \},$$
are global invariant manifolds of $\phi_2$. Furthermore, 
$$\phi_{X_2}(D_2^-)=D_2^+ \textrm{ and } \phi_{Y_2}(D_2^+)=D_2^-.$$

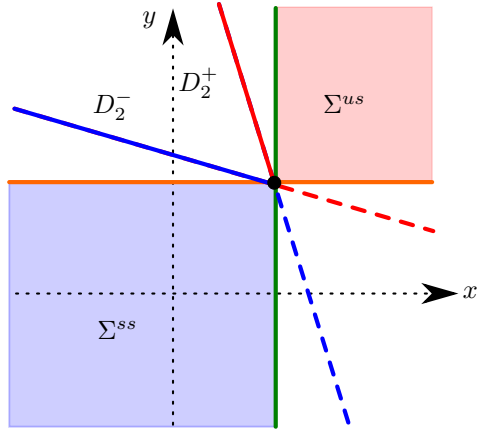
\begin{figure}[h!]
	\centering
	\bigskip
	\begin{overpic}[width=6cm]{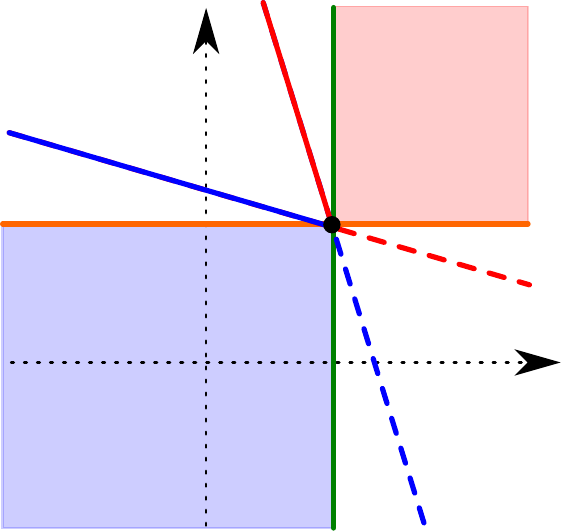}
		\put(101,29){{\footnotesize $x$}}		
		\put(30,90){{\footnotesize $y$}}		
		\put(38,75){{\footnotesize $D_2^+$}}	
		\put(19,70){{\footnotesize $D_2^-$}}	
		\put(70,70){{\footnotesize $\s^{us}$}}	
		\put(20,20){{\footnotesize $\s^{ss}$}}														
	\end{overpic}
	\bigskip
	\caption{Illustration of $\mathcal{N}_2\cap \s$.}\label{s4fig}
\end{figure}

These facts ensure the existence of a nonsmooth diabolo $\mathcal{N}_2$ which intersects $\s$ in $D_2^{\pm}$. See Figure \ref{s4fig}.

\subsection{A Filippov System with a Cross Shaped Switching Manifold}

Consider a new switching manifold $\Pi=g^{-1}(0)$, where 
\begin{equation}\label{g}
g(x,y,z)=-(y+17/10 x-5/2). 
\end{equation}

Let $W^u_1(0,0,0)$ be the unstable branch of the nonsmooth diabolo $\mathcal{N}_1$ of $Z_1$ at the origin. In this case, $W^u_1(0,0,0)\cap\{z>0\}$ is given by the parametrized set
\begin{equation}
\mathcal{W}^{u,+}_1=\{(\ag(-1+\sqrt{2}/2)\ag+t,\ag-t,-(\ag-t)^2/2+\ag^2/2);\ 0\leq t\leq 2\ag \textrm{ and }\ag\geq 0\}.
\end{equation}

Similarly, denoting the stable branch of the nonsmooth diabolo $\mathcal{N}_2$ of $Z_2$ at $(2,2,0)$ by $W^s_2(2,2,0)$, we obtain that $W^s_2(2,2,0)\cap\{z>0 \}$ is given by 
\begin{equation}
\mathcal{W}^{s,+}_2=\{(2+(-1-\sqrt{2/3})\ag+t,2+\ag-t,(2+\ag)^2/2-(2+\ag-t)^2/2-2t);\ 0\leq t\leq 2\ag \textrm{ and }\ag\geq 0\}.
\end{equation}

A straight computation shows that 
\begin{equation}
\mathcal{C}^u_+=\mathcal{W}^{u,+}_1\cap \Pi=\{\cg^u_+(\ag);\ 25/191 (-14 + 17 \sqrt{2})\leq\ag\leq 25/191 (14 + 17 \sqrt{2})  \},
\end{equation}
where
\begin{equation}
\cg^u_+(\ag)=(-5/7(-5+\sqrt{2}\ag),1/14(-50+17\sqrt{2}\ag),1/196(-1250+850\sqrt{2}\ag-191\ag^2)).
\end{equation}

Also
\begin{equation}
\mathcal{C}^s_-=\mathcal{W}^{s,+}_2\cap \Pi=\{\cg^s_+(\ag);\ 29/431 (-21 + 17 \sqrt{6})\leq\ag\leq 29/431 (21 + 17 \sqrt{6})  \},
\end{equation}
where
\begin{equation}
\cg^s_+(\ag)=(5/21(-9+2\sqrt{6}\ag),1/21(129-17\sqrt{6}\ag),1/294(-2523+986\sqrt{6}\ag-431\ag^2)).
\end{equation}

\begin{figure}[h!]
	\centering
	\bigskip
	\begin{overpic}[width=7.5cm]{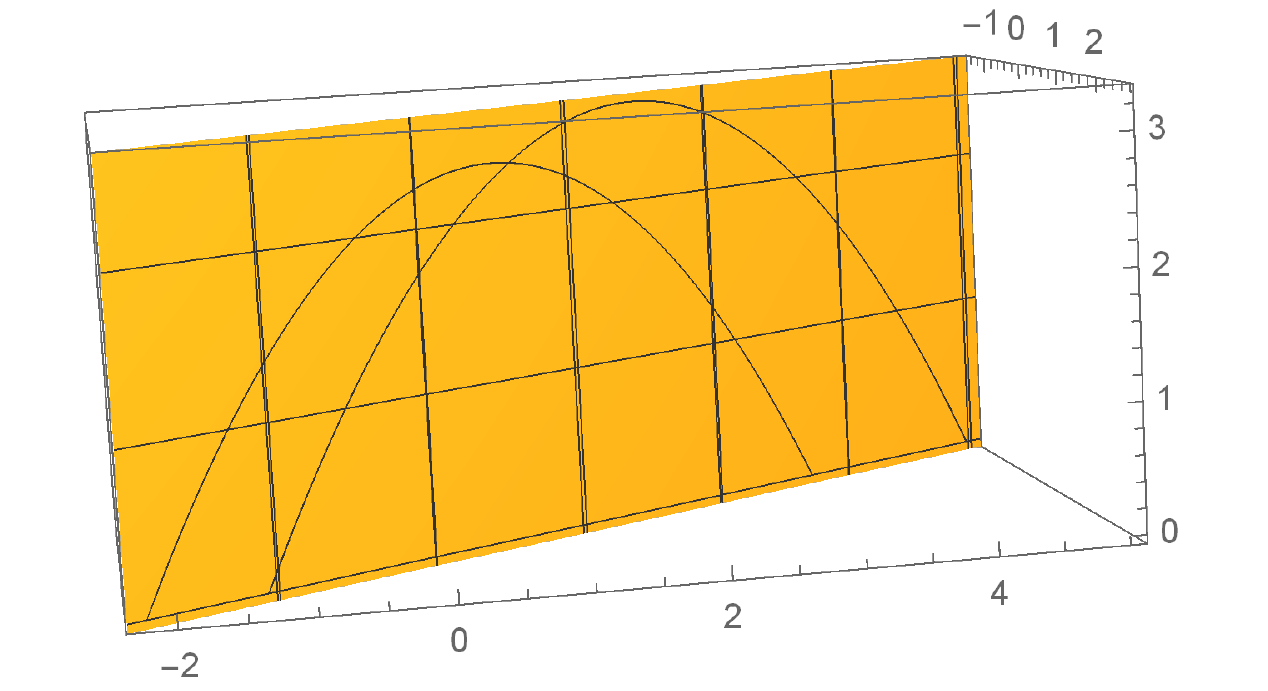}
		\put(1,35){\footnotesize $\Pi$}
		\put(25,38){\footnotesize $\mathcal{C}^s_+$}
		\put(60,44){\footnotesize $\mathcal{C}^u_+$}
		\put(38,43){\footnotesize $p_+^*$}	
	\end{overpic}
	\bigskip
	\caption{Sketch of $\mathcal{C}^u_+$ and $\mathcal{C}^s_-$. }\label{paafig}	
\end{figure}
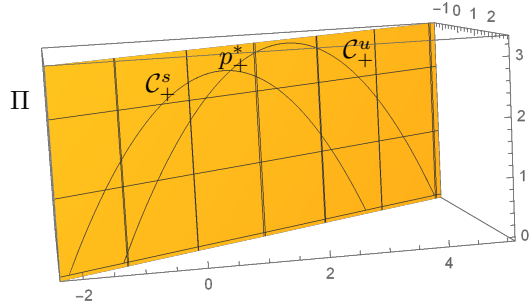

In addition
\begin{equation}\label{pstar}
\mathcal{C}^u_+\cap\mathcal{C}^s_-=p_+^* =\left({-5/49 (-67 + 8 \sqrt{51}), 
	1/49 (-447 + 68 \sqrt{51}), (-330577 + 48248 \sqrt{51})/4802}\right),
\end{equation}
and $\mathcal{C}^u_+\pitchfork \mathcal{C}^s_-$ at $p_+^*$. See Figure \ref{paafig}.

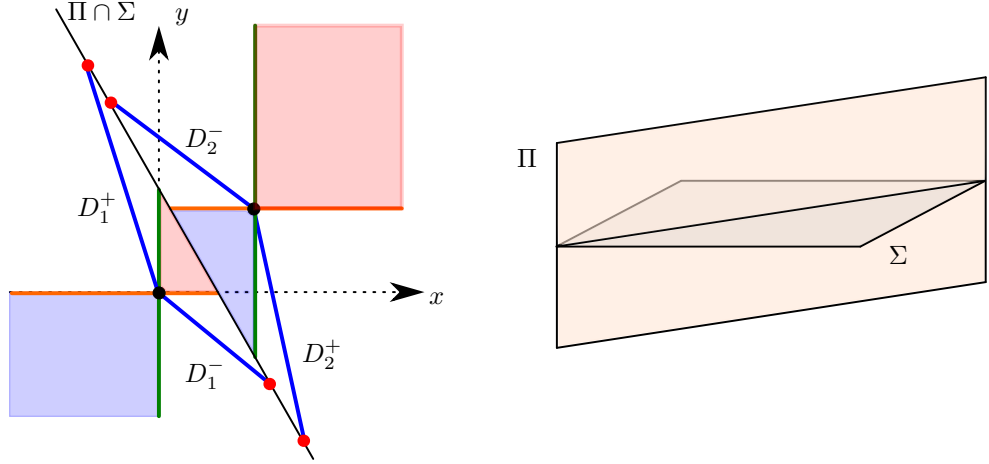
\begin{figure}[h!]
	\centering
	\bigskip
	\begin{overpic}[width=13cm]{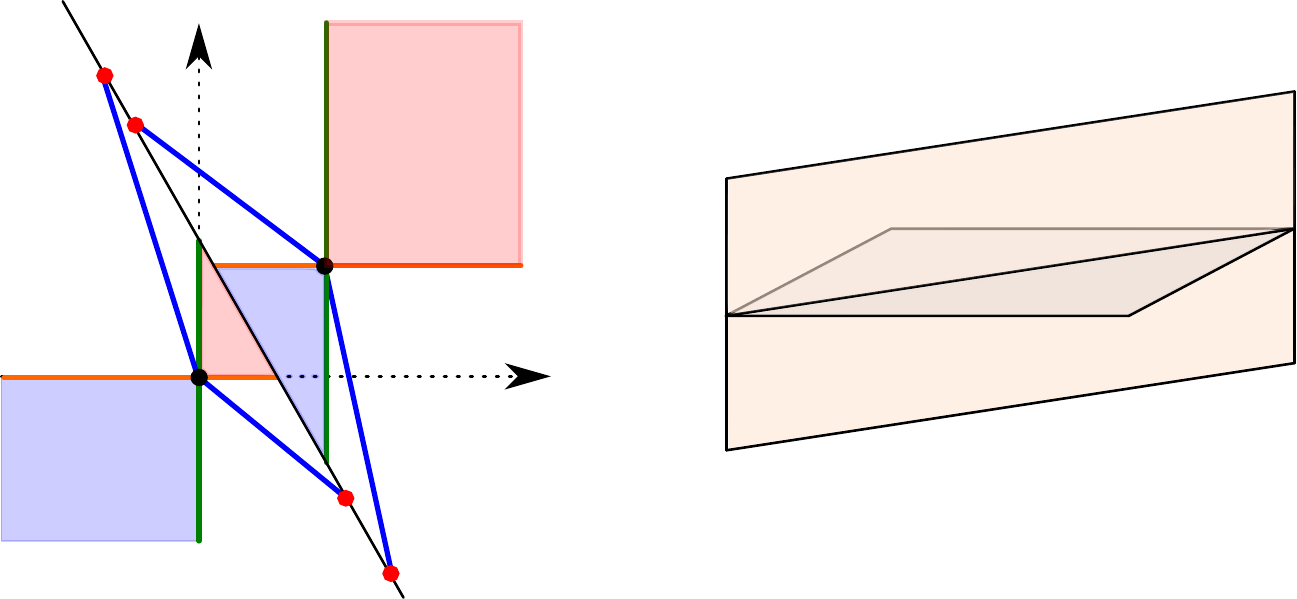}
		\put(43,16){{\footnotesize $x$}}
		\put(17,45){{\footnotesize $y$}}
		\put(6,45){{\footnotesize $\Pi\cap\s$}}		
		\put(18,32){{\footnotesize $D_2^-$}}
		\put(30,10){{\footnotesize $D_2^+$}}
		\put(7,25){{\footnotesize $D_1^+$}}
		\put(18,8){{\footnotesize $D_1^-$}}	
		\put(52,30){{\footnotesize $\Pi$}}	
		\put(90,20){{\footnotesize $\s$}}										
	\end{overpic}
	\bigskip
	\caption{$(a)$ Switching manifold $\s$ associated to $Z_0$. $(b)$ Cross-shaped switching manifold $\s\cup\Pi$.}	\label{s6fig}
\end{figure}

Also, notice that the vector fields $X_1,X_2,Y_1$ and $Y_2$ are transverse to $\Pi$ at every point of $\Pi$. Thus, the piecewise smooth system

\begin{equation}
Z_0(x,y,z)=\left\{
\begin{array}{lcl}
Z_1(x,y,z),& \textrm{ if }g(x,y,z)<0,\vspace{0.2cm}\\
Z_2(x,y,z),& \textrm{ if }g(x,y,z)>0.
\end{array}
\right.
\end{equation}
with a cross-shaped switching manifold has an isolated crossing orbit connecting the $T$-singularities $(0,0,0)$ and $(2,2,0)$ of $Z_0$ passing through $p_+^*$ given in \eqref{pstar}. See Figure \ref{s6fig}.

Analogous conclusions can be shown for $z<0$. In this case, $Z_0$ has another isolated crossing orbit connecting the $T$-singularities $(0,0,0)$ and $(2,2,0)$ of $Z_0$ passing through a point $p_-^*\in\Pi$ contained in $z<0$.
\subsection{A Filippov system presenting fold-fold connections}

Consider the $\mathcal{C}^1$-regularization function
\begin{equation}
\p(x)=\left\{\begin{array}{cl}
-1,& \textrm{if }x<-1,\\
\sin(\pi/2 x),& \textrm{if }|x|\leq 1,\\
1,& \textrm{if }x>1.
\end{array}\right.
\end{equation}

Thus, regularizing the Filippov systems $(X_1,X_2,\Pi)$ and $(Y_1,Y_2,\Pi)$ with respect to the switching manifold $\Pi$, we obtain
\begin{equation}
\mathcal{X}_\e(x,y,z)=\dfrac{X_2(x,y,z)+X_1(x,y,z)}{2}+\p\left(\dfrac{g(x,y,z)}{\e}\right)\dfrac{X_2(x,y,z)-X_1(x,y,z)}{2}
\end{equation}
and
\begin{equation}
\mathcal{Y}_\e(x,y,z)=\dfrac{Y_2(x,y,z)+Y_1(x,y,z)}{2}+\p\left(\dfrac{g(x,y,z)}{\e}\right)\dfrac{Y_2(x,y,z)-Y_1(x,y,z)}{2},
\end{equation}
where $X_1,Y_1$ are given by \eqref{z1}, $X_2,Y_2$ are given by \eqref{z2} and $g$ is given by \eqref{g}.

Thus, for $\e>0$, $\mathcal{Z}_\e=(\mathcal{X}_\e,\mathcal{Y}_\e)$ is a Filippov system with switching manifold $\s=\{z=0\}$ and $\mathcal{Z}_\e\in\Omega^1$.

\begin{remark}
	If we consider a regularizing function $\p$ which is of class $\Cr$, then $\mathcal{Z}_{\e}\in\Or$. 
\end{remark}

It follows that $\mathcal{Z}_\e$ has two (stable) $T$-singularities at $p_1=(0,0,0)$ and $p_2=(2,2,0)$. Now, since the invariant manifolds $W_1^u(p_1)$ and $W_2^s(p_2)$ intersect $\Pi$ transversally in two topological circles $\mathcal{C}^u$ and $\mathcal{C}^s$, we have that:

\begin{enumerate}
	\item The unstable crossing invariant manifold $W^u_\e(p_1)$ of $\mathcal{Z}_\e$ at $p_1$ intersects the transversal section $\Pi_{-\e}=\{g(x,y,z)=-\e\}$ in a topological circle $\mathcal{C}^u_{\e}$;
	\item The stable crossing invariant manifold $W^s_\e(p_2)$ of $\mathcal{Z}_\e$ at $p_2$ intersects the transversal section $\Pi_{\e}=\{g(x,y,z)=\e\}$ in a topological circle $\mathcal{C}^s_{\e}$.	
\end{enumerate}

Consider a small annulus $D^u_{\e}$ around $\mathcal{C}^u_{\e}$ contained in the plane $\Pi_{-\e}$. Now, $\mathcal{X}_{\e}$ and $\mathcal{Y}_{\e}$ are transverse to $\Pi_{-\e}$, $\Pi_{\e}$ and $z=0$  and there are no singularities of $\mathcal{Z}_{\e}$ in the cylindrical region delimited by $D^u_{\e}$ and the regularization zone $R_{\e}=\{ (x,y,z);\ |g(x,y,z)|<\e  \}$. It means that the flow of $\mathcal{Z}_{\e}$ is tubular inside this region (it has only crossing orbits which goes from $D^u_{\e}$ to $\Pi_{\e}$.

Therefore, $W^u_\e(p_1)$ extends itself in the regularization zone through crossing orbits of $\mathcal{Z}_{\e}$, and it intersects $\Pi_{\e}$ in a topological circle $\widehat{\mathcal{C}^u_{\e}}$.

Since $\mathcal{C}^u$ and $\mathcal{C}^s$ intersect transversally at the points $p_{-}^*$ and $p_{+}^*$, it follows that, for $\e>0$ sufficiently small, $\widehat{\mathcal{C}^u_{\e}}$ and $\mathcal{C}^s_{\e}$ intersect themselves transversally at two points $q_{1}(\e)\in\s^-$ and $q_2(\e)\in\s^+$.

It means that, there exists $\e_0>0$ sufficiently small such that, for each $\e\leq\e_0$, the one-parameter family $\mathcal{Z}_{\e}\in\Omega^1$ has two robust fold-fold connections between the $T$-singularities $p_1$ and $p_2$.

\section{Further Directions}\label{further_sec}

In this paper, we have presented a global phenomenon involving connections at a T-singularity $p$ of a system $Z$, which generates a chaotic behavior in the crossing orbits of $Z$. 

There are some models of Filippov systems which can present non-deterministic chaos at a stable T- singularity mainly due to the behavior of the sliding vector field at such point (see \cite{CJ2}). Therefore, one should investigate how the sliding dynamics of a Filippov systems presenting a robust T-chain interacts with the hyperbolic invariant set $\Lambda$ associated to the Smale horseshoe of the first return map $\Phi_0$ found in this paper.  It is an arduous task which might bring new ways towards the comprehension of chaos for Filippov systems.

The notion of chaos in general Filippov systems is still poorly understood taking in account the richness of the dynamics generated by discontinuities. Most of works exploring this subject (see \cite{BCE16,Nshil15, G18,NPV17} for instance) relies on the use of classical approaches to characterize chaotic behavior, as it was partially done in this paper. Nevertheless, a generalization of the concept of chaos is called for in this approach. In our point of view, some hidden and unknown objects can appear and a classical setting may not be sufficient to study them.

Based on these lines, we think that a general interesting problem is to provide a chaotic notion for Filippov systems which takes into account all the non-standard properties of these systems, such as sliding phenomena, non-uniqueness of solutions and $\s$-singularities.

\appendix

\section{Proof of Lemma \ref{extlema}}\label{appendixA}
We prove only item $(ii)$, since item $(iii)$ is proved in an analogous way and item $(iv)$ is a direct consequence of items $(ii)$ and $(iii)$. Item $(i)$ will follows from the construction. Consider the notation introduced in Sections \ref{tchainssec} and \ref{robustsec}.
 
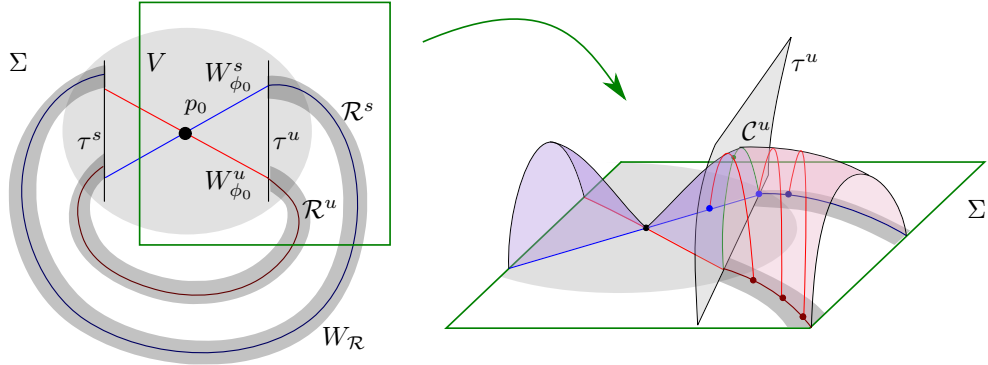
\begin{figure}[h!]
	\centering
	\bigskip
	\begin{overpic}[width=13cm]{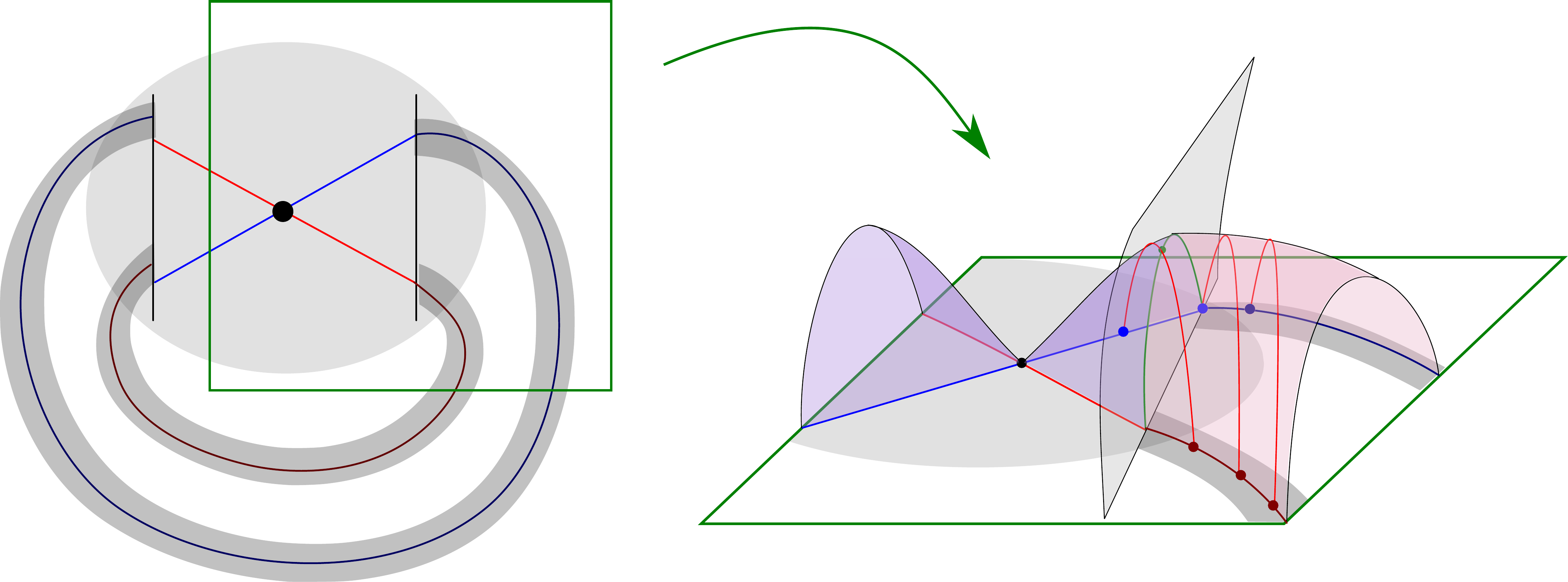}							
		\put(14,30){{\footnotesize $V$}}		
		\put(32,2){{\footnotesize $W_{\mathcal{R}}$}}		
		\put(20,18){{\footnotesize $W^{u}_{\phi_0}$}}
		\put(20,29){{\footnotesize $W^{s}_{\phi_0}$}}		
		\put(18,26){{\scriptsize $p_0$}}	
		\put(7,22){{\footnotesize $\tau^s$}}			
		\put(27,22){{\footnotesize $\tau^u$}}		
		\put(30,15){{\footnotesize $\mathcal{R}^u$}}		
		\put(34,25){{\footnotesize $\mathcal{R}^s$}}		
		\put(80,30){{\footnotesize $\tau^u$}}		
		\put(75,23){{\footnotesize $\mathcal{C}^u$}}		
		\put(0,30){{\footnotesize $\s$}}		
		\put(98,15){{\footnotesize $\s$}}			
	\end{overpic}
	\bigskip
	\caption{Neighborhood $W$ for which the extended first return map $\Phi_{0}$ is defined and behavior of the orbits of $X_0$ at points of $\mathcal{R}^u$.}
	\label{lemmafig}	
\end{figure}

Since $Y_0$ is transverse to $\tau^u$ and $\tau^s$ at $\mathcal{C}^u\cap \overline{M^-}$ and $\widehat{\mathcal{C}^u}\cap \overline{M^-}$, respectively, it follows from $(TC_3)$ condition that, for each point $p\in\mathcal{R}^{u}$, there exists either $q\in \mathcal{R}^s$ or $q\in V$ such that $p$ and $q$ are connected by a unique orbit of $Y_0$ contained in $\overline{M^-}$. Also, if $q\in V$, then such orbit intersects $\mathcal{C}^u$ or $\widehat{\mathcal{C}^u}$. See Figure \ref{lemmafig}.

Therefore, for each $p\in\mathcal{R}^u$, we use the Implicit Function Theorem to define a $\Cr$ map $\Phi^{Y_0}_p: W_p\rightarrow \s$ induced by the orbits of $Y_0$ in such a way that, $x\in W_p$ and $\Phi^{Y_0}_p(x)$ are connected by a unique orbit of $Y_0$ contained in $\overline{M^-}$. An analogous argument shows that the same holds for points of $\mathcal{R}^s$.

Clearly, given $p_1,p_2\in \mathcal{R}^s\cup\mathcal{R}^s$, if $x\in W_{p_1}\cap W_{p_2}$, then $\Phi^{Y_0}_{p_1}(x)=\Phi^{Y_0}_{p_2}(x)$.

From compactness of $\mathcal{R}^u\cup\mathcal{R}^s$, there exist a small neighborhood $W_{\mathcal{R}}\subset \s^c$ of $\mathcal{R}^u\cup \mathcal{R}^s$ and a $\Cr$ map $\Phi_{Y_0}: V\cup W_{\mathcal{R}}\rightarrow V\cup W_{\mathcal{R}}$ induced by orbits of $Y_0$. From construction, we have that $\Phi_{Y_0}$ is an involution and $\Phi_{Y_0}|_{V}=\phi_{Y_0}$. Take $W=V\cup W_{\mathcal{R}}$.

Clearly, the local invariant manifolds $W^{u,s}_{\phi_0}$ are extended through $\mathcal{R}^{u,s}$ to invariant manifolds of $\Phi_0$, respectively.

	\section*{Acknowledgements}
	OMLG is partially supported by the Brazilian CNPq grant 438975/2018-9 and the Brazilian FAPESP grant 2019/01682-4. MAT is partially supported  by the Brazilian CNPq grant 301275/2017-3 and the Brazilian FAPESP grant 2018/ 13481-0.

	  \bibliographystyle{abbrv}
	  \bibliography{references_otavio.bib}
\end{document}